\documentclass{amsart}

\newtheorem{theorem}{Theorem}[section]

\newtheorem{lemma}[theorem]{Lemma}
\newtheorem{proposition}[theorem]{Proposition}

\begin{document}

\title[Mail your comments to wonniepark@postech.ac.kr]{Analytic continuation of a biholomorphic mapping}
\author[Won K. Park]{Won K. Park\thanks{%
E-mail: wonkpark@euclid.postech.ac.kr 
\newline
} \\
.............................\\
Any comments, suggestions, errors to\\
wonniepark@postech.ac.kr}
\address{Department of Mathematics, Postech Pohang, Korea, 790-784}
\maketitle
\tableofcontents

\begin{abstract}
We present a new proof and its generalization of Pinchuk's theorem of the
analytic continuation of a biholomorphic mapping from a strictly
pseudoconvex real-analytic hypersurface to a compact strictly pseudoconvex
real-analytic hypersurface.
\end{abstract}

\addtocounter{section}{-1}

\section{\textbf{Introduction and Preliminary}}

The main purpose of this article is to present new proofs and their
generalizations of the following three theorems concerning the analytic
continuation of a biholomorphic mapping on a strongly pseudoconvex analytic
real hypersurface.

First, we concern Vitushkin's theorem of a germ of a biholomorphic
mapping(cf. \cite{Vi85}) as follows:

\begin{theorem}[Vitushkin]
Let $M,$ $M^{\prime }$ be two strongly pseudoconvex analytic real
hypersurfaces in $\Bbb{C}^{n+1}$ and $p,p^{\prime }$ be points respectively
on $M,M^{\prime }$ such that the germ $M$ at the point $p$ and the germ $%
M^{\prime }$ at the point $p^{\prime }$ are biholomorphically equivalent.
Then there is a positive real number $\delta $ depending only on $%
M,M^{\prime }$ and $p,p^{\prime }$ such that a biholomorphic mapping $\phi $
on an open connected neighborhood $U$ of the point $p$ is analytically
continued to the open ball $B(p;\delta )$ as a biholomorphic mapping
whenever 
\[
\phi (p)=p^{\prime }\quad \text{and}\quad \phi (U\cap M)\subset M^{\prime }.
\]
\end{theorem}

Next, we concern Pinchuk's theorem of the analytic continuation of a
biholomorphic mapping on a spherical analytic real hypersurface(cf. \cite
{Pi78}, \cite{CJ95}) as follows:

\begin{theorem}[Pinchuk, Chern-Ji]
Let $M$ be a spherical strongly pseudoconvex analytic real hypersurface in $%
\Bbb{C}^{n+1}$, $U$ be an open connected neighborhood of a point $p\in M,$
and $\phi $ be a biholomorphic mapping on $U$ such that $\phi \left( U\cap
M\right) \subset S^{2n+1}.$ Then the mapping $\phi $ is analytically
continued along any path in $M$ as a local biholomorphic mapping.
\end{theorem}

Finally, we concern Pinchuk's theorem of the analytic continuation of a
biholomorphic mapping on a nonspherical nondegenerate analytic real
hypersurface(cf. \cite{Pi78}, \cite{Vi85}) as follows:

\begin{theorem}[Pinchuk, Ezhov-Kruzhilin-Vitushkin]
Let $M,M^{\prime }$ be nonspherical strongly pseudoconvex analytic real
hypersurfaces in $\Bbb{C}^{n+1}$ such that $M^{\prime }$ is compact. Suppose
that there are an open connected neighborhood $U$ of a point $p\in M$ and a
biholomorphic mapping $\phi $ on $U$ such that 
\[
\phi \left( M\cap U\right) \subset M^{\prime }.
\]
Then the mapping $\phi $ is analytically continued along any path on $M$ as
a local biholomorphic mapping.
\end{theorem}

In the following subsections, we provide some preliminary results from the
papers \cite{Pa1} and \cite{Pa2}. We have attempted to present the results
of this paper in the 18th DaeWoo Workshop at Hanseo University, Korea. The
short outline of the main results in this article shall appear in the
proceedings of the Daewoo Workshop.

\subsection{Existence and uniqueness theorem}

We take a coordinate systime of $\Bbb{C}^{n}\times \Bbb{C}$ as follows: 
\[
z\equiv \left( z^{1},\cdots ,z^{n}\right) ,\quad w=u+iv\equiv z^{n+1}.
\]
A holomorphic mapping $\phi $ in $\Bbb{C}^{n}\times \Bbb{C}$ consists of $%
(n+1)$ holomorphic functions 
\[
f\equiv (f^{1},\cdots ,f^{n}),\quad g\equiv f^{n+1}.
\]
We keep the notations 
\[
\langle z,z\rangle \equiv z^{1}\overline{z}^{1}+\cdots +z^{e}\overline{z}%
^{e}-\cdots -z^{n}\overline{z}^{n}
\]
and 
\[
\Delta \equiv \frac{\partial ^{2}}{\partial z^{1}\partial \overline{z}^{1}}%
+\cdots +\frac{\partial ^{2}}{\partial z^{e}\partial \overline{z}^{e}}%
-\cdots -\frac{\partial ^{2}}{\partial z^{n}\partial \overline{z}^{n}}.
\]
Then it is known that a nondegenerate analytic hypersurface $M$ is locally
biholomorphic to a real hypersurface of the following form(cf. \cite{CM74}, 
\cite{Pa1}): 
\[
v=\langle z,z\rangle +\sum_{s,t\geq 2}^{\infty }F_{st}\left( z,\overline{z}%
,u\right) 
\]
where 
\[
\Delta F_{22}=\Delta ^{2}F_{23}=\Delta ^{3}F_{33}=0.
\]
We shall denote by $H$ the isotropy subgroup of a real hyperquadric $%
v=\langle z,z\rangle $ such that 
\[
H=\left\{ \left( 
\begin{array}{ccc}
\rho  & 0 & 0 \\ 
-\sqrt{\left| \rho \right| }Ua & \sqrt{\left| \rho \right| }U & 0 \\ 
-r-i\langle a,a\rangle  & 2ia^{\dagger } & 1
\end{array}
\right) :
\begin{array}{l}
\langle Uz,Uz\rangle =\text{\textrm{sign}}(\rho )\langle z,z\rangle ,\quad
a\in \Bbb{C}^{n}, \\ 
\rho \neq 0,\quad \rho ,r\in \Bbb{R} \\ 
a^{\dagger }=\left( \overline{a^{1}},\cdots ,\overline{a^{e}},-\overline{%
a^{e+1}},\cdots ,-\overline{a^{n}}\right) 
\end{array}
\right\} .
\]

\begin{theorem}[Chern-Moser]
\label{exi-uni}Let $M$ be a nondegenerate analytic real hypersurface in $%
\Bbb{C}^{n+1}$ defined near the origin by the equation 
\begin{equation}
v=\langle z,z\rangle +F\left( z,\overline{z},u\right)   \tag*{(0.1)}
\label{def-eq}
\end{equation}
where 
\[
F\left( z,\overline{z},u\right) =o\left( \left| z_{1}\right| ^{2}+\cdots
+\left| z_{n}\right| ^{2}+\left| w\right| ^{2}\right) .
\]
Then, for each element $(U,a,\rho ,r)\in H,$ there exists a unique
biholomorphic mapping $\phi =(f,g)$ near the origin which transforms $M$ to
a real hypersurface in Chern-Moser normal form such that 
\begin{eqnarray*}
\left( \left. \frac{\partial f}{\partial z}\right| _{0}\right)  &=&C,\quad
\left( \left. \frac{\partial f}{\partial w}\right| _{0}\right) =-Ca \\
\rm{Re}\left( \left. \frac{\partial g}{\partial w}\right| _{0}\right) 
&=&\rho ,\quad \rm{Re}\left( \left. \frac{\partial ^{2}g}{\partial w^{2}}%
\right| _{0}\right) =2\rho r
\end{eqnarray*}
where the constants $(U,a,\rho ,r)$ shall be called the initial value of the
normalization $\phi .$
\end{theorem}

We present a brief outline of the proof of Theorem \ref{exi-uni}(cf. \cite
{CM74}, \cite{Pa1}). First of all, we show that there is a biholomorphic
mapping 
\[
\phi _{1}:\left\{ 
\begin{array}{l}
z=z^{*}+D(z^{*},w^{*}) \\ 
w=w^{*}+g(z^{*},w^{*})
\end{array}
\right. 
\]
which transforms the equation \ref{def-eq} to an equation of the form 
\[
v^{*}=F_{11}^{*}\left( z^{*},\overline{z}^{*},u^{*}\right) +\sum_{s,t\geq
2}F_{st}^{*}\left( z^{*},\overline{z}^{*},u^{*}\right) . 
\]
Then we set 
\[
p(u)\equiv D(0,u) 
\]
and we verify that the functions 
\[
D(z,w),\quad g(z,w),\quad F_{st}^{*}\left( z,\overline{z},u\right) 
\]
are uniquely determined by the function $F\left( z,\overline{z},u\right) $
and $p(u)$ whenever we require the following normalizing condition 
\[
\overline{g(0,u)}=-g(0,u). 
\]
Further, the functions 
\[
\left( \left. \frac{\partial ^{\left| I\right| }D}{\partial z^{I}}\right|
_{z=v=0}\right) ,\quad \left( \left. \frac{\partial ^{\left| I\right| }g}{%
\partial z^{I}}\right| _{z=v=0}\right) ,\quad \left( \left. \frac{\partial
^{\left| I\right| +\left| J\right| }F_{st}^{*}}{\partial z^{I}\partial 
\overline{z}^{J}}\right| _{z=\overline{z}=0}\right) 
\]
depend analytically on $u$ and $p(u),$ rationally on $p^{\prime }(u),$ and
polynomially on the higher order derivatives of $p(u)$.

At this point, we need an operator $\mathrm{tr}$ introduced by Chern and
Moser as follows: 
\[
\mathrm{tr}F_{st}^{*}\left( z,\overline{z},u\right) =\frac{1}{st}%
\sum_{\alpha ,\beta =0}^{n}h^{\alpha \beta }(u)\left( \frac{\partial
^{2}F_{st}^{*}}{\partial z^{\alpha }\partial \overline{z}^{\beta }}\right)
\left( z,\overline{z},u\right) 
\]
where 
\[
F_{11}^{*}\left( z,\overline{z},u\right) =\sum_{\alpha ,\beta
=0}^{n}h_{\alpha \beta }(u)z^{\alpha }\overline{z}^{\beta } 
\]
and $\left( h^{\alpha \beta }(u)\right) $ is the inverse matrix of $\left(
h_{\alpha \beta }(u)\right) .$ Then we show that the equation 
\[
(\mathrm{tr})^{2}F_{23}^{*}\left( z,\overline{z},u\right) =0 
\]
is an ordinary differential equation of the function $p(u)$ as follows: 
\[
p^{\prime \prime }=Q(u,p,\overline{p},p^{\prime },\overline{p}^{\prime }). 
\]
Hence, for a given value $p^{\prime }(0)\equiv D_{w}(0,0)\in \Bbb{C}^{n},$
there is a unique biholomorphic mapping $\phi _{1}$ which satisfies the
normalizaing condition and which transforms the equation \ref{def-eq} to an
equation of the following form 
\begin{equation}
v=F_{11}^{*}\left( z,\overline{z},u\right) +\sum_{s,t\geq 2}F_{st}^{*}\left(
z,\overline{z},u\right)  \tag*{(0.2)}  \label{def-eq2}
\end{equation}
where 
\[
(\mathrm{tr})^{2}F_{23}^{*}\left( z,\overline{z},u\right) =0. 
\]

Note that, for a biholomorphic mapping $\phi ,$ $\phi (0)=0,$ near the
origin, there is a unique decomposition 
\[
\phi =\phi _{2}\circ \phi _{1} 
\]
where 
\begin{eqnarray*}
\phi _{1} &:&\left\{ 
\begin{array}{l}
z=z^{*}+D(z^{*},w^{*}) \\ 
w=w^{*}+g(z^{*},w^{*})
\end{array}
\right. \\
\phi _{2} &:&\left\{ 
\begin{array}{l}
z^{*}=\sqrt{\mathrm{sign\{}q^{\prime }(0)\mathrm{\}}q^{\prime }(w)}E(w)z \\ 
w^{*}=q(w)
\end{array}
\right. .
\end{eqnarray*}
and where the function $D(z,w),g(z,w),E(w),q(w)$ are complex analytic such
that 
\begin{gather*}
D(0,0)=D_{z}(0,w)=0,\quad g(0,0)=q(0)=0 \\
\overline{g(0,u)}=-g(0,u),\quad \overline{q(u)}=q(u) \\
\det E(0)\neq 0\quad \det q^{\prime }(0)\neq 0.
\end{gather*}
Let $\phi $ be the biholomorphic mapping which transforms the equation \ref
{def-eq2} to a defining equation satisfying the following condition 
\[
v=\langle z,z\rangle +\sum_{s,t\geq 2}G_{st}\left( z,\overline{z},u\right) 
\]
where 
\[
\Delta ^{2}G_{23}\left( z,\overline{z},u\right) =0. 
\]
Then there is a unique decomposition of a biholomorphic mapping $\phi $ such
that 
\[
\phi =\phi ^{*}\circ \phi _{1} 
\]
where 
\[
\phi ^{*}:\left\{ 
\begin{array}{l}
z^{*}=\sqrt{\mathrm{sign\{}q^{\prime }(0)\mathrm{\}}q^{\prime }(w)}E(w)z \\ 
w^{*}=q(w)
\end{array}
\right. 
\]
and 
\[
F_{11}^{*}\left( z,\overline{z},u\right) =\mathrm{sign\{}q^{\prime }(0)%
\mathrm{\}}\langle E(u)z,E(u)z\rangle . 
\]

Second, we take a matrix valued function $E_{1}(u)$ such that 
\[
F_{11}^{*}\left( z,\overline{z},u\right) =\langle E_{1}(u)z,E_{1}(u)z\rangle
. 
\]
Then there is a biholomorphic mapping 
\[
\phi _{2}:\left\{ 
\begin{array}{l}
z^{*}=E(w)z \\ 
w^{*}=w
\end{array}
\right. 
\]
which transforms the equation \ref{def-eq2} to an equation of the same form 
\[
v^{*}=\langle z^{*},z^{*}\rangle +\sum_{s,t\geq 2}G_{st}\left( z^{*},%
\overline{z}^{*},u^{*}\right) 
\]
where 
\[
\Delta ^{2}G_{23}\left( z,\overline{z},u\right) =0. 
\]
Further, the function $E(u)$ is uniquely determined up to a function $U(u)$
such that 
\[
E(u)=U(u)E_{1}(u) 
\]
where 
\[
\langle U(u)z,U(u)z\rangle =\langle z,z\rangle . 
\]
Then we show that the equation 
\[
\Delta G_{22}\left( z,\overline{z},u\right) =0 
\]
is an ordinary differential equation of the function $U(u)$ as follows: 
\[
U(u)^{-1}U^{\prime }(u)=R(u). 
\]
Hence, for a given value $U(0),$ there is a unique biholomorphic mapping $%
\phi _{2}$ which transforms the equation \ref{def-eq2} to an equation of the
following form 
\begin{equation}
v=\langle z,z\rangle +\sum_{s,t\geq 2}G_{st}\left( z,\overline{z},u\right) 
\tag*{(0.3)}  \label{def-eq3}
\end{equation}
where 
\[
\Delta G_{22}\left( z,\overline{z},u\right) =\Delta ^{2}G_{23}\left( z,%
\overline{z},u\right) =0. 
\]

Third, we show that there is a biholomorphic mapping 
\[
\phi _{3}:\left\{ 
\begin{array}{l}
z^{*}=\sqrt{\mathrm{sign\{}q^{\prime }(0)\mathrm{\}}q^{\prime }(w)}z \\ 
w^{*}=q(w)
\end{array}
\right. 
\]
which transforms the equation \ref{def-eq3} to an equation of the same form 
\[
v^{*}=\langle z^{*},z^{*}\rangle +\sum_{s,t\geq 2}G_{st}^{*}\left( z^{*},%
\overline{z}^{*},u^{*}\right) 
\]
where 
\[
\Delta G_{22}^{*}\left( z,\overline{z},u\right) =\Delta ^{2}G_{23}^{*}\left(
z,\overline{z},u\right) =0.
\]
Then we show that the equation 
\[
\Delta ^{3}G_{33}^{*}\left( z,\overline{z},u\right) =0
\]
is an ordinary differential equation of the function $q(u)$ as follows: 
\[
\frac{q^{\prime \prime \prime }}{3q^{\prime }}-\frac{1}{2}\left( \frac{%
q^{\prime \prime }}{q^{\prime }}\right) ^{2}=\kappa (u).
\]
Hence, for given values $q^{\prime }(0),q^{\prime \prime }(0),$ there is a
unique biholomorphic mapping $\phi _{3}$ which transforms the equation \ref
{def-eq3} to an equation of the following form 
\[
v=\langle z,z\rangle +\sum_{s,t\geq 2}G_{st}^{*}\left( z,\overline{z}%
,u\right) 
\]
where 
\[
\Delta G_{22}\left( z,\overline{z},u\right) =\Delta ^{2}G_{23}\left( z,%
\overline{z},u\right) =\Delta ^{3}G_{33}\left( z,\overline{z},u\right) =0.
\]
Thus the existence and uniqueness of the biholomorphic mapping $\phi $ have
been reduced to the existence and uniqueness of solutions of the ordinary
differential equations, where some constants $U,a,\rho ,r$ appear as the
initial values of the solutions.

In the paper \cite{Pa1}, we have showed that there exist a family of normal
forms as follows: 
\[
v=\frac{1}{4\alpha }\ln \frac{1}{1-4\alpha \langle z,z\rangle }%
+\sum_{s,t\geq 2}F_{st}\left( z,\overline{z},u\right) 
\]
where $\alpha ,\beta \in \Bbb{R}$ and 
\[
\left\{ 
\begin{array}{l}
\Delta F_{22}=\Delta ^{2}F_{23}=0 \\ 
\Delta ^{3}F_{33}=\beta \Delta ^{4}\left( F_{22}\right) ^{2}.
\end{array}
\right. 
\]
In the case of $\alpha =0,$ we assume 
\[
v=\langle z,z\rangle +\sum_{s,t\geq 2}F_{st}\left( z,\overline{z},u\right) . 
\]
The value $(\alpha ,\beta )$ is called the type of normal form. Chern-Moser
normal form is given in the case of $\alpha =\beta =0$ and Moser-Vitushkin
normal form is defined by taking $\alpha \neq 0$ and $\beta =0.$

Then each normalization of a real hypersurface $M$ to a normal form of a
given type $(\alpha ,\beta )$ is determined by constant initial value
parameterized by the local automorphism group $H$ of the following real
hypersurface 
\[
v=\frac{1}{4\alpha }\ln \frac{1}{1-4\alpha \langle z,z\rangle }, 
\]
which is locally biholomorphic to a real hyperquadric.

\begin{theorem}
\label{exi-uni2}Let $M$ be a nondegenerate analytic real hypersurface
defined by 
\[
v=\sum_{k=2}^{\infty }F_{k}\left( z,\overline{z},u\right) .
\]
Then there exist unique natural mappings for each $k\geq 2$ such that 
\begin{eqnarray*}
\nu :\left\{ F_{l}\left( z,\overline{z},u\right) :l\leq k\right\} \times
H\times \Bbb{R}^{2}\longmapsto \left( f_{k-1}\left( z,w\right) ,g_{k}\left(
z,w\right) \right)  \\
\kappa :\left\{ F_{l}\left( z,\overline{z},u\right) :l\leq k\right\} \times
H\times \Bbb{R}^{2}\longmapsto F_{k}^{*}\left( z,\overline{z},u\right) 
\end{eqnarray*}
such that, for a given $\sigma \in H$ and $\alpha ,\beta \in \Bbb{R},$ the
formal series mapping 
\[
\phi =\left( \sum_{k=1}^{\infty }f_{k}\left( z,w\right) ,\sum_{k=2}^{\infty
}g_{k}\left( z,w\right) \right) 
\]
is a biholomorphic normalization of $M$ with initial value $\sigma \in H$
and 
\[
v=\langle z,z\rangle +\sum_{k=4}^{\infty }F_{k}^{*}\left( z,\overline{z}%
,u\right) 
\]
is the defining equation of the real hypersurface $\phi \left( M\right) $ in
normal form of type $(\alpha ,\beta ).$
\end{theorem}

Then we obtain the following theorem

\begin{theorem}
\label{main}Let $M$ be a nondegenerate analytic real hypersurface defined by 
\[
v=\sum_{k=2}^{\infty }F_{2}(z,\overline{z},u)
\]
and $\phi =(\sum_{k}f_{k},\sum_{k}g_{k})$ be a normalization of $M$ such
that the real hypersurface $\phi \left( M\right) $ is defined in normal form
of type $(\alpha ,\beta )$ by the equation 
\[
v=\langle z,z\rangle +\sum_{k=4}^{\infty }F_{k}^{*}(z,\overline{z},u).
\]
Then the functions $f_{k-1},g_{k},F_{k}^{*},$ $k\geq 3,$ are given as a
finite linear combination of finite multiples of the following factors:

\begin{enumerate}
\item[(1)]  the coefficients of the functions $F_{l}$, $l\leq k,$

\item[(2)]  the constants $C,C^{-1},\rho ,\rho ^{-1},a,r,\alpha ,\beta ,$
\end{enumerate}

\noindent where $(C,a,\rho ,r)$ are the initial value of the normalization $%
\phi $ and $\alpha ,\beta $ are the parameters of normal forms.
\end{theorem}

\subsection{Equation of a chain}

In the proof of Theorem \ref{exi-uni}, we have a distinguished curve $\gamma 
$ on $M,$ which is named a chain by E. Cartan \cite{Ca32} and Chern-Moser 
\cite{CM74}. Suppose that there is a nondegenerate analytic real
hypersurface $M$ defined near the origin by 
\[
v=F(z,\overline{z},u),\quad \left. F\right| _{0}=\left. F_{z}\right|
_{0}=\left. F_{\overline{z}}\right| _{0}=0. 
\]
Then there exists an ordinary differential equation 
\begin{equation}
p^{\prime \prime }=Q\left( u,p,\overline{p},p^{\prime },\overline{p}^{\prime
}\right)  \tag*{(0.4)}  \label{ordinary}
\end{equation}
such that a chain $\gamma ,$ passing through the origin $0\in M$, is given
near the origin by the equation 
\[
\gamma :\left\{ 
\begin{array}{l}
z=p(u) \\ 
w=u+iF\left( p(u),\overline{p}(u),u\right)
\end{array}
\right. 
\]
where $p(u)$ is a solution of the ordinary differential equation \ref
{ordinary}.

The explicit form of the equation \ref{ordinary}, which depends on the
function $F\left( z,\overline{z},u\right) ,$ is quite complicate(cf. \cite
{CM74}, \cite{Pa1}). Roughly, the function $Q$ in \ref{ordinary} is given as
follows: 
\begin{equation}
Q(u,p,\overline{p},p^{\prime },\overline{p}^{\prime })=(A_{1}-A_{2}\overline{%
A_{1}}^{-1}\overline{A_{2}})^{-1}(B-A_{2}\overline{A_{1}}^{-1}\overline{B}) 
\tag*{(0.5)}  \label{Q-function}
\end{equation}
where

\begin{enumerate}
\item[(1)]  $A_{1},A_{2},B$ are functions of $u,p,\overline{p},p^{\prime },%
\overline{p}^{\prime }$,

\item[(2)]  $A_{1},A_{2}$ are $n\times n$ matrices respectively given by 
\begin{eqnarray*}
\left[ A_{1}\left( u,p,\overline{p},p^{\prime },\overline{p}^{\prime
}\right) \right] _{\alpha \beta } \\
=\left\{ 2iF_{\beta \overline{\alpha }}+2\left( 1+iF^{\prime }\right)
^{-1}\left( F_{\beta }^{\prime }+iF^{\prime \prime }F_{\beta }\right) F_{%
\overline{\alpha }}\right\} \times  \\
\qquad \qquad \left\{ 1-i\left( 1+iF^{\prime }\right) F_{\gamma }p^{\gamma
\prime }+i\left( 1-iF^{\prime }\right) F_{\overline{\gamma }}p^{\overline{%
\gamma }\prime }+F^{\prime 2}\right\}  \\
-i\left( 1+iF^{\prime }\right) F_{\beta }\times \left\{ 2iF_{\gamma 
\overline{\alpha }}p^{\gamma \prime }+2F^{\prime \prime }F_{\overline{\alpha 
}}+i\left( 1+iF^{\prime }\right) F_{\overline{\alpha }}^{\prime }\right.  \\
\qquad \qquad \left. +2\left( 1+iF^{\prime }\right) ^{-1}F_{\overline{\alpha 
}}\left( F_{\gamma }^{\prime }p^{\gamma \prime }+iF^{\prime \prime
}F_{\gamma }p^{\gamma \prime }+iF^{\prime \prime }F_{\overline{\gamma }}p^{%
\overline{\gamma }\prime }\right) \right\} 
\end{eqnarray*}
and 
\begin{eqnarray*}
\left[ A_{2}\left( u,p,\overline{p},p^{\prime },\overline{p}^{\prime
}\right) \right] _{\alpha \beta } \\
=2iF^{\prime \prime }\left( 1+iF^{\prime }\right) ^{-1}F_{\overline{\alpha }%
}F_{\overline{\beta }}\times  \\
\qquad \qquad \left\{ 1-i\left( 1+iF^{\prime }\right) F_{\gamma }p^{\gamma
\prime }+i\left( 1-iF^{\prime }\right) F_{\overline{\gamma }}p^{\overline{%
\gamma }\prime }+F^{\prime 2}\right\}  \\
+i\left( 1-iF^{\prime }\right) F_{\overline{\beta }}\times \left\{
2iF_{\gamma \overline{\alpha }}p^{\gamma \prime }+2F^{\prime \prime }F_{%
\overline{\alpha }}+i\left( 1+iF^{\prime }\right) F_{\overline{\alpha }%
}^{\prime }\right.  \\
\qquad \qquad \left. +2\left( 1+iF^{\prime }\right) ^{-1}F_{\overline{\alpha 
}}\left( F_{\gamma }^{\prime }p^{\gamma \prime }+iF^{\prime \prime
}F_{\gamma }p^{\gamma \prime }+iF^{\prime \prime }F_{\overline{\gamma }}p^{%
\overline{\gamma }\prime }\right) \right\} 
\end{eqnarray*}
where 
\begin{eqnarray*}
F_{\alpha }=\left( \frac{\partial F}{\partial z^{\alpha }}\right) \left(
p(u),\overline{p}(u),u\right) ,\quad F_{\overline{\beta }}=\left( \frac{%
\partial F}{\partial \overline{z}^{\beta }}\right) \left( p(u),\overline{p}%
(u),u\right)  \\
F^{\prime }=\left( \frac{\partial F}{\partial u}\right) \left( p(u),%
\overline{p}(u),u\right) ,\quad F^{\prime \prime }=\frac{1}{2}\left( \frac{%
\partial ^{2}F}{\partial u^{2}}\right) \left( p(u),\overline{p}(u),u\right) 
\\
F_{\alpha }^{\prime }=\left( \frac{\partial ^{2}F}{\partial z^{\alpha
}\partial u}\right) \left( p(u),\overline{p}(u),u\right) ,\quad F_{\alpha 
\overline{\beta }}=\left( \frac{\partial ^{2}F}{\partial z^{\alpha }\partial 
\overline{z}^{\beta }}\right) \left( p(u),\overline{p}(u),u\right) ,
\end{eqnarray*}

\item[(3)]  $B$ is a $n\times 1$ matrix given by at most cubic polynomial
with respect to $p^{\prime },\overline{p}^{\prime }$ such that $B$ is a
finite linear combination of multiples of the derivatives $p^{\prime },%
\overline{p}^{\prime }$ and the following terms: 
\[
\left( \frac{\partial ^{\left| I\right| +\left| J\right| +m}F}{\partial
z^{I}\partial \overline{z}^{J}\partial u^{m}}\right) \left( p(u),\overline{p}%
(u),u\right) \quad \text{for }\left| I\right| +\left| J\right| +m\leq 5
\]
and 
\[
\left[ \det \left\{ \left( 1-iF^{\prime }\right) ^{2}F_{\alpha \overline{%
\beta }}-i\left( 1+iF^{\prime }\right) F_{\alpha }^{\prime }F_{\overline{%
\beta }}+i\left( 1-iF^{\prime }\right) F_{\overline{\beta }}^{\prime
}F_{\alpha }+2F^{\prime \prime }F_{\alpha }F_{\overline{\beta }}\right\}
\right] ^{-1}.
\]
\end{enumerate}

On the real hyperquadric $v=\langle z,z\rangle ,$ the chain $\gamma $ is
locally given by 
\[
\gamma :\left\{ 
\begin{array}{l}
z=p(u) \\ 
w=u+i\langle p(u),p(u)\rangle
\end{array}
\right. 
\]
where $p(u)$ is a solution of the ordinary differential equation(cf. \cite
{Pa1}): 
\[
p^{\prime \prime }=\frac{2ip^{\prime }\langle p^{\prime },p^{\prime }\rangle
\left( 1+3i\langle p,p^{\prime }\rangle -i\langle p^{\prime },p\rangle
\right) }{\left( 1+i\langle p,p^{\prime }\rangle -i\langle p^{\prime
},p\rangle \right) \left( 1+2i\langle p,p^{\prime }\rangle -2i\langle
p^{\prime },p\rangle \right) }. 
\]
Further, the chain $\gamma $ on a real hyperquadric $v=\langle z,z\rangle $
is necessarily given as an intersection of a complex line(cf. \cite{CM74}, 
\cite{Pa1}).

Then we may define a chain $\gamma $ globally. Let $M$ be a nondegenerate
analytic real hypersurface and $\gamma :(0,1)\rightarrow M$ be an open
connected curve. Then the curve $\gamma $ is called a chain if, for each
point $p\in \gamma ,$ there exist an open neighborhood $U$ of the point $p$
and a biholomorphic mapping $\phi $ which translates the point $p$ to the
origin and transforms $M$ to Chern-Moser normal form such that 
\[
\phi \left( U\cap \gamma \right) \subset \left\{ z=v=0\right\} . 
\]
An alternative definition of a chain $\gamma $ may be given through the
intrinsic geometry of nondegenerate real hypersurfaces(cf. \cite{Ca32}, \cite
{CM74}, \cite{Ta76}).

\section{Nonsingular matrices}

\subsection{A family of nonsingular matrices}

\begin{lemma}
\label{L2}Let $A_{m}$ be a matrix as follows: 
\[
\left( 
\begin{array}{cccccc}
0 & 2 & 0 & \cdots  &  & 0 \\ 
m & 3 & 4 & \ddots  &  &  \\ 
0 & m-1 & 6 & 6 & \ddots  & \vdots  \\ 
\vdots  & \ddots  & \ddots  & \ddots  & \ddots  & 0 \\ 
&  & \ddots  & 2 & 3m-3 & 2m \\ 
0 &  & \cdots  & 0 & 1 & 3m
\end{array}
\right) .
\]
Then the eigenvalues of $A_{m}$ are given by 
\[
\frac{3m}{2}+\frac{m-2s}{2}\sqrt{17}\quad \text{for }s=0,\cdots ,m.
\]
\end{lemma}

\proof
We consider the following system of first order ordinary differential
equations: 
\begin{eqnarray*}
y^{\prime } &=&z \\
z^{\prime } &=&3z+2y.
\end{eqnarray*}
Then the general solutions $y,z$ are given by 
\begin{align}
y(t)& =c_{1}e^{t\lambda _{1}}+c_{2}e^{t\lambda _{2}},  \nonumber \\
z(t)& =c_{1}\lambda _{1}e^{t\lambda _{1}}+c_{2}\lambda _{2}e^{t\lambda _{2}},
\tag*{(1.1)}  \label{sys}
\end{align}
where $\lambda _{1},\lambda _{2},\lambda _{1}\neq \lambda _{2},$ are the two
solutions of the quadratic equation: 
\[
x^{2}-3x-2=0, 
\]
and $c_{1},c_{2}$ are arbitrary real numbers.

We take nonzero constants $c_{1},c_{2}$ so that $y(t),z(t)$ are linear
independent. Then we obtain 
\begin{equation}
\left( 
\begin{array}{l}
e^{t\lambda _{1}} \\ 
e^{t\lambda _{2}}
\end{array}
\right) =\left( 
\begin{array}{cc}
c_{1} & c_{2} \\ 
c_{1}\lambda _{1} & c_{2}\lambda _{2}
\end{array}
\right) ^{-1}\left( 
\begin{array}{l}
y(t) \\ 
z(t)
\end{array}
\right) .  \tag*{(1.2)}  \label{inv}
\end{equation}
We consider a real vector space $V$ generated by the following elements: 
\[
y^{m-s}z^{s}\quad \text{for}\quad s=0,1,\cdots ,m. 
\]
By the equalities \ref{sys} and \ref{inv}, the vector space $V$ is generated
as well by the following elements: 
\begin{equation}
\exp t(s\lambda _{1}+(m-s)\lambda _{2})\quad \text{for}\quad s=0,1,\cdots ,m.
\tag*{(1.3)}  \label{eign}
\end{equation}

We put 
\[
B_{1}=\left( 
\begin{array}{l}
y^{m} \\ 
y^{m-1}z \\ 
\vdots \\ 
yz^{m-1} \\ 
z^{m}
\end{array}
\right) \quad \text{and}\quad B_{2}=\left( 
\begin{array}{l}
e^{tm\lambda _{1}} \\ 
e^{t((m-1)\lambda _{1}+\lambda _{2})} \\ 
\vdots \\ 
e^{t(\lambda _{1}+(m-1)\lambda _{2})} \\ 
e^{tm\lambda _{2}}
\end{array}
\right) . 
\]
Then it is verified that 
\begin{equation}
\frac{dB_{1}}{dt}=\left( 
\begin{array}{cccccc}
0 & 2 & 0 & \cdots &  & 0 \\ 
m & 3 & 4 & \ddots &  &  \\ 
0 & m-1 & 6 & 6 & \ddots & \vdots \\ 
\vdots & \ddots & \ddots & \ddots & \ddots & 0 \\ 
&  & \ddots & 2 & 3m-3 & 2m \\ 
0 &  & \cdots & 0 & 1 & 3m
\end{array}
\right) ^{T}B_{1},  \tag*{(1.4)}  \label{dt}
\end{equation}
and 
\[
\frac{dB_{2}}{dt}=\left( 
\begin{array}{ccccc}
m\lambda _{1} & 0 & \cdots &  & 0 \\ 
0 & (m-1)\lambda _{1}+\lambda _{2} & \ddots &  &  \\ 
& \ddots & \ddots & \ddots & \vdots \\ 
\vdots &  & \ddots & \lambda _{1}+(m-1)\lambda _{2} & 0 \\ 
0 & \cdots &  & 0 & m\lambda _{2}
\end{array}
\right) B_{2}. 
\]
Hence the derivative $\frac{d}{dt}$ is an endomorphism on $V$ and the
vectors in \ref{eign} are the eigenvectors of the endomorphism $\frac{d}{dt}%
. $ Thus the matrix $A_{m}$ has eigenvalues as follows: 
\[
s\lambda _{1}+(m-s)\lambda _{2}\quad \text{for}\quad s=0,1,\cdots ,m 
\]
where 
\[
\lambda _{1},\lambda _{2}=\frac{3\pm \sqrt{17}}{2}. 
\]
This completes the proof.\endproof

\begin{lemma}
\label{Corr1}Let $B_{m}$ be a matrix as follows: 
\[
B_{m}=\left( 
\begin{array}{cccccc}
1 & 2 & 3 & \cdots  & m & m+1 \\ 
m & 7-m & 4 & 0 & \cdots  & 0 \\ 
0 & m-1 & 10-m & 6 & \ddots  & \vdots  \\ 
& \ddots  & \ddots  & \ddots  & \ddots  & 0 \\ 
\vdots  &  & \ddots  & 2 & 2m+1 & 2m \\ 
0 & \cdots  &  & 0 & 1 & 2m+4
\end{array}
\right) .
\]
Then the matrix $B_{m}$ is nonsingular.
\end{lemma}

\proof
We easily verify that 
\begin{equation}
\det B_{m}=\frac{1}{4}\det C_{m}  \tag*{(1.5)}  \label{bc}
\end{equation}
where 
\[
C_{m}=\left( 
\begin{array}{cccccc}
4-m & 2 & 0 &  & \cdots & 0 \\ 
m & 7-m & 4 & \ddots &  & \vdots \\ 
0 & m-1 & 10-m & 6 & \ddots &  \\ 
& \ddots & \ddots & \ddots & \ddots & 0 \\ 
\vdots &  & \ddots & 2 & 2m+1 & 2m \\ 
0 & \cdots &  & 0 & 1 & 2m+4
\end{array}
\right) . 
\]
Note that 
\begin{equation}
C_{m}=A_{m}-(m-4)id_{(m+1)\times (m+1)}.  \tag*{(1.6)}  \label{dt2}
\end{equation}
By Lemma \ref{L2}, the eigenvalues of the matrix $A_{m}$ is given as
follows: 
\[
\frac{3m}{2}+\frac{m-2s}{2}\sqrt{17}\quad \text{for }s=0,\cdots ,m. 
\]
Thus the eigenvalues of the matrix $C_{m}$ is given by 
\[
\frac{m+8}{2}+\frac{m-2s}{2}\sqrt{17}\quad \text{for }s=0,\cdots ,m. 
\]
The matrix $C_{m}$ does not have $0$ as its eigenvalue. Therefore the matrix 
$C_{m}$ is nonsingular, i.e., 
\[
\det C_{m}\neq 0. 
\]
By the relation \ref{bc}, 
\[
\det B_{m}=\frac{1}{4}\det C_{m}\neq 0 
\]
so that the matrix $B_{m}$ is nonsingular. This completes the proof.\endproof

\subsection{Sufficient condition for Nonsingularity}

\begin{lemma}
\label{Lemm1}Let $\Delta (m),$ $m\in \Bbb{N},$ denote the function defined
as follows: 
\[
\Delta (m)=\sum_{k=0}^{m}\frac{\binom{m}{k}}{k\lambda _{1}+(m-k)\lambda
_{2}-(m-4)}\left( \frac{\lambda _{2}}{\lambda _{2}-\lambda _{1}}\right)
^{k}\left( \frac{-\lambda _{1}}{\lambda _{2}-\lambda _{1}}\right) ^{m-k}
\]
where 
\[
\lambda _{1}=\frac{3-\sqrt{17}}{2},\quad \lambda _{2}=\frac{3+\sqrt{17}}{2}.
\]
Then 
\[
\frac{\det E_{m}(m+1)}{\det E_{m}(m)}=\Delta (m)^{-1}
\]
where 
\[
E_{m}(m+1)=\left( 
\begin{array}{ccccc}
4-m & 2 & 0 & \cdots  & 0 \\ 
m & 7-m & 4 & \ddots  & \vdots  \\ 
0 & \ddots  & \ddots  & \ddots  & 0 \\ 
\vdots  & \ddots  & 2 & 2m+1 & 2m \\ 
0 & \cdots  & 0 & 1 & 2m+4
\end{array}
\right) 
\]
and 
\[
E_{m}(m)=\left( 
\begin{array}{ccccc}
7-m & 4 & 0 & \cdots  & 0 \\ 
m-1 & 10-m & 6 & \ddots  & \vdots  \\ 
0 & \ddots  & \ddots  & \ddots  & 0 \\ 
\vdots  & \ddots  & 2 & 2m+1 & 2m \\ 
0 & \cdots  & 0 & 1 & 2m+4
\end{array}
\right) .
\]
\end{lemma}

\proof
We easily see that 
\[
\varepsilon =\frac{\det E_{m}(m+1)}{\det E_{m}(m)} 
\]
if and only if the following matrix is singular: 
\[
\left( 
\begin{array}{cccccc}
4-m-\varepsilon & 2 & 0 &  & \cdots & 0 \\ 
m & 7-m & 4 & \ddots &  & \vdots \\ 
0 & m-1 & 10-m & 6 & \ddots &  \\ 
& \ddots & \ddots & \ddots & \ddots & 0 \\ 
\vdots &  & \ddots & 2 & 2m+1 & 2m \\ 
0 & \cdots &  & 0 & 1 & 2m+4
\end{array}
\right) . 
\]
Then, by the equalities \ref{dt} and \ref{dt2}, there are constants $c_{s},$ 
$s=0,\cdots ,m-1,$ which are not all zero and satisfy the following
equality: 
\begin{equation}
\sum_{s=0}^{m-1}c_{s}\frac{d}{dt}\left( y^{s}z^{m-s}e^{-t(m-4)}\right) =%
\frac{d}{dt}\left( y^{m}e^{-t(m-4)}\right) -\varepsilon y^{m}e^{-t(m-4)} 
\tag*{(1.7)}  \label{dtt}
\end{equation}
whenever 
\[
\varepsilon =\frac{\det E_{m}(m+1)}{\det E_{m}(m)}. 
\]
By the expression \ref{sys}, we obtain 
\begin{eqnarray*}
y^{m}e^{-t(m-4)} &=&(c_{1}e^{t\lambda _{1}}+c_{2}e^{t\lambda
_{2}})^{m}e^{-t(m-4)} \\
&=&\sum_{k=0}^{m}\binom{m}{k}c_{1}^{k}c_{2}^{m-k}e^{tk\lambda
_{1}+t(m-k)\lambda _{2}}e^{-t(m-4)} \\
&=&\frac{d}{dt}\left\{ \sum_{k=0}^{m}\binom{m}{k}\frac{%
c_{1}^{k}c_{2}^{m-k}e^{tk\lambda _{1}+t(m-k)\lambda _{2}}e^{-t(m-4)}}{%
k\lambda _{1}+(m-k)\lambda _{2}-(m-4)}\right\} .
\end{eqnarray*}
By using the expression \ref{inv}, we obtain 
\begin{equation}
y^{m}e^{-t(m-4)}=\frac{d}{dt}\left\{ \sum_{k=0}^{m}\frac{\binom{m}{k}\left( 
\frac{\lambda _{2}y-z}{\lambda _{2}-\lambda _{1}}\right) ^{k}\left( \frac{%
-\lambda _{1}y+z}{\lambda _{2}-\lambda _{1}}\right) ^{m-k}e^{-t(m-4)}}{%
k\lambda _{1}+(m-k)\lambda _{2}-(m-4)}\right\} .  \tag*{(1.8)}  \label{dtt2}
\end{equation}
Because the derivative $\frac{d}{dt}$ is an isomorphism, the equalities \ref
{dtt} and \ref{dtt2} yields 
\begin{equation}
\sum_{s=0}^{m-1}c_{s}y^{s}z^{m-s}=y^{m}-\varepsilon \sum_{k=0}^{m}\frac{%
\binom{m}{k}\left( \frac{\lambda _{2}y-z}{\lambda _{2}-\lambda _{1}}\right)
^{k}\left( \frac{-\lambda _{1}y+z}{\lambda _{2}-\lambda _{1}}\right) ^{m-k}}{%
k\lambda _{1}+(m-k)\lambda _{2}-(m-4)}.  \tag*{(1.9)}  \label{house}
\end{equation}
We easily see that the equality \ref{house} is satisfied by some constants $%
c_{s}$ only if we have the following equality: 
\[
1-\varepsilon \sum_{k=0}^{m}\frac{\binom{m}{k}\left( \frac{\lambda _{2}}{%
\lambda _{2}-\lambda _{1}}\right) ^{k}\left( \frac{-\lambda _{1}}{\lambda
_{2}-\lambda _{1}}\right) ^{m-k}}{k\lambda _{1}+(m-k)\lambda _{2}-(m-4)}=0. 
\]
Thus we have 
\[
\varepsilon =\Delta (m)^{-1}. 
\]
This completes the proof.\endproof

\begin{lemma}
\label{Lemm2}Let $B_{m}(2),$ $m\geq 3,$ and $B_{m}(3),$ $m\geq 4,$ be
matrices as follows: 
\[
B_{m}(2)=\left( 
\begin{array}{cccccc}
2 & 3 & 4 & \cdots  & m & m+1 \\ 
m-1 & 7-m & 6 & 0 & \cdots  & 0 \\ 
0 & m-2 & 10-m & 8 & \ddots  & \vdots  \\ 
& \ddots  & \ddots  & \ddots  & \ddots  & 0 \\ 
\vdots  &  & \ddots  & 2 & 2m+1 & 2m \\ 
0 & \cdots  &  & 0 & 1 & 2m+4
\end{array}
\right) 
\]
and 
\[
B_{m}(3)=\left( 
\begin{array}{cccccc}
3 & 4 & 5 & \cdots  & m & m+1 \\ 
m-2 & 10-m & 8 & 0 & \cdots  & 0 \\ 
0 & m-3 & 13-m & 10 & \ddots  & \vdots  \\ 
& \ddots  & \ddots  & \ddots  & \ddots  & 0 \\ 
\vdots  &  & \ddots  & 2 & 2m+1 & 2m \\ 
0 & \cdots  &  & 0 & 1 & 2m+4
\end{array}
\right) .
\]
Then 
\begin{eqnarray*}
\det B_{m}(2)\neq 0\quad \text{if and only if}\quad \Delta (m)^{-1}\neq 4, \\
\det B_{m}(3)\neq 0\quad \text{if and only if}\quad \Delta (m)^{-1}\neq -%
\frac{4}{3}(m-3).
\end{eqnarray*}
\end{lemma}

\proof
We easily verify that 
\[
\det B_{m}(2)=\frac{1}{4}\det C_{m}(2) 
\]
where 
\[
C_{m}(2)=\left( 
\begin{array}{cccccc}
9-m & 4 & 0 &  & \cdots & 0 \\ 
m-1 & 10-m & 6 & \ddots &  & \vdots \\ 
0 & m-2 & 13-m & 8 & \ddots &  \\ 
& \ddots & \ddots & \ddots & \ddots & 0 \\ 
\vdots &  & \ddots & 2 & 2m+1 & 2m \\ 
0 & \cdots &  & 0 & 1 & 2m+4
\end{array}
\right) . 
\]
Note that 
\[
\det C_{m}(2)=0 
\]
if and only if there are numbers $c_{1},\cdots ,c_{m-1}$ satisfying 
\[
\left( 
\begin{array}{cccccc}
9-m & 4 & 0 &  & \cdots & 0 \\ 
m-1 & 10-m & 6 & \ddots &  & \vdots \\ 
0 & m-2 & 13-m & 8 & \ddots &  \\ 
& \ddots & \ddots & \ddots & \ddots & 0 \\ 
\vdots &  & \ddots & 2 & 2m+1 & 2m \\ 
0 & \cdots &  & 0 & 1 & 2m+4
\end{array}
\right) \left( 
\begin{array}{c}
1 \\ 
c_{1} \\ 
\\ 
\vdots \\ 
\\ 
c_{m-1}
\end{array}
\right) =0. 
\]
Then we easily see 
\[
\left( 
\begin{array}{cccccc}
-m & 2 & 0 &  & \cdots & 0 \\ 
m & 7-m & 4 & \ddots &  & \vdots \\ 
0 & m-1 & 10-m & 6 & \ddots &  \\ 
& \ddots & \ddots & \ddots & \ddots & 0 \\ 
\vdots &  & \ddots & 2 & 2m+1 & 2m \\ 
0 & \cdots &  & 0 & 1 & 2m+4
\end{array}
\right) \left( 
\begin{array}{c}
\frac{2}{m} \\ 
1 \\ 
c_{1} \\ 
\vdots \\ 
\\ 
c_{m-1}
\end{array}
\right) =0 
\]
so that 
\[
\det \left( 
\begin{array}{cccccc}
-m & 2 & 0 &  & \cdots & 0 \\ 
m & 7-m & 4 & \ddots &  & \vdots \\ 
0 & m-1 & 10-m & 6 & \ddots &  \\ 
& \ddots & \ddots & \ddots & \ddots & 0 \\ 
\vdots &  & \ddots & 2 & 2m+1 & 2m \\ 
0 & \cdots &  & 0 & 1 & 2m+4
\end{array}
\right) =0. 
\]
Hence, by Lemma \ref{Lemm1}, we verify that 
\[
\det B_{m}(2)\neq 0 
\]
if and only if 
\[
\Delta (m)^{-1}\neq 4. 
\]

For the case of $B_{m}(3),$ we easily verify that 
\[
\det B_{m}(3)=\frac{1}{4}\det C_{m}(3) 
\]
where 
\[
C_{m}(3)=\left( 
\begin{array}{cccccc}
14-m & 6 & 0 &  & \cdots & 0 \\ 
m-2 & 13-m & 8 & \ddots &  & \vdots \\ 
0 & m-3 & 16-m & 10 & \ddots &  \\ 
& \ddots & \ddots & \ddots & \ddots & 0 \\ 
\vdots &  & \ddots & 2 & 2m+1 & 2m \\ 
0 & \cdots &  & 0 & 1 & 2m+4
\end{array}
\right) . 
\]
Note that 
\[
\det C_{m}(3)=0 
\]
if and only if there are numbers $c_{1},\cdots ,c_{m-2}$ satisfying 
\[
\left( 
\begin{array}{cccccc}
14-m & 6 & 0 &  & \cdots & 0 \\ 
m-2 & 13-m & 8 & \ddots &  & \vdots \\ 
0 & m-3 & 16-m & 10 & \ddots &  \\ 
& \ddots & \ddots & \ddots & \ddots & 0 \\ 
\vdots &  & \ddots & 2 & 2m+1 & 2m \\ 
0 & \cdots &  & 0 & 1 & 2m+4
\end{array}
\right) \left( 
\begin{array}{c}
1 \\ 
c_{1} \\ 
\\ 
\vdots \\ 
\\ 
c_{m-2}
\end{array}
\right) =0. 
\]
Then we easily see 
\[
\left( 
\begin{array}{cccccc}
\frac{m}{3} & 2 & 0 &  & \cdots & 0 \\ 
m & 7-m & 4 & \ddots &  & \vdots \\ 
0 & m-1 & 10-m & 6 & \ddots &  \\ 
& \ddots & \ddots & \ddots & \ddots & 0 \\ 
\vdots &  & \ddots & 2 & 2m+1 & 2m \\ 
0 & \cdots &  & 0 & 1 & 2m+4
\end{array}
\right) \left( 
\begin{array}{c}
-\frac{24}{m(m-1)} \\ 
\frac{4}{m-1} \\ 
1 \\ 
c_{1} \\ 
\vdots \\ 
c_{m-2}
\end{array}
\right) =0 
\]
so that 
\[
\det \left( 
\begin{array}{cccccc}
\frac{m}{3} & 2 & 0 &  & \cdots & 0 \\ 
m & 7-m & 4 & \ddots &  & \vdots \\ 
0 & m-1 & 10-m & 6 & \ddots &  \\ 
& \ddots & \ddots & \ddots & \ddots & 0 \\ 
\vdots &  & \ddots & 2 & 2m+1 & 2m \\ 
0 & \cdots &  & 0 & 1 & 2m+4
\end{array}
\right) =0. 
\]
Hence, by Lemma \ref{Lemm1}, we verify that 
\[
\det B_{m}(3)\neq 0 
\]
if and only if 
\[
\Delta (m)^{-1}\neq -\frac{4}{3}(m-3). 
\]
This completes the proof.\endproof

\subsection{Estimates}

\begin{lemma}
\label{17}For any two positive integers $p,q,$ the inequality 
\[
\left| \sqrt{17}-\frac{p}{q}\right| >\frac{2}{17q^{2}}
\]
is satisfied$.$
\end{lemma}

\proof
Note that 
\[
\left| \sqrt{17}-\frac{p}{q}\right| =\frac{\left| p^{2}-17q^{2}\right| }{%
\left| \sqrt{17}+\frac{p}{q}\right| q^{2}}\geq \frac{1}{\left| \sqrt{17}+%
\frac{p}{q}\right| q^{2}}. 
\]
Let $c$ be a positive real number. Then we consider integer pairs $(p,q)$
such that 
\[
\left| \sqrt{17}-\frac{p}{q}\right| \leq \frac{1}{c}, 
\]
which yields 
\begin{eqnarray*}
\left| \sqrt{17}+\frac{p}{q}\right| &\leq &2\sqrt{17}+\left| \sqrt{17}-\frac{%
p}{q}\right| \\
&\leq &2\sqrt{17}+\frac{1}{c}.
\end{eqnarray*}
Thus the inequality 
\begin{equation}
\left| \sqrt{17}-\frac{p}{q}\right| >\frac{1}{cq^{2}}  \nonumber
\label{liou}
\end{equation}
is satisfied for all integer pair $(p,q),$ $q\geq 1,$ by any positive real
number $c$ satisfying 
\[
c>2\sqrt{17}+\frac{1}{c}, 
\]
i.e., 
\[
c>\sqrt{17}+\sqrt{18}=8.3657\ldots . 
\]
This completes the proof.\endproof

\begin{lemma}
\label{f1}Let $F_{1}(m)$ be a function of $m\in \Bbb{N}$ defined by 
\[
F_{1}(m)\equiv 192m^{3}\binom{m}{\left[ 0.7m\right] }\left( \frac{\lambda
_{2}}{\lambda _{2}-\lambda _{1}}\right) ^{\left[ 0.7m\right] }\left( \frac{%
-\lambda _{1}}{\lambda _{2}-\lambda _{1}}\right) ^{m-\left[ 0.7m\right] }.
\]
Then 
\[
F_{1}(k)\leq F_{1}(m)
\]
whenever 
\[
m\geq 100\quad \text{and}\quad k\geq m+11.
\]
\end{lemma}

\proof
We easily verify that 
\[
\frac{F_{1}(m+1)}{F_{1}(m)}=\left( \frac{\lambda _{2}}{\lambda _{2}-\lambda
_{1}}\right) \left( 1+\frac{1}{m}\right) ^{3}\frac{m+1}{\left[ 0.7m\right] +1%
} 
\]
whenever 
\[
\left[ 0.7m\right] \neq \left[ 0.7m+0.7\right] 
\]
and that 
\[
\frac{F_{1}(m+1)}{F_{1}(m)}=\left( \frac{-\lambda _{1}}{\lambda _{2}-\lambda
_{1}}\right) \left( 1+\frac{1}{m}\right) ^{3}\frac{m+1}{m-\left[ 0.7m\right]
+1} 
\]
whenever 
\[
\left[ 0.7m\right] =\left[ 0.7m+0.7\right] . 
\]
Then we obtain the following estimates: 
\[
\frac{F_{1}(m+1)}{F_{1}(m)}\leq \frac{10}{7}\left( \frac{\lambda _{2}}{%
\lambda _{2}-\lambda _{1}}\right) \left( 1+\frac{1}{m}\right) ^{4} 
\]
whenever 
\[
\left[ 0.7m\right] \neq \left[ 0.7m+0.7\right] 
\]
and that 
\[
\frac{F_{1}(m+1)}{F_{1}(m)}\leq \frac{10}{3}\left( \frac{-\lambda _{1}}{%
\lambda _{2}-\lambda _{1}}\right) \left( 1+\frac{1}{m}\right) ^{4} 
\]
whenever 
\[
\left[ 0.7m\right] =\left[ 0.7m+0.7\right] . 
\]
Hence we obtain 
\begin{eqnarray*}
\frac{F_{1}(m+11)}{F_{1}(m)} &=&\frac{F_{1}(m+11)}{F_{1}(m+10)}\times \cdots
\times \frac{F_{1}(m+1)}{F_{1}(m)} \\
&\leq &\frac{10^{10}}{3^{3}7^{7}}\left( \frac{\lambda _{2}}{\lambda
_{2}-\lambda _{1}}\right) ^{7}\left( \frac{-\lambda _{1}}{\lambda
_{2}-\lambda _{1}}\right) ^{3}\left( 1+\frac{11}{m}\right) ^{4}.
\end{eqnarray*}
Note that 
\[
\frac{10}{3}\left( \frac{-\lambda _{1}}{\lambda _{2}-\lambda _{1}}\right)
\leq \frac{10}{7}\left( \frac{\lambda _{2}}{\lambda _{2}-\lambda _{1}}%
\right) 
\]
and 
\begin{eqnarray*}
0.7\times 1 &=&0.7,\quad 0.7\times 2=1.4, \\
0.7\times 3 &=&2.1,\quad 0.7\times 4=2.8, \\
0.7\times 5 &=&3.5,\quad 0.7\times 6=4.2, \\
0.7\times 7 &=&4.9,\quad 0.7\times 8=5.6, \\
0.7\times 9 &=&6.3,\quad 0.7\times 0=0.
\end{eqnarray*}
Thus we have the following estimates: 
\begin{eqnarray*}
\frac{F_{1}(m+12)}{F_{1}(m)} &\leq &A^{8}B^{3}\left( 1+\frac{12}{m}\right)
^{4},\quad \frac{F_{1}(m+13)}{F_{1}(m)}\leq A^{8}B^{4}\left( 1+\frac{13}{m}%
\right) ^{4} \\
\frac{F_{1}(m+14)}{F_{1}(m)} &\leq &A^{8}B^{5}\left( 1+\frac{14}{m}\right)
^{4},\quad \frac{F_{1}(m+15)}{F_{1}(m)}\leq A^{9}B^{5}\left( 1+\frac{15}{m}%
\right) ^{4} \\
\frac{F_{1}(m+16)}{F_{1}(m)} &\leq &A^{9}B^{6}\left( 1+\frac{16}{m}\right)
^{4},\quad \frac{F_{1}(m+17)}{F_{1}(m)}\leq A^{9}B^{7}\left( 1+\frac{17}{m}%
\right) ^{4} \\
\frac{F_{1}(m+18)}{F_{1}(m)} &\leq &A^{10}B^{7}\left( 1+\frac{18}{m}\right)
^{4},\quad \frac{F_{1}(m+19)}{F_{1}(m)}\leq A^{10}B^{8}\left( 1+\frac{19}{m}%
\right) ^{4} \\
\frac{F_{1}(m+20)}{F_{1}(m)} &\leq &A^{10}B^{9}\left( 1+\frac{20}{m}\right)
^{4},\quad \frac{F_{1}(m+21)}{F_{1}(m)}\leq A^{10}B^{10}\left( 1+\frac{21}{m}%
\right) ^{4}
\end{eqnarray*}
where 
\[
A=\frac{10}{7}\left( \frac{\lambda _{2}}{\lambda _{2}-\lambda _{1}}\right)
,\quad B=\frac{10}{3}\left( \frac{-\lambda _{1}}{\lambda _{2}-\lambda _{1}}%
\right) . 
\]
Then the straight forward computation yields 
\begin{eqnarray*}
65 &\geq &\max \left\{ \frac{11}{A^{-\frac{7}{4}}B^{-\frac{3}{4}}-1},\text{%
\quad }\frac{12}{A^{-\frac{8}{4}}B^{-\frac{3}{4}}-1},\text{\quad }\frac{13}{%
A^{-\frac{8}{4}}B^{-\frac{4}{4}}-1},\right. \\
&&\qquad \left. \frac{14}{A^{-\frac{8}{4}}B^{-\frac{5}{4}}-1},\text{\quad }%
\frac{15}{A^{-\frac{9}{4}}B^{-\frac{5}{4}}-1},\text{\quad }\frac{16}{A^{-%
\frac{9}{4}}B^{-\frac{6}{4}}-1},\right. \\
&&\qquad \left. \frac{17}{A^{-\frac{9}{4}}B^{-\frac{7}{4}}-1},\text{\quad }%
\frac{18}{A^{-\frac{10}{4}}B^{-\frac{7}{4}}-1},\text{\quad }\frac{19}{A^{-%
\frac{10}{4}}B^{-\frac{8}{4}}-1},\right. \\
&&\qquad \left. \frac{20}{A^{-\frac{10}{4}}B^{-\frac{9}{4}}-1},\text{\quad }%
\frac{21}{A^{-\frac{10}{4}}B^{-\frac{10}{4}}-1}\right\}
\end{eqnarray*}
so that 
\[
F_{1}(k)\leq F_{1}(m) 
\]
whenever 
\[
m\geq 65\quad \text{and}\quad k\geq m+11. 
\]
This completes the proof.\endproof

\begin{lemma}
\label{nonsingular}Let $\Delta (m)$ be the function defined in Lemma \ref
{Lemm1}. Then 
\[
\Delta (m)^{-1}\neq 4,-\frac{4}{3}(m-3)\quad \text{for all }m.
\]
Thus the matrices $B_{m}(2),$ $m\geq 3,$ and $B_{m}(3),$ $m\geq 4,$ are
nonsingular.
\end{lemma}

\proof
We define a function $\delta _{m}$ as follows: 
\[
\Delta (m)^{-1}=(4-m)(1-\delta _{m}) 
\]
so that 
\begin{align}
& (4-m)\Delta (m)  \nonumber \\
& =(1-\delta _{m})^{-1}  \nonumber \\
& =\sum_{k=0}^{m}\frac{\binom{m}{k}\left( \frac{\lambda _{2}}{\lambda
_{2}-\lambda _{1}}\right) ^{k}\left( \frac{-\lambda _{1}}{\lambda
_{2}-\lambda _{1}}\right) ^{m-k}}{1-\frac{k}{m}\frac{m\lambda _{1}}{m-4}-(1-%
\frac{k}{m})\frac{m\lambda _{2}}{m-4}}  \tag*{(1.10)}  \label{average}
\end{align}
where 
\begin{eqnarray*}
\lambda _{1} &=&\frac{3-\sqrt{17}}{2}=-0.5615\ldots \\
\lambda _{2} &=&\frac{3+\sqrt{17}}{2}=3.5615\ldots \\
\frac{\lambda _{2}}{\lambda _{2}-\lambda _{1}} &=&\frac{3+\sqrt{17}}{2\sqrt{%
17}}=0.8638\ldots \\
\frac{-\lambda _{1}}{\lambda _{2}-\lambda _{1}} &=&\frac{-3+\sqrt{17}}{2%
\sqrt{17}}=0.1361\ldots .
\end{eqnarray*}
Note that 
\[
\sum_{k=0}^{m}\binom{m}{k}\left( \frac{\lambda _{2}}{\lambda _{2}-\lambda
_{1}}\right) ^{k}\left( \frac{-\lambda _{1}}{\lambda _{2}-\lambda _{1}}%
\right) ^{m-k}=\left( \frac{\lambda _{2}-\lambda _{1}}{\lambda _{2}-\lambda
_{1}}\right) ^{m}=1. 
\]
Thus the summation \ref{average} is an average of the function 
\begin{equation}
\frac{1}{1-\frac{m\lambda _{1}}{m-4}\frac{k}{m}-\frac{m\lambda _{2}}{m-4}(1-%
\frac{k}{m})}=\frac{m-4}{m\sqrt{17}}\left( \frac{k}{m}-\frac{1}{2}-\frac{m+8%
}{2m\sqrt{17}}\right) ^{-1}  \nonumber
\end{equation}
under the binary distribution 
\[
\binom{m}{k}\left( \frac{\lambda _{2}}{\lambda _{2}-\lambda _{1}}\right)
^{k}\left( \frac{-\lambda _{1}}{\lambda _{2}-\lambda _{1}}\right) ^{m-k}. 
\]
In the equation \ref{average}, the function 
\begin{equation}
\frac{1}{1-\frac{m\lambda _{1}}{m-4}X-\frac{m\lambda _{2}}{m-4}(1-X)} 
\tag*{(1.11)}  \label{singular}
\end{equation}
has a singular point at the value 
\[
X=\frac{m+8+m\sqrt{17}}{2m\sqrt{17}}. 
\]
But we easily see that 
\[
\frac{k}{m}\neq \frac{m+8+m\sqrt{17}}{2m\sqrt{17}} 
\]
for each $k=0,\cdots ,m.$

Then 
\begin{eqnarray*}
&&(1-\delta _{m})^{-1}-1 \\
&=&\delta _{m}(1-\delta _{m})^{-1} \\
&=&\sum_{k=0}^{m}\frac{\frac{k}{m}-\frac{1}{2}-\frac{3}{2\sqrt{17}}}{\frac{k%
}{m}-\frac{1}{2}-\frac{m+8}{2m\sqrt{17}}}\cdot \binom{m}{k}\left( \frac{%
\lambda _{2}}{\lambda _{2}-\lambda _{1}}\right) ^{k}\left( \frac{-\lambda
_{1}}{\lambda _{2}-\lambda _{1}}\right) ^{m-k}.
\end{eqnarray*}
For $m\geq 13,$ we have 
\begin{equation}
\frac{1}{2}+\frac{m+8}{2m\sqrt{17}}\leq 0.7\leq \frac{1}{2}+\frac{3}{2\sqrt{%
17}}  \tag*{(1.12)}  \label{inequality}
\end{equation}
so that 
\begin{eqnarray*}
&&\sum_{k=0}^{m}\frac{\left| \frac{k}{m}-\frac{1}{2}-\frac{3}{2\sqrt{17}}%
\right| }{\left| \frac{k}{m}-\frac{1}{2}-\frac{m+8}{2m\sqrt{17}}\right| }%
\cdot \binom{m}{k}\left( \frac{\lambda _{2}}{\lambda _{2}-\lambda _{1}}%
\right) ^{k}\left( \frac{-\lambda _{1}}{\lambda _{2}-\lambda _{1}}\right)
^{m-k} \\
&=&\sum_{k=0}^{\left[ 0.7m\right] }\frac{\left| \frac{k}{m}-\frac{1}{2}-%
\frac{3}{2\sqrt{17}}\right| }{\left| \frac{k}{m}-\frac{1}{2}-\frac{m+8}{2m%
\sqrt{17}}\right| }\cdot \binom{m}{k}\left( \frac{\lambda _{2}}{\lambda
_{2}-\lambda _{1}}\right) ^{k}\left( \frac{-\lambda _{1}}{\lambda
_{2}-\lambda _{1}}\right) ^{m-k} \\
&&+\sum_{k=\left[ 0.7m\right] +1}^{m}\frac{\left| \frac{k}{m}-\frac{1}{2}-%
\frac{3}{2\sqrt{17}}\right| }{\left| \frac{k}{m}-\frac{1}{2}-\frac{m+8}{2m%
\sqrt{17}}\right| }\cdot \binom{m}{k}\left( \frac{\lambda _{2}}{\lambda
_{2}-\lambda _{1}}\right) ^{k}\left( \frac{-\lambda _{1}}{\lambda
_{2}-\lambda _{1}}\right) ^{m-k}
\end{eqnarray*}
where the first summation contains the singular terms as $m\rightarrow
\infty .$

By Lemma \ref{17}, we have the following estimate: 
\begin{align}
& \frac{1}{\left| \frac{k}{m}-\frac{1}{2}-\frac{m+8}{2m\sqrt{17}}\right| } 
\nonumber \\
& =\frac{34m}{(m+8)\left| \sqrt{17}-\frac{17(2k-m)}{m+8}\right| }  \nonumber
\\
& <17^{2}m(m+8).  \tag*{(1.13)}  \label{estimate}
\end{align}
By the inequality \ref{inequality} and the estimate \ref{estimate}, we
obtain 
\begin{align}
& \sum_{k=0}^{m}\frac{\left| \frac{k}{m}-\frac{1}{2}-\frac{3}{2\sqrt{17}}%
\right| }{\left| \frac{k}{m}-\frac{1}{2}-\frac{m+8}{2m\sqrt{17}}\right| }%
\cdot \binom{m}{k}\left( \frac{\lambda _{2}}{\lambda _{2}-\lambda _{1}}%
\right) ^{k}\left( \frac{-\lambda _{1}}{\lambda _{2}-\lambda _{1}}\right)
^{m-k}  \nonumber \\
& \leq \frac{17\sqrt{17}(3+\sqrt{17})}{2}m(m+8)\sum_{k=0}^{\left[
0.7m\right] }\binom{m}{k}\left( \frac{\lambda _{2}}{\lambda _{2}-\lambda _{1}%
}\right) ^{k}\left( \frac{-\lambda _{1}}{\lambda _{2}-\lambda _{1}}\right)
^{m-k}  \nonumber \\
& +\left( \frac{1}{5}-\frac{m+8}{2m\sqrt{17}}\right) ^{-1}\sum_{k=\left[
0.7m\right] +1}^{m}\left| \frac{k}{m}-\frac{1}{2}-\frac{3}{2\sqrt{17}}\right|
\nonumber \\
& \hspace{2in}\times \binom{m}{k}\left( \frac{\lambda _{2}}{\lambda
_{2}-\lambda _{1}}\right) ^{k}\left( \frac{-\lambda _{1}}{\lambda
_{2}-\lambda _{1}}\right) ^{m-k}.  \tag*{(1.14)}  \label{home}
\end{align}
Note that the binary distribution 
\begin{equation}
\binom{m}{k}\left( \frac{\lambda _{2}}{\lambda _{2}-\lambda _{1}}\right)
^{k}\left( \frac{-\lambda _{1}}{\lambda _{2}-\lambda _{1}}\right) ^{m-k} 
\tag*{(1.15)}  \label{binary}
\end{equation}
increases up to the average 
\[
\frac{k}{m}=\frac{1}{2}+\frac{3}{2\sqrt{17}}\geq 0.7. 
\]
Thus we obtain, for $m\geq 100,$ 
\begin{eqnarray*}
&&\frac{17\sqrt{17}(3+\sqrt{17})m(m+8)}{2}\sum_{k=0}^{\left[ 0.7m\right] }%
\binom{m}{k}\left( \frac{\lambda _{2}}{\lambda _{2}-\lambda _{1}}\right)
^{k}\left( \frac{-\lambda _{1}}{\lambda _{2}-\lambda _{1}}\right) ^{m-k} \\
&\leq &192m^{3}\binom{m}{\left[ 0.7m\right] }\left( \frac{\lambda _{2}}{%
\lambda _{2}-\lambda _{1}}\right) ^{\left[ 0.7m\right] }\left( \frac{%
-\lambda _{1}}{\lambda _{2}-\lambda _{1}}\right) ^{m-\left[ 0.7m\right] }
\end{eqnarray*}
by the following inequality 
\[
\frac{17\sqrt{17}(3+\sqrt{17})}{2}\left[ 0.7m\right] m(m+8)\leq
192m^{3}\quad \text{for }m\geq 100. 
\]
By numerical computation, we obtain 
\begin{eqnarray*}
F_{1}(100) &=&2114.7\ldots \\
F_{1}(200) &=&1.5207\ldots \\
F_{1}(300) &\leq &5.3215\times 10^{-4}.
\end{eqnarray*}
Then, by Lemma \ref{f1}, we obtain the following numerical estimate: 
\begin{equation}
F_{1}(m)\leq 5.33\times 10^{-4}\quad \text{for }m\geq 400.  \tag*{(1.16)}
\label{f11}
\end{equation}

For the second part of the inequality \ref{home}, we have the following
estimate 
\begin{eqnarray*}
&&\sum_{k=0.7m}^{m}\left| \frac{k}{m}-\frac{1}{2}-\frac{3}{2\sqrt{17}}%
\right| \binom{m}{k}\left( \frac{\lambda _{2}}{\lambda _{2}-\lambda _{1}}%
\right) ^{k}\left( \frac{-\lambda _{1}}{\lambda _{2}-\lambda _{1}}\right)
^{m-k} \\
&\leq &\sqrt{\sum_{k=0}^{m}\left( \frac{k}{m}-\frac{1}{2}-\frac{3}{2\sqrt{17}%
}\right) ^{2}\binom{m}{k}\left( \frac{\lambda _{2}}{\lambda _{2}-\lambda _{1}%
}\right) ^{k}\left( \frac{-\lambda _{1}}{\lambda _{2}-\lambda _{1}}\right)
^{m-k}} \\
&=&\sqrt{\frac{2}{17m}}.
\end{eqnarray*}
We easily verify that 
\[
F_{2}(m)\equiv \left( \frac{1}{5}-\frac{m+8}{2m\sqrt{17}}\right) ^{-1}\sqrt{%
\frac{2}{17m}}\searrow 0\quad \text{as }m\rightarrow \infty . 
\]
By numerical computation, we obtain 
\begin{align}
F_{2}(400)& =0.2247\ldots  \nonumber \\
F_{2}(600)& =0.1815\ldots  \nonumber \\
F_{2}(800)& =0.1564\ldots .  \tag*{(1.17)}  \label{f12}
\end{align}
Note that we have the following estimate 
\[
\left| \delta _{m}\right| \leq \frac{\left| F(m)\right| }{1-\left|
F(m)\right| } 
\]
whenever 
\[
F(m)\equiv F_{1}(m)+F_{2}(m)<1. 
\]
Thus we obtain 
\[
\left| \delta _{m}\right| \leq 0.2 
\]
whenever 
\[
\left| F(m)\right| \leq \frac{1}{6}=0.1666\ldots . 
\]
Therefore, by the numerical result in \ref{f11} and \ref{f12}, 
\[
\left| \delta _{m}\right| \leq 0.2\quad \text{for all }m\geq 800. 
\]
Hence, it suffices to compute the numerical value of the function $\Delta
(m) $ up to $m\leq 800.$

Indeed, by numerical computation up to $m=800$, we can check the tendency of
the function $\delta _{m}$ as follows: 
\begin{equation}
\left| \delta _{m}\right| \leq 0.2  \tag*{(1.18)}  \label{computer1}
\end{equation}
for $m\geq 30$ so that 
\begin{equation}
\Delta (m)^{-1}\leq -20.  \nonumber
\end{equation}
Thus we easily see 
\[
\Delta (m)^{-1}\neq 4,-\frac{4}{3}(m-3)\quad \text{for }m\geq 30. 
\]
Then we need to check 
\[
\Delta (m)^{-1}=\frac{\det E_{m}(m+1)}{\det E_{m}(m)}\neq 4,-\frac{4}{3}%
(m-3)\quad \text{for }1\leq m\leq 29, 
\]
or, equivalently, 
\begin{eqnarray*}
\eta (m) &=&\frac{\det E_{m}(m)}{\det E_{m}(m-1)}=\frac{2m}{4-m-\Delta
(m)^{-1}} \\
&\neq &-2\text{ or }6\quad \text{for }1\leq m\leq 29.
\end{eqnarray*}
We may compute the value $\eta (m)$ for $1\leq m\leq 29$ by using the
following recurrence relation 
\[
\det E_{m}(s+1)=(2m+4-3s)\det E_{m}(s)-2s(m-s+1)\det E_{m}(s-1) 
\]
for $s=2,\cdots ,m,$ with the initial values 
\[
\det E_{m}(1)=2m+4\quad \text{and}\quad \det E_{m}(2)=4(m+1)^{2} 
\]
where 
\[
E_{m}(s+1)=\left( 
\begin{array}{ccccc}
2m-3s+4 & 2(m-s+1) & 0 & \cdots & 0 \\ 
s & 2m-3s+7 & 2(m-s+2) & \ddots & \vdots \\ 
0 & \ddots & \ddots & \ddots & 0 \\ 
\vdots & \ddots & 2 & 2m+1 & 2m \\ 
0 & \cdots & 0 & 1 & 2m+4
\end{array}
\right) . 
\]
Then we obtain 
\begin{equation}
\begin{tabular}{rrcrrr}
$\eta (1)=$ & $2.66\ldots $ & $\hspace{0.3cm}\eta (11)=$ & $-7.14\ldots $ & $%
\hspace{0.3cm}\eta (21)=$ & $-1.65\ldots $ \\ 
$\eta (2)=$ & $4.5\hspace{0.55cm}\;$ & $\hspace{0.3cm}\eta (12)=$ & $%
39.86\ldots $ & $\hspace{0.3cm}\eta (22)=$ & $-11.21\ldots $ \\ 
$\eta (3)=$ & $2.75\;\quad $ & $\hspace{0.3cm}\eta (13)=$ & $-1.96\ldots $ & 
$\hspace{0.3cm}\eta (23)=$ & $-31.81\ldots $ \\ 
$\eta (4)=$ & $0.36\ldots $ & $\hspace{0.3cm}\eta (14)=$ & $-9.84\ldots $ & $%
\hspace{0.3cm}\eta (24)=$ & $-7.43\ldots $ \\ 
$\eta (5)=$ & $-5.24\ldots $ & $\hspace{0.3cm}\eta (15)=$ & $19.83\ldots $ & 
$\hspace{0.3cm}\eta (25)=$ & $-14.96\ldots $ \\ 
$\eta (6)=$ & $12.05\ldots $ & $\hspace{0.3cm}\eta (16)=$ & $-4.61\ldots $ & 
$\hspace{0.3cm}\eta (26)=$ & $118.96\ldots $ \\ 
$\eta (7)=$ & $0.64\ldots $ & $\hspace{0.3cm}\eta (17)=$ & $-13.56\ldots $ & 
$\hspace{0.3cm}\eta (27)=$ & $-12.12\ldots $ \\ 
$\eta (8)=$ & $-5.36\ldots $ & $\hspace{0.3cm}\eta (18)=$ & $6.45\ldots $ & $%
\hspace{0.3cm}\eta (28)=$ & $-19.25\ldots $ \\ 
$\eta (9)=$ & $35.53\ldots $ & $\hspace{0.3cm}\eta (19)=$ & $-7.76\ldots $ & 
$\hspace{0.3cm}\eta (29)=$ & $-3.97\ldots $ \\ 
$\eta (10)=$ & $-0.11\ldots $ & $\hspace{0.3cm}\eta (20)=$ & $-19.12\ldots $
& $\hspace{0.3cm}\eta (30)=$ & $-16.30\ldots $%
\end{tabular}
\tag*{(1.19)}  \label{computer2}
\end{equation}
This completes the proof.\endproof

\section{Local automorphism group of a real hypersurface}

\subsection{Polynomial Identities}

We shall use the following notations: 
\begin{eqnarray*}
O(k+1) &=&\sum_{s+2t\geq k+1}O\left( \left| z\right| ^{s}\left| w\right|
^{t}\right) \\
O_{\times }(k+1) &=&\left( O(k),\cdots ,O(k),O(k+1)\right) .
\end{eqnarray*}

\begin{lemma}
\label{B}Let $M$ be a nondegenerate analytic real hypersurface defined near
the origin by 
\begin{equation}
v=F(z,\bar{z},u),\quad \left. F\right| _{0}=\left. dF\right| _{0}=0 
\nonumber
\end{equation}
and $\phi $ be a biholomorphic mapping near the origin such that the
transformed real hypersurface $\phi (M)$ is defined by the equation 
\begin{equation}
v=\langle z,z\rangle +F^{*}(z,\bar{z},u)+O(k+1).  \nonumber  \label{3.7}
\end{equation}
Suppose that the equation 
\[
v=\langle z,z\rangle +F^{*}(z,\bar{z},u)
\]
is in normal form. Then there is a normalization $\varphi $ of $M$ such that 
\[
\varphi =\phi +O_{\times }(k+1).
\]
Further, suppose that the normalization $\varphi $ transforms $M$ to a real
hypersurface $M^{\prime }$ in normal form defined by 
\[
v=\langle z,z\rangle +F^{\prime }(z,\bar{z},u).
\]
Then 
\[
F^{\prime }(z,\bar{z},u)=F^{*}(z,\bar{z},u)+O(k+1).
\]
\end{lemma}

In the paper \cite{Pa2}, we have given the proof of Lemma \ref{B}.

\begin{lemma}
\label{Lemma-a}Let $M$ be a real hypersurface in normal form defined by the
equation 
\[
v=\langle z,z\rangle +\sum_{\min (s,t)\geq 2}F_{st}(z,\overline{z},u)+O(l+2),
\]
where, for all complex number $\mu ,$%
\[
F_{st}(\mu z,\mu \overline{z},\mu ^{2}u)=\mu ^{l}F_{st}(z,\overline{z},u).
\]
Let $\phi $ be a normalization of $M$ with initial value $(id_{n\times
n},a,1,0)\in H$ such that $\phi $ transforms $M$ to a real hypersurface in
normal form defined by the equation 
\[
v=\langle z,z\rangle +F^{*}(z,\overline{z},u).
\]
Then 
\[
F^{*}(z,\overline{z},u)=\sum_{\min (s,t)\geq 2}F_{st}(z,\overline{z}%
,u)+O(l+2)
\]
if and only if 
\begin{align}
-2i(\langle z,a\rangle -\langle a,z\rangle )\sum_{\min (s,t)\geq 2}F_{st}(z,%
\overline{z},u)  \nonumber \\
+\sum_{\min (s-1,t)\geq 2}\sum_{\alpha }\left( \frac{\partial F_{st}}{%
\partial z^{\alpha }}\right) (z,\overline{z},u)a^{\alpha }(u+i\langle
z,z\rangle )  \nonumber \\
+2i\sum_{t\geq 2}\sum_{\alpha }\left( \frac{\partial F_{2t}}{\partial
z^{\alpha }}\right) (z,\overline{z},u)a^{\alpha }\langle z,z\rangle  
\nonumber \\
+\sum_{\min (s,t-1)\geq 2}\sum_{\alpha }\left( \frac{\partial F_{st}}{%
\partial \overline{z}^{\alpha }}\right) (z,\overline{z},u)\overline{a}%
^{\alpha }(u-i\langle z,z\rangle )  \nonumber \\
-2i\sum_{s\geq 2}\sum_{\alpha }\left( \frac{\partial F_{s2}}{\partial 
\overline{z}^{\alpha }}\right) (z,\overline{z},u)\overline{a}^{\alpha
}\langle z,z\rangle   \nonumber \\
+\frac{i}{2}\sum_{\min (s,t)\geq 2}\left( \frac{\partial F_{st}}{\partial u}%
\right) (z,\overline{z},u)\left\{ \langle z,a\rangle (u+i\langle z,z\rangle
)-\langle a,z\rangle (u-i\langle z,z\rangle )\right\}   \nonumber \\
+G_{l+1}(z,\overline{z},u)  \nonumber \\
=0  \tag*{(2.1)}  \label{lem-a}
\end{align}
where, for $l=2k-1,$%
\begin{eqnarray*}
G_{l+1}(z,\overline{z},u)=\frac{g}{2}\{(k-1)\langle z,z\rangle
+iu\}(u+i\langle z,z\rangle )^{k-1} \\
+\frac{g}{2}\{(k-1)\langle z,z\rangle -iu\}(u-i\langle z,z\rangle )^{k-1}
\end{eqnarray*}
and, for $l=2k,$%
\begin{eqnarray*}
G_{l+1}(z,\overline{z},u)=\langle \kappa ,z\rangle (u+i\langle z,z\rangle
)^{k}+\langle z,\kappa \rangle (u-i\langle z,z\rangle )^{k} \\
+2ik\langle z,z\rangle \langle z,\kappa \rangle (u+i\langle z,z\rangle
)^{k-1} \\
-2ik\langle z,z\rangle \langle \kappa ,z\rangle (u-i\langle z,z\rangle
)^{k-1} \\
-\langle z,\kappa \rangle (u+i\langle z,z\rangle )^{k}-\langle \kappa
,z\rangle (u-i\langle z,z\rangle )^{k}
\end{eqnarray*}
and 
\begin{eqnarray*}
\langle \kappa ,z\rangle =\frac{u^{3-k}}{4k(k-1)(n+1)(n+2)}\left\{
\sum_{\alpha }a^{\alpha }\Delta ^{2}\left( \frac{\partial F_{33}}{\partial
z^{\alpha }}\right) (z,\overline{z},u)\right.  \\
\hspace{2in}\left. +\sum_{\alpha }\overline{a}^{\alpha }\Delta ^{2}\left( 
\frac{\partial F_{24}}{\partial \overline{z}^{\alpha }}\right) (z,\overline{z%
},u)\right\}  \\
g=\frac{u^{4-k}}{2k(k-1)(k-2)n(n+1)(n+2)}\left\{ \sum_{\alpha }a^{\alpha
}\Delta ^{3}\left( \frac{\partial F_{43}}{\partial z^{\alpha }}\right) (z,%
\overline{z},u)\right.  \\
\hspace{2.2in}\left. +\sum_{\alpha }\overline{a}^{\alpha }\Delta ^{3}\left( 
\frac{\partial F_{34}}{\partial \overline{z}^{\alpha }}\right) (z,\overline{z%
},u)\right\} .
\end{eqnarray*}
\end{lemma}

\proof
For the initial value $(id_{n\times n},a,1,0)\in H,$ we have the following
decomposition(cf. \cite{Pa1}): 
\[
\phi =E\circ \psi 
\]
where $E$ is a normalization with identity initial value and 
\[
\psi :\left\{ 
\begin{array}{c}
z^{*}=\frac{z-aw}{1+2i\langle z,a\rangle -i\langle a,a\rangle w} \\ 
w^{*}=\frac{\rho w}{1+2i\langle z,a\rangle -i\langle a,a\rangle w}
\end{array}
\right. . 
\]
The mapping $\psi $ transforms $M$ to a real hypersurface $M^{\prime }$
defined up to $O(l+2)$ by the equation 
\begin{eqnarray*}
v &=&\langle z,z\rangle +F_{l}(z,z,u) \\
&&-2i(\langle z,a\rangle -\langle a,z\rangle )F_{l}(z,\overline{z},u) \\
&&+\sum_{\alpha }\left( \frac{\partial F_{l}}{\partial z^{\alpha }}\right)
(z,\overline{z},u)a^{\alpha }(u+i\langle z,z\rangle ) \\
&&+\sum_{\alpha }\left( \frac{\partial F_{l}}{\partial \overline{z}^{\alpha }%
}\right) (z,\overline{z},u)\overline{a}^{\alpha }(u-i\langle z,z\rangle ) \\
&&+\frac{i}{2}\left( \frac{\partial F_{l}}{\partial u}\right) (z,\overline{z}%
,u)\{\langle z,a\rangle (u+i\langle z,z\rangle )-\langle a,z\rangle
(u-i\langle z,z\rangle )\} \\
&&+O(l+2)
\end{eqnarray*}
where 
\[
F_{l}(z,z,u)=\sum_{\min (s,t)\geq 2}F_{st}(z,\overline{z},u). 
\]
By virtue of Lemma \ref{B}, we normalize $M^{\prime }$ up to $O(l+2)$ by a
mapping $h=(f,g)$ satisfying 
\begin{eqnarray*}
\left( \left. \frac{\partial f}{\partial z}\right| _{0}\right)
&=&id_{n\times n},\quad \left( \left. \frac{\partial f}{\partial w}\right|
_{0}\right) =0, \\
\Re \left( \left. \frac{\partial g}{\partial w}\right| _{0}\right)
&=&1,\quad \Re \left( \left. \frac{\partial ^{2}g}{\partial w^{2}}\right|
_{0}\right) =0,
\end{eqnarray*}
so that we obtain 
\begin{align}
F^{*}(z,\overline{z},u)& =F_{l}(z,\overline{z},u)-2i(\langle z,a\rangle
-\langle a,z\rangle )\sum_{\min (s,t)\geq 2}F_{st}(z,\overline{z},u) 
\nonumber \\
& +\sum_{\min (s-1,t)\geq 2}\sum_{\alpha }\left( \frac{\partial F_{st}}{%
\partial z^{\alpha }}\right) (z,\overline{z},u)a^{\alpha }(u+i\langle
z,z\rangle )  \nonumber \\
& +2i\sum_{t\geq 2}\sum_{\alpha }\left( \frac{\partial F_{2t}}{\partial
z^{\alpha }}\right) (z,\overline{z},u)a^{\alpha }\langle z,z\rangle 
\nonumber \\
& +\sum_{\min (s,t-1)\geq 2}\sum_{\alpha }\left( \frac{\partial F_{st}}{%
\partial \overline{z}^{\alpha }}\right) (z,\overline{z},u)\overline{a}%
^{\alpha }(u-i\langle z,z\rangle )  \nonumber \\
& -2i\sum_{s\geq 2}\sum_{\alpha }\left( \frac{\partial F_{s2}}{\partial 
\overline{z}^{\alpha }}\right) (z,\overline{z},u)\overline{a}^{\alpha
}\langle z,z\rangle  \nonumber \\
& +\frac{i}{2}\sum_{\min (s,t)\geq 2}\left( \frac{\partial F_{st}}{\partial u%
}\right) (z,\overline{z},u)\{\langle z,a\rangle (u+i\langle z,z\rangle ) 
\nonumber \\
& \hspace{5cm}-\langle a,z\rangle (u-i\langle z,z\rangle )\}  \nonumber \\
& +G_{l+1}(z,\overline{z},u)+O(l+2).  \tag*{(2.2)}  \label{sourse}
\end{align}
where, for $l=2k-1,$%
\begin{eqnarray*}
G_{l+1}(z,\overline{z},u) &=&\langle \chi z,z\rangle (u+i\langle z,z\rangle
)^{k-1}+\langle z,\chi z\rangle (u-i\langle z,z\rangle )^{k-1} \\
&&-\frac{g}{2i}(u+i\langle z,z\rangle )^{k}+\frac{g}{2i}(u-i\langle
z,z\rangle )^{k}
\end{eqnarray*}
and, for $l=2k,$%
\begin{eqnarray*}
G_{l+1}(z,\overline{z},u) &=&\langle \kappa ,z\rangle (u+i\langle z,z\rangle
)^{k}+\langle z,\kappa \rangle (u-i\langle z,z\rangle )^{k} \\
&&+2ik\langle z,z\rangle \langle z,\kappa \rangle (u+i\langle z,z\rangle
)^{k-1} \\
&&-2ik\langle z,z\rangle \langle \kappa ,z\rangle (u-i\langle z,z\rangle
)^{k-1} \\
&&-\langle z,\kappa \rangle (u+i\langle z,z\rangle )^{k}-\langle \kappa
,z\rangle (u-i\langle z,z\rangle )^{k}.
\end{eqnarray*}
Here the constants $\chi ,g,\kappa $ satisfy the conditions 
\begin{gather*}
\langle \chi z,z\rangle +\langle z,\chi z\rangle =kg\langle z,z\rangle , \\
g\in \Bbb{R},\quad \kappa \in \Bbb{C}^{n},
\end{gather*}
and they are uniquely determined by the following conditions: 
\[
\Delta F_{22}^{*}=\Delta ^{2}F_{23}^{*}=\Delta ^{3}F_{33}^{*}=0, 
\]
where 
\begin{eqnarray*}
F^{*}(z,\overline{z},u) &=&\sum_{\min (s,t)\geq 2}F_{st}(z,\overline{z},u) \\
&&+\sum_{\min (s,t)\geq 2}F_{st}^{*}(z,\overline{z},u)+O(l+2)
\end{eqnarray*}
and, for all complex number $\mu ,$%
\[
F_{st}^{*}(\mu z,\mu \overline{z},\mu ^{2}u)=\mu ^{l+1}F_{st}^{*}(z,%
\overline{z},u). 
\]
Indeed, from the equality \ref{sourse}, we obtain 
\begin{eqnarray*}
F_{22}^{*}(z,\overline{z},u) &=&2(k-1)i\langle \chi z,z\rangle \langle
z,z\rangle u^{k-2}-k(k-1)ig\langle z,z\rangle ^{2}u^{k-2} \\
F_{23}^{*}(z,\overline{z},u) &=&-2k(k-1)\langle \kappa ,z\rangle \langle
z,z\rangle ^{2}u^{k-2}+2i\langle a,z\rangle F_{22}(z,\overline{z},u) \\
&&+2i\sum_{\alpha }\left( \frac{\partial F_{22}}{\partial z^{\alpha }}%
\right) (z,\overline{z},u)a^{\alpha }\langle z,z\rangle -\frac{i}{2}\left( 
\frac{\partial F_{22}}{\partial u}\right) (z,\overline{z},u)\langle
a,z\rangle u \\
&&+\sum_{\alpha }\left( \frac{\partial F_{33}}{\partial z^{\alpha }}\right)
(z,\overline{z},u)a^{\alpha }u+\sum_{\alpha }\left( \frac{\partial F_{24}}{%
\partial \overline{z}^{\alpha }}\right) (z,\overline{z},u)\overline{a}%
^{\alpha }u \\
F_{33}^{*}(z,\overline{z},u) &=&-\frac{k(k-1)(k-2)}{3}g\langle z,z\rangle
^{3}u^{k-3} \\
&&-2i\langle z,a\rangle F_{23}(z,\overline{z},u)+2i\langle a,z\rangle
F_{32}(z,\overline{z},u) \\
&&+\sum_{\alpha }\left( \frac{\partial F_{43}}{\partial z^{\alpha }}\right)
(z,\overline{z},u)a^{\alpha }u+i\sum_{\alpha }\left( \frac{\partial F_{32}}{%
\partial z^{\alpha }}\right) (z,\overline{z},u)a^{\alpha }\langle z,z\rangle
\\
&&+\sum_{\alpha }\left( \frac{\partial F_{34}}{\partial \overline{z}^{\alpha
}}\right) (z,\overline{z},u)\overline{a}^{\alpha }u-i\sum_{\alpha }\left( 
\frac{\partial F_{23}}{\partial \overline{z}^{\alpha }}\right) (z,\overline{z%
},u)\overline{a}^{\alpha }\langle z,z\rangle \\
&&+\frac{i}{2}\left( \frac{\partial F_{23}}{\partial u}\right) (z,\overline{z%
},u)\langle z,a\rangle u-\frac{i}{2}\left( \frac{\partial F_{32}}{\partial u}%
\right) (z,\overline{z},u)\langle a,z\rangle u.
\end{eqnarray*}
Hence we obtain 
\begin{eqnarray*}
\Delta F_{22}^{*}(z,\overline{z},u) &=&2(k-1)(n+2)i\langle \chi z,z\rangle
u^{k-2}+2(k-1)i\mathrm{Tr(}\chi \mathrm{)}\langle z,z\rangle u^{k-2} \\
&&-2k(k-1)(n+1)ig\langle z,z\rangle u^{k-2} \\
\Delta ^{2}F_{22}^{*}(z,\overline{z},u) &=&4(k-1)(n+1)i\mathrm{Tr(}\chi 
\mathrm{)}u^{k-2}-2k(k-1)n(n+1)igu^{k-2} \\
\Delta ^{2}F_{23}^{*}(z,\overline{z},u) &=&-4k(k-1)(n+1)(n+2)\langle \kappa
,z\rangle u^{k-2} \\
&&+\sum_{\alpha }ua^{\alpha }\Delta ^{2}\left( \frac{\partial F_{33}}{%
\partial z^{\alpha }}\right) (z,\overline{z},u)+\sum_{\alpha }u\overline{a}%
^{\alpha }\Delta ^{2}\left( \frac{\partial F_{24}}{\partial \overline{z}%
^{\alpha }}\right) (z,\overline{z},u) \\
\Delta ^{3}F_{33}^{*}(z,\overline{z},u) &=&-2k(k-1)(k-2)n(n+1)(n+2)gu^{k-3}
\\
&&+\sum_{\alpha }ua^{\alpha }\Delta ^{3}\left( \frac{\partial F_{43}}{%
\partial z^{\alpha }}\right) (z,\overline{z},u)+\sum_{\alpha }u\overline{a}%
^{\alpha }\Delta ^{3}\left( \frac{\partial F_{34}}{\partial \overline{z}%
^{\alpha }}\right) (z,\overline{z},u)
\end{eqnarray*}
Note that the condition $\Delta F_{22}^{*}=0$ yields 
\begin{equation}
2\langle \chi z,z\rangle =kg\langle z,z\rangle .  \tag*{(2.3)}
\label{sourse4}
\end{equation}
The condition $\Delta ^{2}F_{23}^{*}=\Delta ^{3}F_{33}^{*}=0,$ with the
equality \ref{sourse4}, uniquely determines the constants $\chi ,\kappa ,g.$

Then we easily see that 
\[
F^{*}(z,\overline{z},u)=\sum_{\min (s,t)\geq 2}F_{st}(z,\overline{z}%
,u)+O(l+2) 
\]
if and only if 
\[
\sum_{\min (s,t)\geq 2}F_{st}^{*}(z,\overline{z},u)=0. 
\]
This completes the proof.\endproof

Note that, for odd integer $l,$%
\[
G_{l+1}(z,\overline{z},u)=\sum_{\min (s,t)\geq 2,s=t}G_{st}(z,\overline{z}%
,u) 
\]
and, for even integer $l,$%
\[
G_{l+1}(z,\overline{z},u)=\sum_{\min (s,t)\geq 2,s=t\pm 1}G_{st}(z,\overline{%
z},u). 
\]

\begin{lemma}
Let $M$ be a real hypersurface in normal form defined by the equation 
\[
v=\langle z,z\rangle +F_{l}(z,\overline{z},u)+O(l+2)
\]
where 
\[
F_{l}(z,\overline{z},u)=\sum_{s\geq 2}F_{ss}(z,\overline{z},u)
\]
and, for all complex number $\mu $, 
\[
F_{l}(\mu z,\mu \overline{z},\mu ^{2}u)=\mu ^{l}F_{l}(z,\overline{z},u).
\]
Let $\phi $ be a normalization with initial value $(id_{n\times n},a,1,0)\in
H$ such that $\phi $ transforms $M$ to a real hypersurface $M^{\prime }$ in
normal form defined by the equation 
\[
v=\langle z,z\rangle +F^{*}(z,\overline{z},u).
\]
Suppose that 
\[
F^{*}(z,\overline{z},u)=F_{l}(z,\overline{z},u)+O(l+2).
\]
Then there is an identity, for each integer $s,$ $3\leq s\leq k$, as
follows: 
\[
\sum_{s=2}^{k}i^{k-s}\langle z,z\rangle ^{k-s}\sum_{\alpha }a^{\alpha }\frac{%
\partial }{\partial z^{\alpha }}\left\{ \langle z,z\rangle \left( \frac{%
F_{ss}(z,\overline{z},u)}{u^{k-s}}\right) \right\} =-(2i)^{k-1}\langle
\kappa ,z\rangle \langle z,z\rangle ^{k}
\]
where 
\[
l=2k.
\]
\end{lemma}

\proof
By the condition 
\[
F_{l}(z,\overline{z},u)=\sum_{s\geq 2}F_{ss}(z,\overline{z},u), 
\]
the identity \ref{lem-a} in Lemma \ref{Lemma-a} comes to 
\begin{align}
& 2i(\langle z,a\rangle -\langle a,z\rangle )\sum_{s\geq 2}F_{ss}(z,%
\overline{z},u)  \nonumber \\
& -2i\langle z,z\rangle \sum_{\alpha }\left\{ \left( \frac{\partial F_{22}}{%
\partial z^{\alpha }}\right) (z,\overline{z},u)a^{\alpha }-\left( \frac{%
\partial F_{22}}{\partial \overline{z}^{\alpha }}\right) (z,\overline{z},u)%
\overline{a}^{\alpha }\right\}  \nonumber \\
& -\sum_{s\geq 3}\sum_{\alpha }\left\{ \left( \frac{\partial F_{ss}}{%
\partial z^{\alpha }}\right) (z,\overline{z},u)a^{\alpha }(u+i\langle
z,z\rangle )\right.  \nonumber \\
& \hspace{3cm}\left. +\left( \frac{\partial F_{ss}}{\partial \overline{z}%
^{\alpha }}\right) (z,\overline{z},u)\overline{a}^{\alpha }(u-i\langle
z,z\rangle )\right\}  \nonumber \\
& -\frac{i}{2}\sum_{s\geq 2}\left( \frac{\partial F_{ss}}{\partial u}\right)
(z,\overline{z},u)\{\langle z,a\rangle (u+i\langle z,z\rangle )-\langle
a,z\rangle (u-i\langle z,z\rangle )\}  \nonumber \\
& =\langle \kappa ,z\rangle (u+i\langle z,z\rangle )^{k}+\langle z,\kappa
\rangle (u-i\langle z,z\rangle )^{k}  \nonumber \\
& +2ik\langle z,z\rangle \langle z,\kappa \rangle (u+i\langle z,z\rangle
)^{k-1}-2ik\langle z,z\rangle \langle \kappa ,z\rangle (u-i\langle
z,z\rangle )^{k-1}  \nonumber \\
& -\langle z,\kappa \rangle (u+i\langle z,z\rangle )^{k}-\langle \kappa
,z\rangle (u-i\langle z,z\rangle )^{k}  \nonumber \\
& =\langle \kappa ,z\rangle \sum_{t=2}^{k}\{1+(-1)^{t}(2t-1)\}\binom{k}{t}%
u^{k-t}(i\langle z,z\rangle )^{t}  \nonumber \\
& +\langle z,\kappa \rangle \sum_{t=2}^{k}\{1+(-1)^{t}(2t-1)\}\binom{k}{t}%
u^{k-t}(-i\langle z,z\rangle )^{t}  \tag*{(2.4)}  \label{lem-b}
\end{align}
Then, by Lemma \ref{Lemma-a}, the constant $\kappa $ is given by 
\begin{equation}
\langle \kappa ,z\rangle =\frac{u^{2-k}}{4k(k-1)(n+1)(n+2)}\sum_{\alpha
}ua^{\alpha }\Delta ^{2}\left( \frac{\partial F_{33}}{\partial z^{\alpha }}%
\right) (z,\overline{z},u).  \nonumber  \label{b5}
\end{equation}

By collecting functions of type $(m+2,m+3)$ for $m=0,\cdots ,k-2$ in the
identity \ref{lem-b}, we obtain the following identities for each integer $%
s, $ $3\leq s\leq k$: 
\begin{align}
& \langle z,z\rangle \left\{ \langle a,z\rangle \left( \frac{\partial
F_{s-1,s-1}}{\partial u}\right) (z,\overline{z},u)+2i\sum_{\alpha }a^{\alpha
}\left( \frac{\partial F_{ss}}{\partial z^{\alpha }}\right) (z,\overline{z}%
,u)\right\}  \nonumber \\
& =4i\langle a,z\rangle F_{ss}(z,\overline{z},u)+4i\langle z,z\rangle
\sum_{\alpha }a^{\alpha }\left( \frac{\partial F_{ss}}{\partial z^{\alpha }}%
\right) (z,\overline{z},u)  \nonumber \\
& +2\langle \kappa ,z\rangle \{1+(-1)^{s}(2s-1)\}\binom{k}{s}%
u^{k-s}(i\langle z,z\rangle )^{s}  \nonumber \\
& -iu\left\{ \langle a,z\rangle \left( \frac{\partial F_{ss}}{\partial u}%
\right) (z,\overline{z},u)+2i\sum_{\alpha }a^{\alpha }\left( \frac{\partial
F_{s+1,s+1}}{\partial z^{\alpha }}\right) (z,\overline{z},u)\right\} 
\tag*{(2.5)}  \label{6.12}
\end{align}
and 
\begin{align}
& 4i\langle a,z\rangle F_{22}(z,\overline{z},u)+4i\langle z,z\rangle
\sum_{\alpha }a^{\alpha }\left( \frac{\partial F_{22}}{\partial z^{\alpha }}%
\right) (z,\overline{z},u)  \nonumber \\
& -4k(k-1)\langle \kappa ,z\rangle u^{k-2}\langle z,z\rangle ^{2}  \nonumber
\\
& -iu\left\{ \langle a,z\rangle \left( \frac{\partial F_{22}}{\partial u}%
\right) (z,\overline{z},u)+2i\sum_{\alpha }a^{\alpha }\left( \frac{\partial
F_{33}}{\partial z^{\alpha }}\right) (z,\overline{z},u)\right\}  \nonumber \\
& =0.  \tag*{(2.6)}  \label{6.11}
\end{align}
In the equality \ref{6.12}, we assume 
\[
F_{k+1,k+1}(z,\overline{z},u)=0. 
\]
From the equalities \ref{6.12} and \ref{6.11}, we obtain the following
recurrence relation: 
\begin{eqnarray*}
A(s) &=&iu^{-1}\langle z,z\rangle A(s-1) \\
&&+4u^{-1}\sum_{\alpha }a^{\alpha }\frac{\partial }{\partial z^{\alpha }}%
\left\{ \langle z,z\rangle F_{ss}(z,\overline{z},u)\right\} \\
&&-2i\langle \kappa ,z\rangle \{1+(-1)^{s}(2s-1)\}\binom{k}{s}%
u^{k-s-1}(i\langle z,z\rangle )^{s}
\end{eqnarray*}
for $s=2,\cdots ,k,$ and 
\[
A(1)=A(k)=0. 
\]
Thus we obtain the following identity: 
\begin{align}
& \sum_{s=2}^{k}i^{k-s}\langle z,z\rangle ^{k-s}\sum_{\alpha }a^{\alpha }%
\frac{\partial }{\partial z^{\alpha }}\left\{ \langle z,z\rangle \left( 
\frac{F_{ss}(z,\overline{z},u)}{u^{k-s}}\right) \right\}  \nonumber \\
& =\frac{i^{k+1}}{2}\langle \kappa ,z\rangle \langle z,z\rangle
^{k}\sum_{s=2}^{k}\{1+(-1)^{s}(2s-1)\}\binom{k}{s}  \nonumber \\
& =-(2i)^{k-1}\langle \kappa ,z\rangle \langle z,z\rangle ^{k}.  \tag*{(2.7)}
\label{a.20}
\end{align}
This completes the proof.\endproof

\begin{lemma}
\label{divided}Suppose that the functions $F_{ss}(z,\overline{z},u),$ $%
s=2,\cdots ,k,$ satisfy the equalities \ref{6.12}, where $l=2k$. Then that
the polynomial 
\[
F_{ss}(z,\overline{z},u),\quad s=\max (k-m,2),\cdots ,k-1,
\]
is divided by $\langle z,z\rangle ^{m-k+s}$ whenever 
\[
a\neq 0
\]
and $F_{kk}(z,\overline{z},u)$ is divided by $\langle z,z\rangle ^{m}$ for $%
0\leq m\leq k.$
\end{lemma}

\proof
The equality \ref{6.12} yields, for $s=3,\cdots ,k,$%
\begin{eqnarray*}
&&(k-s+1)\langle z,z\rangle \langle a,z\rangle F_{s-1,s-1}(z,\overline{z},u)
\\
&=&-i(k-s-4)\langle a,z\rangle F_{ss}(z,\overline{z},u)+2i\langle z,z\rangle
\sum_{\alpha }a^{\alpha }\left( \frac{\partial F_{ss}}{\partial z^{\alpha }}%
\right) (z,\overline{z},u) \\
&&+2\langle \kappa ,z\rangle \langle z,z\rangle ^{s}\{1+(-1)^{s}(2s-1)\}%
\binom{k}{s}i^{s}u^{k-s} \\
&&+2u\sum_{\alpha }a^{\alpha }\left( \frac{\partial F_{s+1,s+1}}{\partial
z^{\alpha }}\right) (z,\overline{z},u).
\end{eqnarray*}
Since $\langle a,z\rangle $ is not a devisor of $\langle z,z\rangle ,$ this
equality yields the desired result. This completes the proof.\endproof

\begin{lemma}
Let $M$ be a real hypersurface in normal form defined by the equation 
\[
v=\langle z,z\rangle +F_{l}(z,\overline{z},u)+O(l+2)
\]
where 
\[
F_{l}(z,\overline{z},u)=\sum_{\min (s,t)\geq 2}F_{st}(z,\overline{z},u)
\]
and, for all complex number $\mu $, 
\[
F_{l}(\mu z,\mu \overline{z},\mu ^{2}u)=\mu ^{l}F_{l}(z,\overline{z},u).
\]
Let $\phi $ be a normalization with initial value $(id_{n\times n},a,1,0)\in
H$ such that $\phi $ transforms $M$ to a real hypersurface $M^{\prime }$ in
normal form defined by the equation 
\[
v=\langle z,z\rangle +F^{*}(z,\overline{z},u).
\]
Suppose that 
\[
F^{*}(z,\overline{z},u)=F_{l}(z,\overline{z},u)+O(l+2)
\]
and the function $F_{l}(z,\overline{z},u)$ contains a nonzero function $%
F_{st}(z,\overline{z},u)$ of type $(s,t),$ $s\neq t.$ Then there is an
identity, for each integer $s,$ $3\leq s\leq p$, as follows: 
\begin{equation}
\sum_{s=2}^{p}i^{p-s}\langle z,z\rangle ^{p-s}\sum_{\alpha }a^{\alpha }\frac{%
\partial }{\partial z^{\alpha }}\left\{ \langle z,z\rangle \left( \frac{%
F_{s,l-2p+s}(z,\overline{z},u)}{u^{p-s}}\right) \right\} =0  \nonumber
\end{equation}
where 
\[
l-2p=\max \left\{ \left| t-s\right| :F_{st}(z,\overline{z},u)\neq 0\right\} .
\]
\end{lemma}

\proof
We easily verify that $p$ is an integer satisfying 
\begin{equation}
2\leq p\leq \left[ \frac{l-1}{2}\right] .  \tag*{(2.8)}  \label{4.56}
\end{equation}
By collecting functions of type $(s,t)$ satisfying 
\[
t-s=l-2p+1 
\]
in the identity \ref{lem-a} in Lemma \ref{Lemma-a}, we obtain the following
identities for each integer $s,$ $3\leq s\leq p$: 
\begin{align}
& \langle z,z\rangle \left\{ \langle a,z\rangle \left( \frac{\partial
F_{s-1,l-2p+s-1}}{\partial u}\right) (z,\overline{z},u)+2i\sum_{\alpha
}a^{\alpha }\left( \frac{\partial F_{s,l-2p+s}}{\partial z^{\alpha }}\right)
(z,\overline{z},u)\right\}  \nonumber \\
& =4i\langle a,z\rangle F_{s,l-2p+s}(z,\overline{z},u)+4i\langle z,z\rangle
\sum_{\alpha }a^{\alpha }\left( \frac{\partial F_{s,l-2p+s}}{\partial
z^{\alpha }}\right) (z,\overline{z},u)  \nonumber \\
& -iu\left\{ \langle a,z\rangle \left( \frac{\partial F_{s,l-2p+s}}{\partial
u}\right) (z,\overline{z},u)+2i\sum_{\alpha }a^{\alpha }\left( \frac{%
\partial F_{s+1,l-2p+s+1}}{\partial z^{\alpha }}\right) (z,\overline{z}%
,u)\right\}  \tag*{(2.9)}  \label{4.58}
\end{align}
and 
\begin{align}
& 4i\langle a,z\rangle F_{2,l-2p+2}(z,\overline{z},u)+4i\langle z,z\rangle
\sum_{\alpha }a^{\alpha }\left( \frac{\partial F_{2,l-2p+2}}{\partial
z^{\alpha }}\right) (z,\overline{z},u)  \nonumber \\
& -iu\left\{ \langle a,z\rangle \left( \frac{\partial F_{2,l-2p+2}}{\partial
u}\right) (z,\overline{z},u)+2i\sum_{\alpha }a^{\alpha }\left( \frac{%
\partial F_{3,l-2p+3}}{\partial z^{\alpha }}\right) (z,\overline{z}%
,u)\right\}  \nonumber \\
& =0.  \tag*{(2.10)}  \label{4.55}
\end{align}
In the equality \ref{4.58}, we assume 
\[
F_{p+1,l-p+1}(z,\overline{z},u)=0. 
\]
From the equalities \ref{4.58} and \ref{4.55}, we obtain the following
recurrence relation: 
\begin{eqnarray*}
A(s) &=&iu^{-1}\langle z,z\rangle A(s-1) \\
&&+4u^{-1}\sum_{\alpha }a^{\alpha }\frac{\partial }{\partial z^{\alpha }}%
\left\{ \langle z,z\rangle F_{s,l-2p+s}(z,\overline{z},u)\right\}
\end{eqnarray*}
for $s=2,\cdots ,p,$ and 
\[
A(1)=A(p)=0. 
\]
Thus we obtain the following identity: 
\begin{equation}
\sum_{s=2}^{p}i^{p-s}\langle z,z\rangle ^{p-s}\sum_{\alpha }a^{\alpha }\frac{%
\partial }{\partial z^{\alpha }}\left\{ \langle z,z\rangle \left( \frac{%
F_{s,l-2p+s}(z,\overline{z},u)}{u^{p-s}}\right) \right\} =0.  \tag*{(2.11)}
\label{a.19}
\end{equation}
This completes the proof.\endproof

\begin{lemma}
\label{divided2}Suppose that the functions $F_{st}(z,\overline{z},u)$
satisfy the equalities \ref{4.58}. Then the polynomial 
\[
F_{s,l-2p+s}(z,\overline{z},u),\quad s=\max (p-m,2),\cdots ,p-1,
\]
is divided by $\langle z,z\rangle ^{m-p+s}$ whenever 
\[
a\neq 0
\]
and $F_{p,l-p}(z,\overline{z},u)$ is divided by $\langle z,z\rangle ^{m}$
for $1\leq m\leq p.$
\end{lemma}

\proof
The equality \ref{4.58} yields, for $s=3,\cdots ,p,$%
\begin{eqnarray*}
&&(p-s+1)\langle z,z\rangle \langle a,z\rangle F_{s-1,l-2p+s-1}(z,\overline{z%
},u) \\
&=&-(p-s-4)i\langle a,z\rangle F_{s,l-2p+s}(z,\overline{z},u)+2i\langle
z,z\rangle \sum_{\alpha }a^{\alpha }\left( \frac{\partial F_{s,l-2p+s}}{%
\partial z^{\alpha }}\right) (z,\overline{z},u) \\
&&+2u\sum_{\alpha }a^{\alpha }\left( \frac{\partial F_{s+1,l-2p+s+1}}{%
\partial z^{\alpha }}\right) (z,\overline{z},u).
\end{eqnarray*}
Since $\langle a,z\rangle $ is not a divisor of $\langle z,z\rangle ,$ this
equality yields the desired result. This completes the proof.\endproof

\subsection{Injectivity of a Linear Mapping}

\begin{lemma}
\label{li-in}Let $l$ be a positive integer$\geq 4$ and $F_{2,l-2}(z,%
\overline{z},0)$ be a nonzero function of type $(2,l-2).$ Then the following
functions 
\[
H_{\alpha }(z,\overline{z},0)=\frac{\partial }{\partial z^{\alpha }}\left\{
\langle z,z\rangle F_{2,l-2}(z,\overline{z},0)\right\} \quad \text{for }%
\alpha =1,.\cdots ,n,
\]
are linearly independent.
\end{lemma}

\proof
Suppose the functions $H_{1}(z,\overline{z},0),\cdots ,H_{n}(z,\overline{z}%
,0)$ are linearly dependent over $\Bbb{C}$. Then there is a nonzero vector $%
a=(a^{\alpha })\in \Bbb{C}^{n}$ such that 
\[
\sum_{\alpha }a^{\alpha }\frac{\partial }{\partial z^{\alpha }}\{\langle
z,z\rangle F_{2,l-2}(z,\overline{z},0)\}=0. 
\]
Then we obtain 
\begin{equation}
\langle a,z\rangle F_{2,l-2}(z,\overline{z},0)=\langle z,z\rangle
\sum_{\alpha }a^{\alpha }\left( \frac{\partial F_{2,l-2}}{\partial z^{\alpha
}}\right) (z,\overline{z},0).  \tag*{(2.12)}  \label{b1}
\end{equation}
Note that $\langle a,z\rangle $ is not a devisor of $\langle z,z\rangle $
whenever $a\neq 0.$ Otherwise there would be a vector $b\in \Bbb{C}^{n}$ so
that 
\[
\langle z,z\rangle =\langle a,z\rangle \langle z,b\rangle . 
\]
This is a contradiction to the fact that the hermitian form $\langle
z,z\rangle $ is nondegenerate. Hence the polynomial 
\[
F_{2,l-2}(z,\overline{z},0) 
\]
is divided by $\langle z,z\rangle $ so that there is a polynomial $%
G_{1,l-3}(z,\overline{z},0)$ of type $(1,l-3)$ as follows: 
\[
F_{2,l-2}(z,\overline{z},0)=\langle z,z\rangle G_{1,l-3}(z,\overline{z},0). 
\]
Then the equality \ref{b1} comes to 
\[
2\langle a,z\rangle G_{1,l-3}(z,\overline{z},0)=\langle z,z\rangle
\sum_{\alpha }a^{\alpha }\left( \frac{\partial G_{1,l-3}}{\partial z^{\alpha
}}\right) (z,\overline{z},0). 
\]
Note that $G_{1,l-3}(z,\overline{z},0)$ is divided by $\langle z,z\rangle $
as well so that there is a polynomial $G_{0,l-4}(z,\overline{z},0)$ as
follows: 
\[
F_{2,l-2}(z,\overline{z},0)=\langle z,z\rangle ^{2}G_{0,l-4}(z,\overline{z}%
,0). 
\]
Then the equality \ref{b1} comes to 
\begin{equation}
\langle a,z\rangle G_{0,l-4}(z,\overline{z},0)=0.  \tag*{(2.13)}  \label{wed}
\end{equation}
Note that $\langle a,z\rangle \neq 0$ unless $a=0.$ Thus the equality \ref
{wed} yields 
\[
F_{2,l-2}(z,\overline{z},0)=0. 
\]
This is a contradiction to the assumption $F_{2,l-2}(z,\overline{z},0)\neq
0. $ This completes the proof.\endproof

\begin{lemma}
\label{Theo1}Suppose that 
\[
F_{l}(z,\overline{z},u)=\sum_{\min (s,t)\geq 2}F_{st}(z,\overline{z},u)
\]
where 
\[
F_{l}(\mu z,\mu \overline{z},\mu ^{2}u)=\mu ^{l}F_{l}(z,\overline{z},u)
\]
and 
\[
\Delta F_{22}=\Delta ^{2}F_{23}=\Delta ^{3}F_{33}=0.
\]
Then the linear mapping 
\[
a\longmapsto H_{l+1}(z,\overline{z},u;a)
\]
is injective, where 
\begin{eqnarray*}
H_{l+1}(z,\overline{z},u;a)\equiv -2i(\langle z,a\rangle -\langle a,z\rangle
)\sum_{\min (s,t)\geq 2}F_{st}(z,\overline{z},u) \\
+\sum_{\min (s-1,t)\geq 2}\sum_{\alpha }\left( \frac{\partial F_{st}}{%
\partial z^{\alpha }}\right) (z,\overline{z},u)a^{\alpha }(u+i\langle
z,z\rangle ) \\
+2i\sum_{t\geq 2}\sum_{\alpha }\left( \frac{\partial F_{2t}}{\partial
z^{\alpha }}\right) (z,\overline{z},u)a^{\alpha }\langle z,z\rangle  \\
+\sum_{\min (s,t-1)\geq 2}\sum_{\alpha }\left( \frac{\partial F_{st}}{%
\partial \overline{z}^{\alpha }}\right) (z,\overline{z},u)\overline{a}%
^{\alpha }(u-i\langle z,z\rangle ) \\
-2i\sum_{s\geq 2}\sum_{\alpha }\left( \frac{\partial F_{s2}}{\partial 
\overline{z}^{\alpha }}\right) (z,\overline{z},u)\overline{a}^{\alpha
}\langle z,z\rangle  \\
+\frac{i}{2}\sum_{\min (s,t)\geq 2}\left( \frac{\partial F_{st}}{\partial u}%
\right) (z,\overline{z},u)\{\langle z,a\rangle (u+i\langle z,z\rangle ) \\
\hspace{5cm}-\langle a,z\rangle (u-i\langle z,z\rangle )\} \\
+G_{l+1}(z,\overline{z},u)
\end{eqnarray*}
and $G_{l+1}(z,\overline{z},u)$ is the function given in Lemma \ref{Lemma-a}.
\end{lemma}

\proof
First, we assume that $l=2k$ and 
\[
F_{l}(z,\overline{z},u)=\sum_{s=2}^{k}F_{ss}(z,\overline{z},u). 
\]
Suppose that $a\neq 0$ and $F_{kk}(z,\overline{z},0)$ is divided by $\langle
z,z\rangle ^{m}$ for an integer $0\leq m\leq k.$ Then, by Lemma \ref{divided}%
, there are polynomials 
\[
G_{k-m,k-m}^{s}(z,\overline{z},0),\quad s=\max (k-m,2),\cdots ,k, 
\]
of type $(k-m,k-m)$ satisfying 
\[
\frac{F_{ss}(z,\overline{z},u)}{u^{k-s}}=i^{s-k}\langle z,z\rangle
^{m-k+s}G_{k-m,k-m}^{s}(z,\overline{z},0), 
\]
for 
\[
\max (k-m,2)\leq s\leq k. 
\]
Then from the equality \ref{a.20} we obtain 
\begin{eqnarray*}
&&\sum_{s=\max (k-m,2)}^{k}(m-k+s+1)\langle a,z\rangle \langle z,z\rangle
^{m}G_{k-m,k-m}^{s}(z,\overline{z},0) \\
&=&-\sum_{s=\max (k-m,2)}^{k}\langle z,z\rangle ^{m+1}\sum_{\alpha
}a^{\alpha }\left( \frac{\partial G_{k-m,k-m}^{s}}{\partial z^{\alpha }}%
\right) (z,\overline{z},0) \\
&&-\sum_{2\leq s\leq k-m-1}i^{k-s}\langle z,z\rangle ^{k-s}\sum_{\alpha
}a^{\alpha }\frac{\partial }{\partial z^{\alpha }}\{\langle z,z\rangle
F_{ss}(z,\overline{z},0)\} \\
&&-(2i)^{k-1}\langle \kappa ,z\rangle \langle z,z\rangle ^{k}.
\end{eqnarray*}
Hence there are polynomials $A(z,\overline{z};m),$ $1\leq m\leq k,$ such
that 
\begin{equation}
\sum_{s=\max (k-m,2)}^{k}(m-k+s+1)G_{k-m,k-m}^{s}(z,\overline{z},0)=\langle
z,z\rangle A(z,\overline{z};m).  \tag*{(2.14)}  \label{C1}
\end{equation}
The polynomial $A(z,\overline{z};m)$ for $m=k$ is given by 
\begin{align}
A(z,\overline{z};k)& =-(2i)^{k-1}e\langle z,z\rangle ^{k}  \nonumber \\
\langle \kappa ,z\rangle & =e\langle a,z\rangle  \tag*{(2.15)}  \label{D1}
\end{align}
for some constant $e.$ From the equality \ref{6.12}, we obtain for $s=\max
(k-m,2)+1,\cdots ,k,$%
\begin{eqnarray*}
&&(k-s+1)\langle a,z\rangle \langle z,z\rangle ^{m-k+s}G_{k-m,k-m}^{s-1}(z,%
\overline{z},0) \\
&&+(4+2m-3k+3s)\langle a,z\rangle \langle z,z\rangle
^{m-k+s}G_{k-m,k-m}^{s}(z,\overline{z},0) \\
&&+2(m-k+s+1)\langle a,z\rangle \langle z,z\rangle
^{m-k+s}G_{k-m,k-m}^{s+1}(z,\overline{z},0) \\
&=&-2\langle z,z\rangle ^{m-k+s+1}\sum_{\alpha }a^{\alpha }\left( \frac{%
\partial G_{k-m,k-m}^{s}}{\partial z^{\alpha }}\right) (z,\overline{z},0) \\
&&-2\langle z,z\rangle ^{m-k+s+1}\sum_{\alpha }a^{\alpha }\left( \frac{%
\partial G_{k-m,k-m}^{s+1}}{\partial z^{\alpha }}\right) (z,\overline{z},0)
\\
&&+2\langle \kappa ,z\rangle \{1+(-1)^{s}(2s-1)\}\binom{k}{s}(i\langle
z,z\rangle )^{s}.
\end{eqnarray*}
Thus there are polynomials $B^{s-1}(z,\overline{z};m),$ $s=\min
(k-m,2)+1,\cdots ,k,$ such that 
\begin{align}
& (k-s+1)G_{k-m,k-m}^{s-1}(z,\overline{z},0)  \nonumber \\
& +(4-3k+3s+2m)G_{k-m,k-m}^{s}(z,\overline{z},0)  \nonumber \\
& +2(1-k+s+m)G_{k-m,k-m}^{s+1}(z,\overline{z},0)  \nonumber \\
& =\langle z,z\rangle B^{s-1}(z,\overline{z};m).  \tag*{(2.16)}  \label{C2}
\end{align}
The polynomial $B^{s-1}(z,\overline{z};m)$ for $m=k$ is given by 
\begin{align}
B^{s-1}(z,\overline{z};k)& =2i^{s}e\{1+(-1)^{s}(2s-1)\}\binom{k}{s}\langle
z,z\rangle ^{k-m-1},  \nonumber \\
\langle \kappa ,z\rangle & =e\langle a,z\rangle .  \tag*{(2.17)}  \label{D2}
\end{align}
Hence from the equalities \ref{C1} and \ref{C2} we obtain for $k-m\geq 2$%
\begin{equation}
B_{m}\left( 
\begin{array}{l}
G_{k-m,k-m}^{k-m}(z,\overline{z},0) \\ 
G_{k-m,k-m}^{k-m+1}(z,\overline{z},0) \\ 
\vdots \\ 
\\ 
G_{k-m,k-m}^{k}(z,\overline{z},0)
\end{array}
\right) =\langle z,z\rangle \left( 
\begin{array}{l}
A(z,\overline{z};m) \\ 
B^{k-m}(z,\overline{z};m) \\ 
B^{k-m+1}(z,\overline{z};m) \\ 
\vdots \\ 
B^{k-1}(z,\overline{z};m)
\end{array}
\right)  \tag*{(2.18)}  \label{qqq}
\end{equation}
where 
\[
B_{m}=\left( 
\begin{array}{cccccc}
1 & 2 & 3 & \cdots & m & m+1 \\ 
m & 7-m & 4 & 0 & \cdots & 0 \\ 
0 & m-1 & 10-m & 6 & \ddots & \vdots \\ 
& \ddots & \ddots & \ddots & \ddots & 0 \\ 
\vdots &  & \ddots & 2 & 2m+1 & 2m \\ 
0 & \cdots &  & 0 & 1 & 2m+4
\end{array}
\right) . 
\]
By Lemma \ref{Corr1}, the equality \ref{qqq} implies that the function $%
G_{k-m,k-m}^{k}(z,\overline{z},0)$ is divided by $\langle z,z\rangle $ for
all $m\leq k-2.$ Hence the polynomial $F_{kk}(z,\overline{z},0)$ is divided
by $\langle z,z\rangle ^{k-1}$ whenever $a\neq 0.$

Thus $F_{22}(z,\overline{z},u)$ is divided by $\langle z,z\rangle .$ Then
the condition $\Delta F_{22}=0$ implies 
\[
F_{22}(z,\overline{z},u)=i^{2-k}u^{k-2}\langle z,z\rangle G_{11}^{2}(z,%
\overline{z},0)=0. 
\]
Then the equalities \ref{C1} and \ref{C2} yield 
\begin{equation}
B_{k-1}(2)\left( 
\begin{array}{l}
0 \\ 
G_{11}^{3}(z,\overline{z},0) \\ 
\vdots \\ 
\\ 
G_{11}^{k}(z,\overline{z},0)
\end{array}
\right) =\langle z,z\rangle \left( 
\begin{array}{l}
d_{1} \\ 
d_{2} \\ 
d_{3} \\ 
\vdots \\ 
d_{k-1}
\end{array}
\right)  \tag*{(2.19)}  \label{kookoo}
\end{equation}
where $d_{1},\cdots ,d_{k-1}$ are constants and 
\[
B_{k-1}(2)=\left( 
\begin{array}{cccccc}
2 & 3 & 4 & \cdots & k-1 & k \\ 
k-2 & 11-k & 6 & 0 & \cdots & 0 \\ 
0 & k-3 & 14-k & 8 & \ddots & \vdots \\ 
& \ddots & \ddots & \ddots & \ddots & 0 \\ 
\vdots &  & \ddots & 2 & 2k-1 & 2(k-1) \\ 
0 & \cdots &  & 0 & 1 & 2k+2
\end{array}
\right) . 
\]
By Lemma \ref{Lemm2} and Lemma \ref{nonsingular}, the equality \ref{kookoo}
implies that the function $G_{11}^{k}(z,\overline{z},0)$ is divided by $%
\langle z,z\rangle .$ Hence the polynomial $F_{kk}(z,\overline{z},0)$ is
divided by $\langle z,z\rangle ^{k}$ whenever $a\neq 0.$

Thus we obtain 
\begin{eqnarray*}
F_{22}(z,\overline{z},u) &=&0, \\
F_{ss}(z,\overline{z},u) &=&c_{s}\langle z,z\rangle ^{s}\quad \text{for all }%
s=3,\cdots ,k
\end{eqnarray*}
where $c_{s}$ are constant real numbers. By the way, by Lemma \ref{Lemma-a},
the constant $\kappa $ is given by 
\begin{equation}
\langle \kappa ,z\rangle =\frac{u^{2-k}}{4k(k-1)(n+1)(n+2)}\sum_{\alpha
}ua^{\alpha }\Delta ^{2}\left( \frac{\partial F_{33}}{\partial z^{\alpha }}%
\right) (z,\overline{z},u).  \nonumber
\end{equation}
Because of the condition $\Delta ^{3}F_{33}=0,$ we obtain 
\[
F_{33}(z,\overline{z},u)=0\quad \text{and}\quad \kappa =0 
\]
whenever $F_{33}(z,\overline{z},u)$ is divided by $\langle z,z\rangle ^{3}$.
Therefore, we have 
\[
c_{3}=\kappa =0. 
\]
Thus the equalities \ref{D1} and \ref{D2} yield 
\[
A(z,\overline{z};k)=B^{s-1}(z,\overline{z};k)=0 
\]
for all $s=3,\cdots ,k.$ Then the equalities \ref{C1} and \ref{C2} yield 
\[
\left( 
\begin{array}{cccccc}
3 & 4 & 5 & \cdots & k & k+1 \\ 
k-2 & 13-k & 8 & 0 & \cdots & 0 \\ 
0 & k-3 & 16-k & 10 & \ddots & \vdots \\ 
& \ddots & \ddots & \ddots & \ddots & 0 \\ 
\vdots &  & \ddots & 2 & 2k+1 & 2k \\ 
0 & \cdots &  & 0 & 1 & 2k+4
\end{array}
\right) \left( 
\begin{array}{l}
0 \\ 
0 \\ 
c_{4} \\ 
\vdots \\ 
c_{k-1} \\ 
c_{k}
\end{array}
\right) =0. 
\]
Hence we obtain 
\[
c_{4}=\cdots =c_{k}=0. 
\]
This is a contradiction to the assumption 
\[
F_{l}(z,\overline{z},u)\neq 0. 
\]
Thus we ought to have $a=0.$

Assume that $F_{l}(z,\overline{z},u)$ contains a function $F_{st}(z,%
\overline{z},u)$ of type $(s,t),s\neq t,$ so that 
\begin{eqnarray}
l-2p &=&\max \left\{ \left| t-s\right| :F_{st}(z,\overline{z},u)\neq
0\right\}  \nonumber \\
2 &\leq &p\leq \left[ \frac{l-1}{2}\right] ,  \nonumber
\end{eqnarray}
where 
\[
F_{l}(z,\overline{z},u)=\sum_{\min (s,t)\geq 2}F_{st}(z,\overline{z},u). 
\]
Suppose that $p=2.$ Then the equalities \ref{4.58} and \ref{4.55} reduce to 
\[
4i\langle a,z\rangle F_{2,l-2}(z,\overline{z},u)+4i\langle z,z\rangle
\sum_{\alpha }a^{\alpha }\left( \frac{\partial F_{2,l-2}}{\partial z^{\alpha
}}\right) (z,\overline{z},u)=0, 
\]
where 
\[
F_{2,l-2}(z,\overline{z},u)\neq 0. 
\]
Hence we obtain 
\[
\sum_{\alpha }a^{\alpha }\frac{\partial }{\partial z^{\alpha }}\left\{
\langle z,z\rangle F_{2,l-2}(z,\overline{z},u)\right\} =0. 
\]
By Lemma \ref{li-in}, we obtain $a=0.$

Suppose that 
\[
3\leq p\leq \left[ \frac{l-1}{2}\right] . 
\]
and 
\[
F_{p,l-p}(z,\overline{z},u)=0. 
\]
Then by the equalities \ref{4.58} and \ref{4.55}, there is a integer $m$
such that 
\[
\langle z,z\rangle \langle a,z\rangle \left( \frac{\partial F_{m-1,l-m-1}}{%
\partial u}\right) (z,\overline{z},u)=0, 
\]
where 
\begin{gather*}
3\leq m\leq p, \\
F_{m-1,l-m-1}(z,\overline{z},u)\neq 0.
\end{gather*}
Note that 
\begin{eqnarray*}
\left( \frac{\partial F_{m-1,l-m-1}}{\partial u}\right) (z,\overline{z},u)
&=&(p-m+1)u^{-1}F_{m-1,l-m-1}(z,\overline{z},u) \\
&\neq &0.
\end{eqnarray*}
Thus we obtain $a=0.$

Hence we may assume that 
\[
3\leq p\leq \left[ \frac{l-1}{2}\right] . 
\]
and 
\[
F_{p,l-p}(z,\overline{z},u)\neq 0. 
\]
We claim that $F_{p,l-p}(z,\overline{z},0)$ is divided by $\langle
z,z\rangle ^{p-1}$ whenever $a\neq 0.$ Suppose that $a\neq 0$ and $%
F_{p,l-p}(z,\overline{z},0)$ is divided by $\langle z,z\rangle ^{m}$ for an
integer $m,$ $0\leq m\leq p-2.$ Then, by Lemma \ref{divided2}, there are
polynomials 
\begin{equation}
G_{p-m,l-p-m}^{s}(z,\overline{z},0),\quad s=\max (p-m,2),\cdots ,p, 
\tag*{(2.20)}  \label{g-fun}
\end{equation}
of type $(p-m,l-p-m)$ satisfying 
\[
\frac{F_{s,l-2p+s}(z,\overline{z},u)}{u^{p-s}}=i^{s-p}\langle z,z\rangle
^{m-p+s}G_{p-m,l-p-m}^{s}(z,\overline{z},0), 
\]
for 
\[
\max (p-m,2)\leq s\leq p. 
\]
With the polynomials $G_{p-m,l-p-m}^{s}(z,\overline{z},0)$ in \ref{g-fun},
the equality \ref{a.19} yields 
\begin{eqnarray*}
&&\sum_{s=\max (p-m,2)}^{p}(m-p+s+1)\langle a,z\rangle \langle z,z\rangle
^{m}G_{p-m,l-p-m}^{s}(z,\overline{z},0) \\
&=&-\sum_{s=p-m}^{p}\langle z,z\rangle ^{m+1}\sum_{\alpha }a^{\alpha }\left( 
\frac{\partial G_{p-m,l-p-m}^{s}}{\partial z^{\alpha }}\right) (z,\overline{z%
},0) \\
&&-\sum_{2\leq s\leq p-m-1}i^{p-s}\langle z,z\rangle ^{p-s}\sum_{\alpha
}a^{\alpha }\frac{\partial }{\partial z^{\alpha }}\{\langle z,z\rangle
F_{s,l-2p+s}(z,\overline{z},0)\}.
\end{eqnarray*}
Thus there are polynomials $A(z,\overline{z};m),$ $0\leq m\leq p-2,$ such
that 
\begin{equation}
\sum_{s=p-m}^{p}(m-p+s+1)G_{p-m,l-p-m}^{s}(z,\overline{z},0)=\langle
z,z\rangle A(z,\overline{z};m).  \tag*{(2.21)}  \label{B2}
\end{equation}
From the equality \ref{4.58}, we obtain for $s=p-m+1,\cdots ,p,$%
\begin{eqnarray*}
&&(p-s+1)\langle a,z\rangle \langle z,z\rangle ^{m-p+s}G_{p-m,l-p-m}^{s-1}(z,%
\overline{z},0) \\
&&+(4+2m-3p+3s)\langle a,z\rangle \langle z,z\rangle
^{m-p+s}G_{p-m,l-p-m}^{s}(z,\overline{z},0) \\
&&+2(m-p+s+1)\langle a,z\rangle \langle z,z\rangle
^{m-p+s}G_{p-m,l-p-m}^{s+1}(z,\overline{z},0) \\
&=&-2\langle z,z\rangle ^{m-p+s+1}\sum_{\alpha }a^{\alpha }\left( \frac{%
\partial G_{p-m,l-p-m}^{s}}{\partial z^{\alpha }}\right) (z,\overline{z},0)
\\
&&-2\langle z,z\rangle ^{m-p+s+1}\sum_{\alpha }a^{\alpha }\left( \frac{%
\partial G_{p-m,l-p-m}^{s+1}}{\partial z^{\alpha }}\right) (z,\overline{z}%
,0).
\end{eqnarray*}
Thus there are polynomials $B^{s-1}(z,\overline{z};m),$ $s=p-m+1,\cdots ,p,$
such that 
\begin{align}
& (p-s+1)G_{p-m,l-p-m}^{s-1}(z,\overline{z},0)  \nonumber \\
& +(4-3p+3s+2m)G_{p-m,l-p-m}^{s}(z,\overline{z},0)  \nonumber \\
& +2(1-p+s+m)G_{p-m,l-p-m}^{s+1}(z,\overline{z},0)  \nonumber \\
& =\langle z,z\rangle B^{s-1}(z,\overline{z};m).  \tag*{(2.22)}  \label{B3}
\end{align}
Hence, from the equalities \ref{B2} and \ref{B3}, we obtain 
\begin{equation}
B_{m}\left( 
\begin{array}{l}
G_{p-m,l-p-m}^{p-m}(z,\overline{z},0) \\ 
G_{p-m,l-p-m}^{p-m+1}(z,\overline{z},0) \\ 
G_{p-m,l-p-m}^{p-m+2}(z,\overline{z},0) \\ 
\vdots \\ 
G_{p-m,l-p-m}^{p}(z,\overline{z},0)
\end{array}
\right) =\langle z,z\rangle \left( 
\begin{array}{l}
A(z,\overline{z};m) \\ 
B^{p-m}(z,\overline{z};mu) \\ 
B^{p-m+1}(z,\overline{z};m) \\ 
\vdots \\ 
B^{p-1}(z,\overline{z};m)
\end{array}
\right)  \tag*{(2.23)}  \label{eed}
\end{equation}
where 
\[
B_{m}=\left( 
\begin{array}{cccccc}
1 & 2 & 3 & \cdots & m & m+1 \\ 
m & 7-m & 4 & 0 & \cdots & 0 \\ 
0 & m-1 & 10-m & 6 & \ddots & \vdots \\ 
& \ddots & \ddots & \ddots & \ddots & 0 \\ 
\vdots &  & \ddots & 2 & 2m+1 & 2m \\ 
0 & \cdots &  & 0 & 1 & 2m+4
\end{array}
\right) . 
\]
By Lemma \ref{Corr1}, the equality \ref{eed} implies that the function $%
G_{p-m,l-p-m}^{p}(z,\overline{z},0)$ is divided by $\langle z,z\rangle .$
Hence we prove our claim that $F_{p,l-p}(z,\overline{z},0)$ is divided by $%
\langle z,z\rangle ^{p-1}$ whenever $a\neq 0.$

Then we claim that $F_{p,l-p}(z,\overline{z},0)$ is divided by $\langle
z,z\rangle ^{p}$ whenever $a\neq 0.$ With the polynomials $%
G_{1,l-2p+1}^{s}(z,\overline{z},0)$ in \ref{g-fun}, the equality \ref{a.19}
yields 
\[
\langle a,z\rangle \sum_{s=2}^{p}sG_{1,l-2p+1}^{s}(z,\overline{z}%
,0)=-\langle z,z\rangle \sum_{s=2}^{p}\sum_{\alpha }a^{\alpha }\left( \frac{%
\partial G_{1,l-2p+1}^{s}}{\partial z^{\alpha }}\right) (z,\overline{z},0). 
\]
So there is a polynomial $A(z,\overline{z};p-1)$ of type $(0,l-2p)$ such
that 
\begin{equation}
\sum_{s=2}^{p}sG_{1,l-2p+1}^{s}(z,\overline{z},0)=\langle z,z\rangle A(z,%
\overline{z};p-1).  \tag*{(2.24)}  \label{edd}
\end{equation}
With the polynomials $G_{1,l-2p+1}^{s}(z,\overline{z},0)$ in \ref{g-fun},
the equality \ref{4.58} yields 
\begin{eqnarray*}
&&\langle a,z\rangle \left\{ (p-s+1)G_{1,l-2p+1}^{s-1}(z,\overline{z}%
,0)+(2-p+3s)G_{1,l-2p+1}^{s}(z,\overline{z},0)\right. \\
&&\hspace{3cm}\left. +2sG_{1,l-2p+1}^{s+1}(z,\overline{z},0)\right\} \\
&=&-2\langle z,z\rangle \left\{ \sum_{\alpha }a^{\alpha }\left( \frac{%
\partial G_{1,l-2p+1}^{s}}{\partial z^{\alpha }}\right) (z,\overline{z}%
,0)+\sum_{\alpha }a^{\alpha }\left( \frac{\partial G_{1,l-2p+1}^{s+1}}{%
\partial z^{\alpha }}\right) (z,\overline{z},0)\right\} .
\end{eqnarray*}
Then there are polynomials $B^{s-1}(z,\overline{z};p-1)$ of type $(0,l-2p)$
for $s=3,\cdots ,p$ such that 
\begin{align}
(p-s+1)G_{1,l-2p+1}^{s-1}(z,\overline{z},0)+& (2-p+3s)G_{1,l-2p+1}^{s}(z,%
\overline{z},0)  \nonumber \\
+2sG_{1,l-2p+1}^{s+1}(z,\overline{z},0)& =\langle z,z\rangle B^{s-1}(z,%
\overline{z};p-1).  \tag*{(2.25)}  \label{edd2}
\end{align}
Hence, from the equalities \ref{edd} and \ref{edd2}, we obtain 
\begin{equation}
B_{p-1}(2)\left( 
\begin{array}{l}
G_{1,l-2p+1}^{2}(z,\overline{z},0) \\ 
G_{1,l-2p+1}^{3}(z,\overline{z},0) \\ 
\vdots \\ 
\\ 
G_{1,l-2p+1}^{p}(z,\overline{z},0)
\end{array}
\right) =\langle z,z\rangle \left( 
\begin{array}{l}
A(z,\overline{z};p-1) \\ 
B^{2}(z,\overline{z};p-1) \\ 
B^{3}(z,\overline{z};p-1) \\ 
\vdots \\ 
B^{p-1}(z,\overline{z};p-1)
\end{array}
\right)  \tag*{(2.26)}  \label{ewr}
\end{equation}
where 
\[
B_{p-1}(2)=\left( 
\begin{array}{cccccc}
2 & 3 & 4 & \cdots & p-1 & p \\ 
p-2 & 11-p & 6 & 0 & \cdots & 0 \\ 
0 & p-3 & 14-p & 8 & \ddots & \vdots \\ 
& \ddots & \ddots & \ddots & \ddots & 0 \\ 
\vdots &  & \ddots & 2 & 2p-1 & 2(p-1) \\ 
0 & \cdots &  & 0 & 1 & 2p+2
\end{array}
\right) . 
\]
By Lemma \ref{Lemm2} and Lemma \ref{nonsingular}, the equality \ref{ewr}
implies that the polynomial $G_{1,l-2p+1}^{p}(z,\overline{z},0)$ is divided
by $\langle z,z\rangle .$ Hence we prove our claim that $F_{p,l-p}(z,%
\overline{z},0)$ is divided by $\langle z,z\rangle ^{p}$ whenever $a\neq 0.$

Then with the polynomials $G_{0,l-2p}^{s}(z,\overline{z},0)$ in \ref{g-fun},
the equality \ref{a.19} yields 
\[
\langle a,z\rangle \sum_{s=2}^{p}(s+1)G_{0,l-2p}^{s}(z,\overline{z},0)=0. 
\]
Whenever $a\neq 0,$ we have 
\begin{equation}
\sum_{s=2}^{p}(s+1)G_{1,l-2p+1}^{s}(z,\overline{z},0)=0.  \tag*{(2.27)}
\label{end}
\end{equation}
With the polynomials $G_{0,l-2p}^{s}(z,\overline{z},0)$ in \ref{g-fun}, the
equality \ref{4.58} yields 
\begin{align*}
\langle a,z\rangle & \left\{ (p-s+1)G_{0,l-2p}^{s-1}(z,\overline{z}%
,0)+(4-p+3s)G_{0,l-2p}^{s}(z,\overline{z},0)+\right. \\
& \hspace{5.2cm}\left. 2(s+1)G_{0,l-2p}^{s+1}(z,\overline{z},0)\right\} =0.
\end{align*}
Whenever $a\neq 0,$ we have 
\begin{align}
(p-s+1)& G_{0,l-2p}^{s-1}(z,\overline{z},0)+(4-p+3s)G_{0,l-2p}^{s}(z,%
\overline{z},0)+  \nonumber \\
& \hspace{3cm}2(s+1)G_{0,l-2p}^{s+1}(z,\overline{z},0)=0  \tag*{(2.28)}
\label{end2}
\end{align}
for $s=3,\cdots ,p.$ Hence, from the equalities \ref{end} and \ref{end2}, we
obtain 
\begin{equation}
B_{p}(3)\left( 
\begin{array}{l}
G_{0,l-2p}^{2}(z,\overline{z},0) \\ 
G_{0,l-2p}^{3}(z,\overline{z},0) \\ 
\vdots \\ 
G_{0,l-2p}^{p}(z,\overline{z},0)
\end{array}
\right) =0  \tag*{(2.29)}  \label{eww}
\end{equation}
where 
\[
B_{p}(3)=\left( 
\begin{array}{cccccc}
3 & 4 & 5 & \cdots & p & p+1 \\ 
p-2 & 13-p & 8 & 0 & \cdots & 0 \\ 
0 & p-3 & 16-p & 10 & \ddots & \vdots \\ 
& \ddots & \ddots & \ddots & \ddots & 0 \\ 
\vdots &  & \ddots & 2 & 2p+1 & 2p \\ 
0 & \cdots &  & 0 & 1 & 2p+4
\end{array}
\right) . 
\]
By Lemma \ref{Lemm2} and Lemma \ref{nonsingular}, the equality \ref{eww}
implies 
\[
G_{0,l-2p}^{p}(z,\overline{z},0)=0. 
\]
This is a contradiction to the assumption 
\[
F_{p,l-p}(z,\overline{z},u)=\langle z,z\rangle ^{p}G_{0,l-2p}^{p}(z,%
\overline{z},0)\neq 0. 
\]
Thus we ought to have $a=0$ as well for the case of $3\leq p\leq \left[ 
\frac{l-1}{2}\right] .$ Therefore we obtain $a=0$ whenever $F_{l}(z,%
\overline{z},u)$ contains a nonvanishing term $F_{st}(z,\overline{z},u)$ of
type $(s,t),$ $s\neq t.$

Therefore, we have showed that $a=0$ whenever 
\[
F_{l}(z,\overline{z},u)\neq 0\quad \text{and}\quad H_{l+1}(z,\overline{z}%
,u;a)=0. 
\]
This completes the proof.\endproof

\begin{theorem}
Let $M$ be a real hypersurface in normal form defined by the equation 
\[
v=\langle z,z\rangle +F_{l}(z,\overline{z},u)+O(l+2),
\]
where 
\[
F_{l}(z,\overline{z},u)\neq 0
\]
and, for all complex numbers $\mu ,$%
\[
F_{l}(\mu z,\mu \overline{z},\mu ^{2}u)=\mu ^{l}F_{l}(z,\overline{z},u).
\]
Suppose that there is a normalization $\phi $ of $M$ with initial value $%
(id_{n\times n},a,1,0)\in H$ such that $\phi $ transforms $M$ to a real
hypersurface $M^{\prime }$ defined by the equation 
\[
v=\langle z,z\rangle +F^{*}(z,\overline{z},u)
\]
and 
\[
F^{*}(z,\overline{z},u)=F_{l}(z,\overline{z},u)+O(l+2).
\]
Then the normalization $\phi $ has identity initial value, i.e., $a=0.$
\end{theorem}

\proof
The conclusion follows from Lemma \ref{Lemma-a} and \ref{Theo1}.\endproof

\subsection{Beloshapka-Loboda Theorem}

\begin{lemma}
\label{orbit}Let $M$ be a real hypersurface in normal form and $\phi
_{\sigma _{1}}$ be a normalization of $M$ with initial value $\sigma _{1}\in
H.$ Suppose that $M$ is transformed to $M^{\prime }$ by the normalization $%
\phi _{\sigma _{1}}$ and $\phi _{\sigma _{2}}$ is a normalization of $%
M^{\prime }$ with initial value $\sigma _{2}\in H$. Then 
\[
\phi _{\sigma _{1}}\circ \phi _{\sigma _{2}}=\phi _{\sigma _{1}\sigma _{2}}
\]
where $\phi _{\sigma _{1}\sigma _{2}}$ is a normalization of $M$ with
initial value $\sigma _{1}\sigma _{2}\in H$.
\end{lemma}

In the paper \cite{Pa2}, we have given the proof of Lemma \ref{orbit}.

\begin{lemma}
Let $M$ be a real hypersurface in normal form defined by the equation 
\[
v=\langle z,z\rangle +F_{l}(z,\overline{z},u)+O(l+2),
\]
where 
\[
F_{l}(z,\overline{z},u)\neq 0
\]
and, for all complex numbers $\mu ,$%
\[
F_{l}(\mu z,\mu \overline{z},\mu ^{2}u)=\mu ^{l}F_{l}(z,\overline{z},u).
\]
Suppose that there is a normalization $\phi $ of $M$ such that $\phi (M)$ is
defined by the equation 
\[
v=\langle z,z\rangle +\rho F_{l}(C^{-1}z,\overline{C^{-1}z},\rho
^{-1}u)+O(l+2)
\]
where 
\[
\sigma \left( \phi \right) =(C,a,\rho ,r)\in H.
\]
Then the normalization $\phi $ have the initial value $(C,0,\rho ,r)\in H$,
i.e., $a=0.$
\end{lemma}

\proof
Note that there is a decomposition of $\phi $ as follows(cf. \cite{Pa1}): 
\[
\phi =\phi _{\sigma _{1}}\circ \phi _{\sigma _{2}} 
\]
where $\phi _{\sigma _{1}},\phi _{\sigma _{2}}$ are normalizations with the
initial values $\sigma _{1},\sigma _{2}$ respectively: 
\begin{eqnarray*}
\sigma _{1} &=&(C,0,\rho ,r)\in H, \\
\sigma _{2} &=&(id_{n\times n},a,1,0)\in H.
\end{eqnarray*}
Then, by Lemma \ref{orbit}, we obtain 
\begin{eqnarray*}
\phi _{\sigma _{2}} &=&\phi _{\sigma _{1}^{-1}}\circ \phi _{\sigma
_{1}}\circ \phi _{\sigma _{2}} \\
&=&\phi _{\sigma _{1}^{-1}}\circ \phi
\end{eqnarray*}
where $\phi _{\sigma _{1}^{-1}}$ is a normalization with initial value $%
\sigma _{1}^{-1}\in H$. Further, suppose that $\phi _{\sigma _{2}}(M)$ is
defined by the equation 
\[
v=\langle z,z\rangle +F^{*}(z,\overline{z},u). 
\]
Then we obtain 
\[
F^{*}(z,\overline{z},u)=F_{l}(z,\overline{z},u)+O(l+2). 
\]
Thus, by Lemma \ref{Theo1}, we obtain 
\[
a=0. 
\]
This completes the proof.\endproof

\begin{theorem}[Beloshapka, Loboda, Vitushkin]
\label{ThBL}Let $M$ be an analytic real hypersurface in normal form, which
is not a real hyperquadric, and $H(M)$ be the isotropy subgroup of $M$ at
the origin. Then there are functions 
\[
\rho (U),\quad a(U),\quad r(U)
\]
on the set 
\begin{equation}
\left\{ U:\left( U,a,\rho ,r\right) \in H(M)\subset H\right\}   \nonumber
\end{equation}
such that, for all $(U,a,\rho ,r)\in H(M),$%
\[
a=a(U),\quad \rho =\rho (U),\quad r=r(U).
\]
\end{theorem}

\proof
Suppose that $M$ is defined in normal form by the equation 
\[
v=\langle z,z\rangle +F_{l}(z,\overline{z},u)+F_{l+1}(z,\overline{z}%
,u)+F_{l+2}(z,\overline{z},u)+O(l+3), 
\]
where 
\[
F_{l}(z,\overline{z},u)\neq 0, 
\]
and the integers $l,l+1,l+2$ represent the weight of the functions 
\[
F_{l}(z,\overline{z},u),\quad F_{l+1}(z,\overline{z},u),\quad F_{l+2}(z,%
\overline{z},u). 
\]

Let $\phi _{\sigma }$ be a normalization of $M$ with initial value $\sigma
\in H(M)$. Suppose that the real hypersurface $\phi _{\sigma }(M)$ is
defined near the origin up to weight $l$ by the equation 
\begin{eqnarray*}
v &=&\langle z,z\rangle +\rho F_{l}(C^{-1}z,\overline{C^{-1}z},\rho
^{-1}u)+O(l+1) \\
&=&\langle z,z\rangle +F_{l}(z,\overline{z},u)+O(l+1)
\end{eqnarray*}
where 
\[
\sigma =(C,a,\rho ,r)\in H(M)\subset H. 
\]
Then we have 
\begin{equation}
\left| \rho \right| ^{\frac{l-2}{2}}F_{l}(z,\overline{z},u)=\lambda
F_{l}(U^{-1}z,\overline{U^{-1}z},\lambda u)\neq 0.  \tag*{(2.30)}
\label{rho-abs}
\end{equation}
The relation 
\[
\langle Uz,Uz\rangle =\lambda \langle z,z\rangle ,\quad \lambda =\mathrm{%
sign\{}\rho \mathrm{\}} 
\]
yields 
\begin{equation}
\lambda =\frac{1}{n}\Delta \langle Uz,Uz\rangle =\pm 1.  \tag*{(2.31)}
\label{rho-sig}
\end{equation}
Then we take a value $z,u$ in the equality \ref{rho-abs} such that 
\[
F_{l}(z,\overline{z},u)\in \Bbb{R}\backslash \{0\} 
\]
and define 
\[
\rho _{1}(U)=\left( \frac{\lambda F_{l}(U^{-1}z,\overline{U^{-1}z},\lambda u)%
}{F_{l}(z,\overline{z},u)}\right) ^{\frac{2}{l-2}}. 
\]
By the unique factorization of a polynomial, we have 
\[
\left| \rho \right| =\rho _{1}(U) 
\]
regardless the choice of the value $z,u.$ Hence, by the equality \ref
{rho-sig}, we define 
\[
\rho (U)\equiv \frac{1}{n}\Delta \langle Uz,Uz\rangle \cdot \rho _{1}(U) 
\]
so that 
\begin{equation}
\rho =\rho (U)  \tag*{(2.32)}  \label{rho}
\end{equation}
for all 
\[
\left( U,a,\rho ,r\right) \in H(M). 
\]

Suppose that the real hypersurface $\phi _{\sigma }(M)$ is defined near the
origin up to weight $l+1$ by the equation 
\begin{eqnarray*}
v-\langle z,z\rangle &=&\rho F_{l}(C^{-1}z,\overline{C^{-1}z},\rho
^{-1}u)+F_{l+1}^{*}(z,\overline{z},u)+O(l+2) \\
&=&F_{l}(z,\overline{z},u)+F_{l+1}(z,\overline{z},u)+O(l+2).
\end{eqnarray*}
By using the equality 
\[
\rho F_{l}(C^{-1}z,\overline{C^{-1}z},\rho ^{-1}u)=F_{l}(z,\overline{z},u), 
\]
we obtain 
\[
F_{l+1}^{*}(z,\overline{z},u)=H_{l+1}(z,\overline{z},u;\rho ^{-1}Ca)+\rho
F_{l+1}(C^{-1}z,\overline{C^{-1}z},\rho ^{-1}u) 
\]
where $a^{*}\rightarrow H_{l+1}(z,\overline{z},u;a^{*})$ is the injective
linear mapping in Lemma \ref{Theo1}.

Then the following requirement 
\[
F_{l+1}^{*}(z,\overline{z},u)=F_{l+1}(z,\overline{z},u) 
\]
yields 
\[
H_{l+1}(z,\overline{z},u;a^{*})=F_{l+1}(z,\overline{z},u)-\rho
F_{l+1}(C^{-1}z,\overline{C^{-1}z},\rho ^{-1}u). 
\]
Then, by the equality \ref{rho}, the equality 
\begin{eqnarray*}
H_{l+1}(z,\overline{z},u;a^{*}) &=&F_{l+1}(z,\overline{z},u)-\rho
F_{l+1}(C^{-1}z,\overline{C^{-1}z},\rho ^{-1}u) \\
&=&F_{l+1}(z,\overline{z},u)-\mathrm{sign\{}\rho \left( U\right) \mathrm{\}}%
\left| \rho \left( U\right) \right| ^{\frac{l+3}{2}}F_{l+1}(U^{-1}z,%
\overline{U^{-1}z},\mathrm{sign\{}\rho \left( U\right) \mathrm{\}}u),
\end{eqnarray*}
yields a unique function $a^{*}(U)$ of $U$ satisfying 
\[
a^{*}=\rho ^{-1}Ca=a^{*}(U). 
\]
Hence we obtain a unique function $a(U)$ of $U$ such that 
\begin{align}
a& =a(U)  \nonumber \\
& \equiv \rho (U)\left| \rho (U)\right| ^{-\frac{1}{2}}U^{-1}a^{*}(U) 
\tag*{(2.33)}  \label{a}
\end{align}
for all 
\[
\left( U,a,\rho ,r\right) \in H(M). 
\]

Then we decompose the normalization $\phi _{\sigma }$ as follows: 
\[
\phi _{\sigma }=\phi _{2}\circ \phi _{1}, 
\]
where $\phi _{1},\phi _{2}$ are normalizations with the initial values $%
\sigma _{1},\sigma _{2}$ respectively: 
\[
\sigma _{1}=(id_{n\times n},a,1,0)\quad \text{and\quad }\sigma
_{2}=(C,0,\rho ,r). 
\]
Suppose that the real hypersurface $\phi _{1}(M)$ is defined by the equation 
\[
v=\langle z,z\rangle +F_{l}(z,\overline{z},u)+\tilde{F}_{l+1}(z,\overline{z}%
,u)+\tilde{F}_{l+2}(z,\overline{z},u)+O(l+3) 
\]
where the functions $\tilde{F}_{l+1}(z,\overline{z},u)$ and $\tilde{F}%
_{l+2}(z,\overline{z},u)$ depend of the parameter $a,$ i.e., 
\begin{eqnarray*}
\tilde{F}_{l+1}(z,\overline{z},u) &=&\tilde{F}_{l+1}(z,\overline{z},u;a) \\
\tilde{F}_{l+2}(z,\overline{z},u) &=&\tilde{F}_{l+2}(z,\overline{z},u;a).
\end{eqnarray*}
Then suppose that the real hypersurface $\phi _{\sigma }(M)$ is defined near
the origin up to weight $l+2$ by the equation 
\begin{eqnarray*}
v &=&\langle z,z\rangle +\rho F_{l}(C^{-1}z,\overline{C^{-1}z},\rho
^{-1}u)+\rho \tilde{F}_{l+1}(C^{-1}z,\overline{C^{-1}z},\rho ^{-1}u) \\
&&+\rho \tilde{F}_{l+2}(C^{-1}z,\overline{C^{-1}z},\rho ^{-1}u) \\
&&-\frac{r}{2}\left\{ \sum_{\min (s,t)\geq 2}(l+s+t)uF_{st}(C^{-1}z,%
\overline{C^{-1}z},\rho ^{-1}u)\right. \\
&&\hspace{1.5cm}+\sum_{\min (s,t)\geq 2}2(s-t)i\langle z,z\rangle
F_{st}(C^{-1}z,\overline{C^{-1}z},\rho ^{-1}u) \\
&&\hspace{1.5cm}\left. -\sum_{\min (s,t)\geq 2}2\rho ^{-1}\langle z,z\rangle
^{2}\left( \frac{\partial F_{st}}{\partial u}\right) (C^{-1}z,\overline{%
C^{-1}z},\rho ^{-1}u)\right\} \\
&&+O(l+3) \\
&=&F_{l}(z,\overline{z},u)+F_{l+1}(z,\overline{z},u)+F_{l+2}(z,\overline{z}%
,u)+O(l+3),
\end{eqnarray*}
where 
\[
F_{l}(z,\overline{z},u)=\sum_{\min (s,t)\geq 2}F_{st}(z,\overline{z},u). 
\]
Hence we have the equality 
\begin{align}
& -\frac{r}{2}\left\{ \sum_{\min (s,t)\geq 2}\left( (l+s+t)u+2(s-t)i\langle
z,z\rangle \right) F_{st}(z,\overline{z},u)\right.  \nonumber \\
& \hspace{1.5cm}\left. -\sum_{\min (s,t)\geq 2}2\langle z,z\rangle
^{2}\left( \frac{\partial F_{st}}{\partial u}\right) (z,\overline{z}%
,u)\right\}  \nonumber \\
& =\rho ^{-1}F_{l+2}(Cz,\overline{Cz},\rho u)-\tilde{F}_{l+2}\left( z,%
\overline{z},u;a\right) .  \tag*{(2.34)}  \label{rrr}
\end{align}
Note that $F_{l}(z,\overline{z},u)\neq 0$ implies 
\[
\sum_{\min (s,t)\geq 2}\left\{ \left( (l+s+t)u+2(s-t)i\langle z,z\rangle
\right) F_{st}(z,\overline{z},u)-2\langle z,z\rangle ^{2}\left( \frac{%
\partial F_{st}}{\partial u}\right) (z,\overline{z},u)\right\} \neq 0. 
\]
Otherwise, we would have 
\begin{align*}
& (l+s+t)uF_{st}(z,\overline{z},u)+2(s-t)i\langle z,z\rangle F_{s-1,t-1}(z,%
\overline{z},u) \\
& -2\langle z,z\rangle ^{2}\left( \frac{\partial F_{s-2,t-2}}{\partial u}%
\right) (z,\overline{z},u)=0
\end{align*}
which yields 
\[
F_{st}(z,\overline{z},u)=0\quad \text{for all }s,t. 
\]
From the equalities \ref{rho} and \ref{a}, we have 
\[
\rho =\rho (U),\quad a=a(U) 
\]
for all 
\[
(U,a,\rho ,r)\in H(M). 
\]
Then we take a value $z,u$ in the equality \ref{rrr} such that 
\[
\sum_{\min (s,t)\geq 2}\left\{ \left( (l+s+t)u+2(s-t)i\langle z,z\rangle
\right) F_{st}(z,\overline{z},u)-2\langle z,z\rangle ^{2}\left( \frac{%
\partial F_{st}}{\partial u}\right) (z,\overline{z},u)\right\} \in \Bbb{R}%
\backslash \{0\} 
\]
and define 
\[
r(U)=\frac{-2\left\{ \rho (U)^{-1}\left| \rho (U)\right| ^{\frac{l+2}{2}%
}F_{l+2}(Uz,\overline{Uz},\lambda u)-\tilde{F}_{l+2}\left( z,\overline{z}%
,u;a(U)\right) \right\} }{\sum_{\min (s,t)\geq 2}\left\{ \left(
(l+s+t)u+2(s-t)i\langle z,z\rangle \right) F_{st}(z,\overline{z},u)-2\langle
z,z\rangle ^{2}\left( \frac{\partial F_{st}}{\partial u}\right) (z,\overline{%
z},u)\right\} }. 
\]
By the unique factorization of a polynomial, we have 
\[
r=r(U) 
\]
regardless the choice of the value $z,u.$ Thus the equality \ref{rrr} yields
a unique function $r(U)$ of $U$ satisfying 
\[
r=r(U) 
\]
for all 
\[
\left( U,a,\rho ,r\right) \in H(M). 
\]
This completes the proof.\endproof

\section{Compact local automorphism groups}

\subsection{Compactness}

\begin{lemma}
\label{analytic}Let $M$ be a nondegenerate analytic real hypersurface
defined by 
\[
v=F(z,\overline{z},u),\quad \left. F\right| _{0}=\left. dF\right| _{0}=0, 
\]
and $\phi _{\sigma }$ be a normalization of $M$ with initial value $\sigma
\in H$. Suppose that $\phi _{\sigma }$ transforms $M$ to a real hypersurface
defined by the equation 
\[
v=\langle z,z\rangle +F^{*}(z,\overline{z},u;\sigma ). 
\]
Then the functions $\phi _{\sigma }\left( z,w\right) $ and $F^{*}(z,%
\overline{z},u;\sigma )$ are analytic of 
\[
\sigma =\left( U,a,\rho ,r\right) \in H. 
\]
Further, each coefficient 
\[
\left( \left. \frac{\partial ^{\left| I\right| +l}\phi _{\sigma }}{\partial
z^{\left| I\right| }\partial w^{l}}\right| _{0}\right) \quad \text{and}\quad
\left( \left. \frac{\partial ^{\left| I\right| +\left| J\right| +l}F^{*}}{%
\partial z^{\left| I\right| }\partial \overline{z}^{\left| J\right|
}\partial u^{l}}\right| _{0}\right) 
\]
depends polynomially on the parameters 
\[
C\equiv \sqrt{\left| \rho \right| }U,\quad C^{-1},\quad \rho ,\quad \rho
^{-1},\quad a,\quad r. 
\]
\end{lemma}

In the paper \cite{Pa2}, we have given the proof of Lemma \ref{analytic}.

Let $M$ be a real hypersurface $M$ in normal form. We define the isotropy
subgroup $H(M)$ of $M$ at the origin as follows: 
\[
H(M)=\{\sigma \in H:\phi _{\sigma }(M)=M\} 
\]
where $\phi _{\sigma }$ is a normalization of $M$ with initial value $\sigma
\in H$. By Lemma \ref{analytic}, the group $H$ is homeomorphic to the set of
germs $\phi _{\sigma },$ $\sigma \in H,$ with a topology induced from the
natural compact-open topology. Further, by Lemma \ref{orbit} and Lemma \ref
{analytic}, the group $H(M)$ is isomorphic as Lie group to the local
automorphism group of $M$.

\begin{lemma}
\label{estimation}Let $M$ be a nonspherical analytic real hypersurface and $%
H(M)$ be the isotropy subgroup of $M$ such that there is a real number $%
c\geq 1$ satisfying 
\[
\sup_{\left( U,a,\rho ,r\right) \in H(M)}\left\| U\right\| \leq c<\infty . 
\]
Then there exists a real number $e>0$ satisfying 
\[
\left| a\right| \leq e,\quad e^{-1}\leq \left| \rho \right| \leq e,\quad
\left| r\right| \leq e 
\]
for all elements 
\[
\left( U,a,\rho ,r\right) \in H(M) 
\]
where $e$ may depend on $M$ and $c.$
\end{lemma}

\proof
For the parameter $\rho ,$ we have 
\begin{eqnarray*}
\left| \rho \left( U\right) \right| ^{\frac{l-2}{2}} &=&\frac{\left|
F_{l}(U^{-1}z,\overline{U^{-1}z},\lambda u)\right| }{\left| F_{l}(z,%
\overline{z},u)\right| } \\
\left| \rho \left( U\right) \right| ^{-\frac{l-2}{2}} &=&\frac{\left|
F_{l}(Uz,\overline{Uz},\lambda u)\right| }{\left| F_{l}(z,\overline{z}%
,u)\right| }
\end{eqnarray*}
whenever we take a value $z,u$ satisfying 
\[
F_{l}(z,\overline{z},u)\neq 0. 
\]
Hence we have the following estimate: 
\begin{eqnarray*}
\left| \rho \left( U\right) \right| ^{\frac{l-2}{2}} &\leq &\sup_{U}\frac{%
\left| F_{l}(U^{-1}z,\overline{U^{-1}z},\lambda u)\right| }{\left| F_{l}(z,%
\overline{z},u)\right| } \\
\left| \rho \left( U\right) \right| ^{-\frac{l-2}{2}} &\leq &\sup_{U}\frac{%
\left| F_{l}(Uz,\overline{Uz},\lambda u)\right| }{\left| F_{l}(z,\overline{z}%
,u)\right| }
\end{eqnarray*}
so that 
\[
\left( \sup_{U}\frac{\left| F_{l}(Uz,\overline{Uz},\lambda u)\right| }{%
\left| F_{l}(z,\overline{z},u)\right| }\right) ^{-1}\leq \left| \rho \left(
U\right) \right| ^{\frac{l-2}{2}}\leq \sup_{U}\frac{\left| F_{l}(U^{-1}z,%
\overline{U^{-1}z},\lambda u)\right| }{\left| F_{l}(z,\overline{z},u)\right| 
}. 
\]
Note that there is a real number $d$ depending only on $F_{l}(z,\overline{z}%
,u)$ such that 
\[
\frac{\left| F_{l}(U^{-1}z,\overline{U^{-1}z},\lambda u)\right| }{\left|
F_{l}(z,\overline{z},u)\right| }\leq d_{1}\cdot c^{l} 
\]
where 
\[
c\equiv \sup_{\left( U,a,\rho ,r\right) \in H(M)}\left\| U\right\| \geq 1. 
\]
Thus we obtain 
\[
d_{1}^{-\frac{2}{l-2}}\cdot c^{-\frac{2l}{l-2}}\leq \left| \rho \left(
U\right) \right| \leq d_{1}^{\frac{2}{l-2}}\cdot c^{\frac{2l}{l-2}}. 
\]

For the parameter $a,$ we have 
\begin{eqnarray*}
&&H_{l+1}\left( z,\overline{z},u;a^{*}(U)\right) \\
&=&F_{l+1}(z,\overline{z},u)-\mathrm{sign\{}\rho \left( U\right) \mathrm{\}}%
\left| \rho \left( U\right) \right| ^{\frac{l+3}{2}}F_{l+1}(U^{-1}z,%
\overline{U^{-1}z},\mathrm{sign\{}\rho \left( U\right) \mathrm{\}}u)
\end{eqnarray*}
where 
\[
a^{*}\left( U\right) =\rho \left( U\right) ^{-1}\sqrt{\left| \rho \left(
U\right) \right| }Ua\left( U\right) . 
\]
Since the mapping $a\mapsto H_{l+1}\left( z,\overline{z},u;a\right) $ is
injective and the function $H_{l+1}\left( z,\overline{z},u;a\right) $
depends only on $F_{l}(z,\overline{z},u),$ we have the following estimate: 
\[
\left| a^{*}\left( U\right) \right| \leq d_{2}^{*}\cdot c^{\frac{2(l^{2}+l-1)%
}{l-2}} 
\]
which yields 
\[
\left| a\left( U\right) \right| \leq d_{2}\cdot c^{\frac{2(l^{2}+2l-2)}{l-2}%
} 
\]
where $d_{2}^{*},d_{2}$ depend only on $F_{l}(z,\overline{z},u)$ and $%
F_{l+1}(z,\overline{z},u).$

For the parameter $r,$ we have 
\begin{align}
& -\frac{r\left( U\right) }{2}\left\{ \sum_{\min (s,t)\geq
2}(l+s+t)uF_{st}(z,\overline{z},u)\right.  \nonumber \\
& \hspace{1.5cm}+\sum_{\min (s,t)\geq 2}2(s-t)i\langle z,z\rangle F_{st}(z,%
\overline{z},u)  \nonumber \\
& \hspace{1.5cm}\left. -\sum_{\min (s,t)\geq 2}2\langle z,z\rangle
^{2}\left( \frac{\partial F_{st}}{\partial u}\right) (z,\overline{z}%
,u)\right\}  \nonumber \\
& =\lambda \left| \rho \left( U\right) \right| ^{\frac{l}{2}}F_{l+2}(Uz,%
\overline{Uz},\lambda u)-\tilde{F}_{l+2}\left( z,\overline{z},u;a\right) . 
\nonumber
\end{align}
By Lemma \ref{analytic}, the function $\tilde{F}_{l+2}\left( z,\overline{z}%
,u;a\right) $ depend polynomially, in fact, quadratically, on the parameter $%
a.$ Hence we obtain the following estimate: 
\[
\left| r\left( U\right) \right| \leq d_{3}\cdot c^{L}\quad \text{for some }%
L\in \Bbb{N} 
\]
where $d_{3}$ depends only on $F_{l}\left( z,\overline{z},u\right) ,$ $%
F_{l+1}\left( z,\overline{z},u\right) ,$ and $F_{l+2}\left( z,\overline{z}%
,u\right) .$

Then we take 
\[
e=\max \left\{ d_{1}^{\frac{2}{l-2}}\cdot c^{\frac{2l}{l-2}},\quad
d_{2}\cdot c^{\frac{2(l^{2}+2l-2)}{l-2}},\quad d_{3}\cdot c^{L}\right\} . 
\]
This completes the proof.\endproof

\begin{theorem}
\label{compactness}Let $M$ be a nonspherical analytic real hypersurface in
normal form. Suppose that there is a real number $c\geq 1$ satisfying 
\begin{equation}
\sup_{\left( U,a,\rho ,r\right) \in H(M)}\left\| U\right\| \leq c<\infty . 
\tag*{(3.1)}  \label{fff}
\end{equation}
Then the group $H(M)$ is compact.
\end{theorem}

\proof
By Theorem \ref{ThBL}, the group $H(M)$ is isomorphic to the following
group: 
\begin{equation}
\left\{ U:\left( U,a,\rho ,r\right) \in H(M)\right\} .  \nonumber
\end{equation}
We claim that the group $H(M)$ is closed under the condition \ref{fff}.
Suppose that there is a convergent sequence in $GL(n;\Bbb{C})$ such that 
\[
U_{m}\in \left\{ U:\left( U,a,\rho ,r\right) \in H(M)\right\} \quad \text{%
for all }m\in \Bbb{N} 
\]
and, by the condition \ref{fff}, 
\[
\lim_{m\rightarrow \infty }U_{m}=U\in GL(n;\Bbb{C}). 
\]
Then, by the functions $\rho (U),a(U),r(U)$ in Theorem \ref{ThBL}, we have
the following sequence: 
\[
\left( U_{m},a\left( U_{m}\right) ,\rho \left( U_{m}\right) ,r\left(
U_{m}\right) \right) \in H(M). 
\]
Under the condition \ref{fff}, by Lemma \ref{estimation}, there is a real
number $e>0$ such that 
\[
\left| a\left( U_{m}\right) \right| \leq e,\quad e^{-1}\leq \left| \rho
\left( U_{m}\right) \right| \leq e,\quad \left| r\left( U_{m}\right) \right|
\leq e\quad \text{for all }m. 
\]
Then, by compactness, there is a subsequence $m_{j}$ such that the following
limits exists: 
\begin{eqnarray*}
a &=&\lim_{j\rightarrow \infty }a\left( U_{m_{j}}\right) , \\
\rho &=&\lim_{j\rightarrow \infty }\rho \left( U_{m_{j}}\right) , \\
r &=&\lim_{j\rightarrow \infty }r\left( U_{m_{j}}\right) ,
\end{eqnarray*}
which satisfy 
\[
\left| a\right| \leq e,\quad e^{-1}\leq \left| \rho \right| \leq e,\quad
\left| r\right| \leq e. 
\]

Then we consider the following subset $K$ of $H$ given by 
\begin{eqnarray*}
K &=&\left\{ \left( U,a,\rho ,r\right) \in H:\frac{1}{c}\leq \left\|
U\right\| \leq c\right. . \\
&&\hspace{1in}\left. \left| a\right| \leq e,\quad e^{-1}\leq \left| \rho
\right| \leq e,\quad \left| r\right| \leq e\right\} .
\end{eqnarray*}
Note that the set $K$ is compact and 
\[
\left( U_{m},a\left( U_{m}\right) ,\rho \left( U_{m}\right) ,r\left(
U_{m}\right) \right) \in K\quad \text{for all }m. 
\]
Then, by Lemma \ref{analytic}, for each $\sigma \in K,$ there exist real
numbers $\varepsilon _{\sigma },\delta _{\sigma }>0$ such that all
normalizations 
\[
\phi _{\sigma ^{\prime }},\quad \sigma ^{\prime }\in K\cap \left\{ \tau \in
GL(n;\Bbb{C}):\left\| \tau -\sigma \right\| \leq \varepsilon _{\sigma
}\right\} 
\]
as a power series at the origin converge absolutely and uniformly on the
open ball $B(0;\delta _{\sigma }).$ Notice that the following family of open
sets 
\[
\left\{ \tau \in GL(n+2;\Bbb{C}):\left\| \tau -\sigma \right\| <\varepsilon
_{\sigma }\right\} ,\quad \sigma \in K 
\]
is an open covering of the set $K.$ Since $K$ is compact, there is a finite
subcover, say, 
\[
\left\{ \tau \in GL(n+2;\Bbb{C}):\left\| \tau -\sigma _{j}\right\|
<\varepsilon _{\sigma _{j}}\right\} ,\quad \sigma _{j}\in H(M),\quad
j=1,\cdots ,l. 
\]
Then we set 
\[
\delta =\min_{1\leq j\leq m}\left\{ \delta _{\sigma _{j}}\right\} >0 
\]
so that each normalization $\phi _{\sigma }$, $\sigma \in K,$ as a power
series at the origin converges absolutely and uniformly on the open ball $%
B(0;\delta ).$ Thus, by Motel theorem, the family of normalizations $\phi
_{\sigma },$ $\sigma \in K$, are a normal family on $B(0;\delta ).$

By a standard argument of a normal family, passing to a subsequence of $%
\{m_{j}\},$ if necessary, there is a holomorphic mapping $\phi $ on the open
ball $B(0;\delta )$ such that 
\[
\phi =\lim_{j\rightarrow \infty }\phi _{\sigma _{m_{j}}} 
\]
where 
\[
\sigma _{m_{j}}=\left( U_{m_{j}},a\left( U_{m_{j}}\right) ,\rho \left(
U_{m_{j}}\right) ,r\left( U_{m_{j}}\right) \right) . 
\]
Then, for $\phi =(f,g),$ we have 
\begin{eqnarray*}
\left( \left. \frac{\partial f}{\partial z}\right| _{0}\right)
&=&\lim_{j\rightarrow \infty }\sqrt{\left| \rho (U_{m_{j}})\right| }%
U_{m_{j}}=\sqrt{\left| \rho \right| }U \\
\left( \left. \frac{\partial f}{\partial w}\right| _{0}\right)
&=&-\lim_{j\rightarrow \infty }\sqrt{\left| \rho (U_{m_{j}})\right| }%
U_{m_{j}}a(U_{m_{j}})=\sqrt{\left| \rho \right| }Ua \\
\left( \left. \frac{\partial g}{\partial w}\right| _{0}\right)
&=&\lim_{j\rightarrow \infty }\rho (U_{m_{j}})=\rho \\
\left( \left. \frac{\partial ^{2}g}{\partial w^{2}}\right| _{0}\right)
&=&2\lim_{j\rightarrow \infty }\rho (U_{m_{j}})r(U_{m_{j}})=2\rho r.
\end{eqnarray*}
Note that 
\[
0<\left| \det \phi ^{\prime }\right| =\left| \rho \right| ^{\frac{n+2}{2}%
}\left| \det U\right| <\infty . 
\]
Thus, by Hurwitz theorem, the mapping $\phi $ is a biholomorphic mapping on
the ball $B(0;\delta ).$ Further, notice 
\[
\phi _{\sigma _{m}}\left( M\cap B(0;\delta )\right) \subset M\quad \text{for
all }m\in \Bbb{N}, 
\]
so that 
\[
\phi \left( M\cap B(0;\delta )\right) \subset M. 
\]
Hence the mapping $\phi $ is a biholomorphic automorphism of $M$ with
initial value $\sigma \in H$ such that 
\[
\sigma =\left( U,a,\rho ,r\right) \in H(M). 
\]
Thus we have showed 
\[
U=\lim_{m\rightarrow \infty }U_{m}\in \left\{ U:\left( U,a,\rho ,r\right)
\in H(M)\right\} . 
\]
Then the group 
\[
\left\{ U:\left( U,a,\rho ,r\right) \in H(M)\right\} \subset GL(n+2;\Bbb{C}) 
\]
is closed so that it is a compact Lie group. Therefore, we prove our claim
that the group $H(M)$ is closed. Hence $H(M)$ is a compact Lie group. This
completes the proof.\endproof

\subsection{Theorem of a germ of a biholomorphic mapping}

We study the analytic continuation of a germ of a biholomorphic mapping to a
finite neighborhood(cf. \cite{Vi85}).

\begin{lemma}
\label{L9}Let $M$ be a nonspherical analytic real hypersurface in normal
form and $H(M)$ be the isotropy subgroup of $M$ such that there is a real
number $c\geq 1$ satisfying 
\[
\sup_{\left( U,a,\rho ,r\right) \in H(M)}\left\| U\right\| \leq c<\infty . 
\]
Then there is a real number $\delta >0$ such that all local automorphisms of 
$M,$ $\phi _{\sigma },$ $\sigma \in H(M)$, converge absolutely and uniformly
on the open ball $B(0;\delta )$.
\end{lemma}

\proof
By Lemma \ref{analytic}, for each $\sigma \in H(M),$ there exist real
numbers $\varepsilon _{\sigma },\delta _{\sigma }>0$ such that all
normalizations 
\[
\phi _{\sigma ^{\prime }},\quad \sigma ^{\prime }\in H(M)\cap \left\{ \tau
\in GL(n+2;\Bbb{C}):\left\| \tau -\sigma \right\| \leq \varepsilon _{\sigma
}\right\} 
\]
as a power series at the origin converges absolutely and uniformly on the
open ball $B(0;\delta _{\sigma }).$

Note that the following family 
\[
\left\{ \tau \in GL(n+2;\Bbb{C}):\left\| \tau -\sigma \right\| <\varepsilon
_{\sigma }\right\} ,\quad \sigma \in H(M) 
\]
is an open covering of the set $H(M).$ By Lemma \ref{compactness}, $H(M)$ is
compact. Thus there is a finite subcover, say, 
\[
\left\{ \tau \in GL(n+2;\Bbb{C}):\left\| \tau -\sigma _{j}\right\|
<\varepsilon _{\sigma _{j}}\right\} ,\quad \sigma _{j}\in H(M),\quad
j=1,\cdots ,m. 
\]
Then we take 
\[
\delta =\min_{1\leq j\leq m}\left\{ \delta _{\sigma _{j}}\right\} >0. 
\]
This completes the proof.\endproof

\begin{theorem}[Vitushkin]
\label{ThVi}Let $M,$ $M^{\prime }$ be a nonspherical analytic real
hypersurface and $p,p^{\prime }$ be points respectively of $M,M^{\prime }$
such that the two germs $M$ at $p$ and $M^{\prime }$ at $p^{\prime }$ are
biholomorphically equivalent. Suppose that there is a real number $c\geq 1$
satisfying 
\[
\sup_{\left( U,a,\rho ,r\right) \in H_{p}(M)}\left\| U\right\| \leq c<\infty 
\]
where $H_{p}(M)$ is a local automorphism group of $M$ at the point $p$ in a
normal coordinate. Then there is a real number $\delta >0$ depending only on 
$M$ and $M^{\prime }$ such that each biholomorphic mapping $\phi $ of $M$
near the point $p$ is analytically continued to the open ball $B(p;\delta )$
whenever $\phi (p)=p^{\prime }$ and there is an open neighborhood $U\subset
B(p;\delta )$ of the point $p$ satisfying 
\[
\phi (U\cap M)\subset M^{\prime }. 
\]
\end{theorem}

\proof
We take a biholomorphic mapping $\phi $ of $M$ to $M^{\prime }$ such that 
\[
\phi (p)=p^{\prime } 
\]
and, for an open neighborhood $U$ of the point $p,$%
\[
\phi (U\cap M)\subset M^{\prime }. 
\]
Then we take normalizations $\phi _{1},\phi _{2}$ respectively of $%
M,M^{\prime }$ such that $\phi _{1},\phi _{2}$ translate the points $%
p,p^{\prime }$ to the origin and there exist open neighborhoods $U_{1},U_{2}$
respectively of $p,p^{\prime }$ and a real hypersurface $M^{*}$ in normal
form satisfying 
\begin{eqnarray*}
\phi _{1}(U_{1}\cap M) &\subset &M^{*} \\
\phi _{2}(U_{2}\cap M) &\subset &M^{*}.
\end{eqnarray*}
Then, we obtain a biholomorphic mapping $\phi ^{*}$ defined by 
\[
\phi ^{*}=\phi _{2}\circ \phi \circ \phi _{1}^{-1}. 
\]
Notice that the mapping $\phi ^{*}$ is a local automorphism of $M^{*}$. By
Lemma \ref{L9}, there is a real number $\delta ^{*}>0$ such that the mapping 
$\phi ^{*}$ continues holomorphically to the open ball $B(0;\delta ^{*})$
satisfying 
\begin{eqnarray*}
B(0;\delta ^{*}) &\subset &\phi _{1}\left( U_{1}\right) \\
B(0;\delta ^{*}) &\subset &\phi _{2}\left( U_{2}\right) .
\end{eqnarray*}
Then the mapping 
\[
\phi =\phi _{2}^{-1}\circ \phi _{\sigma ^{*}}\circ \phi _{1} 
\]
is biholomorphically continued to the open set 
\[
U_{1}\cap \phi _{1}^{-1}\left( B(0;\delta ^{*})\right) . 
\]
We take a real number $\delta >0$ such that 
\[
B(p;\delta )\subset U_{1}\cap \phi _{1}^{-1}\left( B(0;\delta ^{*})\right) . 
\]
This completes the proof.\endproof

\subsection{Kruzhilin-Loboda Theorem}

By Lemma \ref{orbit}, we have a $H$-group action on real hypersurfaces in
normal form. Then the orbit structure in normal form may be studied by
examining the isotropy subgroup $H(M)$ for a real hypersurface $M$ in normal
form.

\begin{lemma}
Let $K$ be a subset of $H.$ The necessary and sufficient condition for the
set $K$ to be conjugate to a subset of 
\[
\left\{ \left( U,0,\pm 1,0\right) \equiv \left( 
\begin{array}{ccc}
\pm 1 & 0 & 0 \\ 
0 & U & 0 \\ 
0 & 0 & 1
\end{array}
\right) \in H\right\} 
\]
is given as follows: 
\[
K\subset \left\{ \left( U,a,\pm 1,r\right) \in H\right\} 
\]
and there exist a vector $d\in \Bbb{C}^{n}$ and a real number $e\in \Bbb{R}$
such that 
\begin{align*}
\left( id_{n\times n}-\lambda U^{-1}\right) d=a \\
\left( 1-\lambda \right) e+i\langle d,a\rangle -i\langle a,d\rangle =r
\end{align*}
for all 
\[
\left( U,a,\lambda ,r\right) \in K. 
\]
\end{lemma}

\proof
Each element of $H$ is decomposed as follows: 
\[
\left( 
\begin{array}{ccc}
\rho ^{\prime } & 0 & 0 \\ 
-C^{\prime }a^{\prime } & C^{\prime } & 0 \\ 
-r^{\prime }-i\langle a^{\prime },a^{\prime }\rangle & 2ia^{\prime \dagger }
& 1
\end{array}
\right) =\left( 
\begin{array}{ccc}
\rho ^{\prime } & 0 & 0 \\ 
0 & C^{\prime } & 0 \\ 
0 & 0 & 1
\end{array}
\right) \left( 
\begin{array}{ccc}
1 & 0 & 0 \\ 
-a^{\prime } & id_{n\times n} & 0 \\ 
-r^{\prime }-i\langle a^{\prime },a^{\prime }\rangle & 2ia^{\prime \dagger }
& 1
\end{array}
\right) 
\]
where 
\[
a^{\prime \dagger }z=\langle z,a^{\prime }\rangle . 
\]
Note that 
\[
\left( 
\begin{array}{ccc}
\rho ^{\prime } & 0 & 0 \\ 
0 & C^{\prime } & 0 \\ 
0 & 0 & 1
\end{array}
\right) \left( 
\begin{array}{ccc}
\pm 1 & 0 & 0 \\ 
0 & U & 0 \\ 
0 & 0 & 1
\end{array}
\right) \left( 
\begin{array}{ccc}
\rho ^{\prime } & 0 & 0 \\ 
0 & C^{\prime } & 0 \\ 
0 & 0 & 1
\end{array}
\right) ^{-1}=\left( 
\begin{array}{ccc}
\pm 1 & 0 & 0 \\ 
0 & C^{\prime }UC^{\prime -1} & 0 \\ 
0 & 0 & 1
\end{array}
\right) . 
\]
Thus the straight forward computation yields 
\begin{eqnarray*}
&&\left( 
\begin{array}{ccc}
1 & 0 & 0 \\ 
-a^{\prime } & id_{n\times n} & 0 \\ 
-r^{\prime }-i\langle a^{\prime },a^{\prime }\rangle & 2ia^{\prime \dagger }
& 1
\end{array}
\right) \left( 
\begin{array}{ccc}
\rho & 0 & 0 \\ 
-Ca & C & 0 \\ 
-r-i\langle a,a\rangle & 2ia^{\dagger } & 1
\end{array}
\right) \times \\
&&\hspace{2in}\left( 
\begin{array}{ccc}
1 & 0 & 0 \\ 
-a^{\prime } & id_{n\times n} & 0 \\ 
-r^{\prime }-i\langle a^{\prime },a^{\prime }\rangle & 2ia^{\prime \dagger }
& 1
\end{array}
\right) ^{-1} \\
&=&\left( 
\begin{array}{ccc}
\rho & 0 & 0 \\ 
-Ca^{*} & C & 0 \\ 
-r^{*}-i\langle a^{*},a^{*}\rangle & 2ia^{*\dagger } & 1
\end{array}
\right)
\end{eqnarray*}
where 
\begin{eqnarray*}
a^{*} &=&\rho C^{-1}a^{\prime }+a-a^{\prime } \\
r^{*} &=&\rho r^{\prime }-r^{\prime }+r+i\langle C(a-a^{\prime }),a^{\prime
}\rangle -i\langle a^{\prime },C(a-a^{\prime })\rangle \\
&&+i\langle a,a^{\prime }\rangle -i\langle a^{\prime },a\rangle .
\end{eqnarray*}
Hence the necessary and sufficient condition for the set $K$ to be conjugate
to a subset of 
\[
\left\{ \left( U,0,\pm 1,0\right) \equiv \left( 
\begin{array}{ccc}
\pm 1 & 0 & 0 \\ 
0 & U & 0 \\ 
0 & 0 & 1
\end{array}
\right) \in H\right\} 
\]
is given by 
\[
\left| \rho \right| =1\quad \text{and}\quad a^{*}=r^{*}=0 
\]
for all 
\[
\left( U,a,\rho ,r\right) \in K. 
\]

The equalities $a^{*}=r^{*}=0$ is yields 
\begin{eqnarray*}
(id_{n\times n}-\rho C^{-1})a^{\prime } &=&a \\
(1-\rho )r^{\prime }+i\langle a^{\prime },a\rangle -i\langle a,a^{\prime
}\rangle &=&r.
\end{eqnarray*}
The necessary and sufficient condition is equivalent to the existence of a
vector $a^{\prime }\in \Bbb{C}^{n}$ and $r^{\prime }\in \Bbb{R}$ satisfying 
\begin{eqnarray*}
\left| \rho \right| &=&1 \\
(id_{n\times n}-\rho C^{-1})a^{\prime } &=&a \\
(1-\rho )r^{\prime }+i\langle a^{\prime },a\rangle -i\langle a,a^{\prime
}\rangle &=&r.
\end{eqnarray*}
for all 
\[
\left( U,a,\rho ,r\right) \in K. 
\]
This completes the proof.\endproof

\begin{theorem}[Kruzhilin-Loboda]
Let $M$ be a real hypersurface in normal form and $H(M)$ be the isotropy
group of $M$ such that there is a real number $c\geq 1$ satisfying 
\[
\sup_{\left( U,a,\rho ,r\right) \in H(M)}\left\| U\right\| \leq c<\infty . 
\]
Then there exists an element $\sigma \in H$ satisfying 
\[
\sigma H(M)\sigma ^{-1}\subset \left\{ \left( 
\begin{array}{ccc}
\pm 1 & 0 & 0 \\ 
0 & U & 0 \\ 
0 & 0 & 1
\end{array}
\right) \in H\right\} . 
\]
\end{theorem}

\proof
By Lemma \ref{compactness}, the group 
\[
G=\left\{ U:\left( U,a,\rho ,r\right) \in H(M)\right\} , 
\]
is a compact Lie group. Thus we have a unique Haar measure $\mu $ on $G$
such that 
\[
\int_{V\in G}d\mu \left( V\right) =1. 
\]

Suppose that $M$ is defined by the equation 
\[
v=\langle z,z\rangle +F_{l}(z,\overline{z},u)+F_{l+1}\left( z,\overline{z}%
,u\right) +O(l+2). 
\]
By Theorem \ref{ThBL}, there is a function $\rho (U)$ satisfying 
\[
\rho =\rho (U) 
\]
for all 
\[
\left( U,a,\rho ,r\right) \in H(M). 
\]
Then we have the following identity 
\[
\left| \rho \left( U\right) \right| ^{\frac{l-2}{2}}F_{l}(z,\overline{z},u)=%
\mathrm{sign\{}\rho \left( U\right) \mathrm{\}}F_{l}\left( U^{-1}z,\overline{%
U^{-1}z},\mathrm{sign\{}\rho \left( U\right) \mathrm{\}}u\right) 
\]
which yields 
\[
F_{l}(z,\overline{z},u)=\frac{\int_{G}\left\{ \mathrm{sign\{}\rho \left(
V\right) \mathrm{\}}F_{l}\left( V^{-1}z,\overline{V^{-1}z},\mathrm{sign\{}%
\rho \left( V\right) \mathrm{\}}u\right) \right\} d\mu \left( V\right) }{%
\int_{G}\left| \rho \left( V\right) \right| ^{\frac{l-2}{2}}d\mu \left(
V\right) }. 
\]
Hence we easily see 
\[
\mathrm{sign\{}\rho \left( U\right) \mathrm{\}}F_{l}\left( U^{-1}z,\overline{%
U^{-1}z},\mathrm{sign\{}\rho \left( U\right) \mathrm{\}}u\right) =F_{l}(z,%
\overline{z},u) 
\]
so that 
\[
\left| \rho \left( U\right) \right| ^{\frac{l-2}{2}}F_{l}(z,\overline{z}%
,u)=F_{l}(z,\overline{z},u). 
\]
Thus we have 
\[
\left| \rho \left( U\right) \right| \equiv 1\quad \text{for all }U\in G 
\]
so that 
\[
H(M)\subset \left\{ \left( U,a,\pm 1,r\right) \in H\right\} . 
\]

By Theorem \ref{ThBL}, there is a function $a(U)$ satisfying 
\[
a=a(U) 
\]
for all 
\[
\left( U,a,\pm 1,r\right) \in H(M). 
\]
Then we have the identity 
\begin{eqnarray*}
&&H_{l+1}\left( z,\overline{z},u;\mathrm{sign\{}\rho \left( U\right) \mathrm{%
\}}Ua\left( U\right) \right) \\
&=&F_{l+1}(z,\overline{z},u)-\mathrm{sign\{}\rho \left( U\right) \mathrm{\}}%
F_{l+1}\left( U^{-1}z,\overline{U^{-1}z},\mathrm{sign\{}\rho \left( U\right) 
\mathrm{\}}u\right) .
\end{eqnarray*}
Hence there is a vector $a^{*}$ satisfying 
\begin{eqnarray*}
&&H_{l+1}\left( z,\overline{z},u;a^{*}\right) \\
&=&F_{l+1}(z,\overline{z},u)-\int_{G}\left\{ \mathrm{sign\{}\rho \left(
V\right) \mathrm{\}}F_{l+1}\left( V^{-1}z,\overline{V^{-1}z},\mathrm{sign\{}%
\rho \left( V\right) \mathrm{\}}u\right) \right\} d\mu (V)
\end{eqnarray*}
where 
\[
a^{*}=\int_{G}\mathrm{sign\{}\rho \left( V\right) \mathrm{\}}Va\left(
V\right) d\mu (V). 
\]
Suppose that the normalization $\phi $ of $M$ with initial value 
\[
(id_{n\times n},-a^{*},1,0)\in H 
\]
transforms $M$ to a real hypersurface $M^{\prime }.$ Then $M^{\prime }$ is
defined up to weight $l+1$ by the equation 
\[
v=\langle z,z\rangle +F_{l}(z,\overline{z},u)+F_{l+1}^{*}\left( z,\overline{z%
},u\right) +O(l+2) 
\]
where 
\[
F_{l+1}^{*}\left( z,\overline{z},u\right) =\int_{G}\left\{ \mathrm{sign\{}%
\rho \left( V\right) \mathrm{\}}F_{l+1}\left( V^{-1}z,\overline{V^{-1}z},%
\mathrm{sign\{}\rho \left( V\right) \mathrm{\}}u\right) \right\} d\mu (V). 
\]
We easily see that 
\[
\mathrm{sign\{}\rho \left( U\right) \mathrm{\}}F_{l+1}^{*}\left( U^{-1}z,%
\overline{U^{-1}z},\mathrm{sign\{}\rho \left( U\right) \mathrm{\}}u\right)
=F_{l+1}^{*}\left( z,\overline{z},u\right) . 
\]
Because the linear mapping $a^{*}\mapsto H_{l+1}\left( z,\overline{z}%
,u;a^{*}\right) $ is injective, we obtain 
\[
H(M^{\prime })\subset \left\{ \left( U,0,\pm 1,r\right) \in H\right\} . 
\]

Suppose that $M^{\prime }$ is defined up to weight $l+2$ by the equation 
\[
F(M^{\prime })=F_{l}(z,\overline{z},u)+F_{l+1}^{*}\left( z,\overline{z}%
,u\right) +F_{l+2}\left( z,\overline{z},u\right) +O(l+3). 
\]
By Theorem \ref{ThBL}, there is a function $r(U)$ satisfying 
\[
r=r(U) 
\]
for all 
\[
\left( U,0,\pm 1,r\right) \in H(M^{\prime }). 
\]
Then we have the following identity 
\begin{align*}
& -\frac{r\left( U\right) }{2}\left\{ \sum_{\min (s,t)\geq
2}(l+s+t)uF_{st}\left( z,\overline{z},u\right) \right. \\
& \hspace{2cm}\left. +\sum_{\min (s,t)\geq 2}2(s-t)i\langle z,z\rangle
F_{st}\left( z,\overline{z},u\right) \right. \\
& \hspace{2cm}\left. -\sum_{\min (s,t)\geq 2}2\langle z,z\rangle ^{2}\left( 
\frac{\partial F_{st}}{\partial u}\right) \left( z,\overline{z},u\right)
\right\} \\
& =\mathrm{sign\{}\rho \left( U\right) \mathrm{\}}F_{l+2}\left( Uz,\overline{%
Uz},\mathrm{sign\{}\rho \left( U\right) \mathrm{\}}u\right) -F_{l+2}\left( z,%
\overline{z},u\right)
\end{align*}
where 
\[
F_{l}(z,\overline{z},u)=\sum_{\min (s,t)\geq 2}F_{st}(z,\overline{z},u). 
\]
Hence there is a real number $r^{*}$ satisfying 
\begin{align*}
& -\frac{r^{*}}{2}\left\{ \sum_{\min (s,t)\geq 2}(l+s+t)uF_{st}\left( z,%
\overline{z},u\right) \right. \\
& \hspace{2cm}\left. +\sum_{\min (s,t)\geq 2}2(s-t)i\langle z,z\rangle
F_{st}\left( z,\overline{z},u\right) \right. \\
& \hspace{2cm}\left. -\sum_{\min (s,t)\geq 2}2\langle z,z\rangle ^{2}\left( 
\frac{\partial F_{st}}{\partial u}\right) \left( z,\overline{z},u\right)
\right\} \\
& =\int_{G}\mathrm{sign\{}\rho \left( V\right) \mathrm{\}}F_{l+2}\left( Vz,%
\overline{Vz},\mathrm{sign\{}\rho \left( V\right) \mathrm{\}}u\right) d\mu
\left( V\right) -F_{l+2}\left( z,\overline{z},u\right)
\end{align*}
where 
\[
r^{*}=\int_{G}r\left( V\right) d\mu \left( V\right) . 
\]
Suppose that the normalization $\phi ^{\prime }$ of $M^{\prime }$ with
initial value 
\[
\left( id,0,1,r^{*}\right) \in H 
\]
transforms $M^{\prime }$ to a real hypersurface $M^{\prime \prime }.$ Then $%
M^{\prime \prime }$ is defined up to weight $l+2$ by the equation 
\[
v=\langle z,z\rangle +F_{l}(z,\overline{z},u)+F_{l+1}^{*}\left( z,\overline{z%
},u\right) +F_{l+2}^{*}\left( z,\overline{z},u\right) +O(l+3) 
\]
where 
\[
F_{l+2}^{*}\left( z,\overline{z},u\right) =\int_{G}\mathrm{sign\{}\rho
\left( V\right) \mathrm{\}}F_{l+2}\left( Vz,\overline{Vz},\mathrm{sign\{}%
\rho \left( V\right) \mathrm{\}}u\right) d\mu \left( V\right) . 
\]
We easily see that 
\[
\mathrm{sign\{}\rho \left( U\right) \mathrm{\}}F_{l+2}^{*}\left( Uz,%
\overline{Uz},\mathrm{sign\{}\rho \left( U\right) \mathrm{\}}u\right)
=F_{l+2}^{*}\left( z,\overline{z},u\right) 
\]
which yields 
\[
H(M^{\prime \prime })\subset \left\{ \left( U,0,\pm 1,0\right) \in H\right\}
. 
\]

Then we take a normalization $\phi _{\sigma }$ with initial value $\sigma
\in H$ such that 
\[
\phi _{\sigma }(M)=M^{\prime \prime }. 
\]
Then, by Lemma \ref{orbit}, we obtain 
\[
\sigma H(M)\sigma ^{-1}=H(M^{\prime \prime }). 
\]
This completes the proof.\endproof

\section{Analytic continuation of a normalizing mapping}

\subsection{Chains on a spherical real hypersurface}

By Theorem \ref{exi-uni}, each biholomorphic automorphism of the real
hyperquadric 
\[
v=\langle z,z\rangle 
\]
is uniquely given by a composition of an affine mapping 
\begin{align}
z^{*}& =z+b  \nonumber \\
w^{*}& =w+2i\langle z,b\rangle +c+i\langle b,b\rangle  \tag*{(4.1)}
\label{line}
\end{align}
and a fractional linear mapping: 
\begin{equation}
\phi =\phi _{\sigma }:\left\{ 
\begin{array}{c}
z^{*}=\frac{C(z-aw)}{1+2i\langle z,a\rangle -w(r+i\langle a,a\rangle )} \\ 
w^{*}=\frac{\rho w}{1+2i\langle z,a\rangle -w(r+i\langle a,a\rangle )}
\end{array}
\right.  \tag*{(4.2)}  \label{2.1}
\end{equation}
where 
\[
b\in \Bbb{C}^{n},\quad c\in \Bbb{R} 
\]
and the constants $\sigma =(C,a,\rho ,r)$ satisfy 
\begin{gather*}
a\in \Bbb{C}^{n},\quad \rho \neq 0,\quad \rho ,r\in \Bbb{R}, \\
C\in GL(n;\Bbb{C}),\quad \langle Cz,Cz\rangle =\rho \langle z,z\rangle .
\end{gather*}
Note that the local automorphism $\phi $ decomposes to 
\[
\phi =\varphi \circ \psi , 
\]
where 
\begin{equation}
\psi :\left\{ 
\begin{array}{c}
z^{*}=\frac{z-aw}{1+2i\langle z,a\rangle -i\langle a,a\rangle w} \\ 
w^{*}=\frac{w}{1+2i\langle z,a\rangle -i\langle a,a\rangle w}
\end{array}
\right. \quad \text{and}\quad \varphi :\left\{ 
\begin{array}{c}
z^{*}=\frac{Cz}{1-rw} \\ 
w^{*}=\frac{\rho w}{1-rw}
\end{array}
\right. .  \nonumber  \label{2.2}
\end{equation}

\begin{lemma}
\label{trans}Let $M$ be the real hyperquadric $v=\langle z,z\rangle .$ Then
the intersection of the real hyperquadric $M$ by a complex line $l$ is given
by a point, a curve $\gamma $, or the complex line $l$ itself. If the
intersection is a curve $\gamma $, then $\gamma $ is transversal to the
complex tangent hyperplane at every point of $\gamma .$
\end{lemma}

\proof
Let $(\kappa ,\chi )\in \Bbb{C}^{n}\times \Bbb{C}$ be a point of the real
hyperquadric $v=\langle z,z\rangle $ such that 
\[
\Im \chi =\langle \kappa ,\kappa \rangle . 
\]
Then a complex line $l$ passing through the point $(\kappa ,\chi )$ is given
by 
\[
\left\{ (\kappa ,\chi )+e(\mu ,\nu ):e\in \Bbb{C}\right\} 
\]
for some nonzero vector $(\mu ,\nu )\in \Bbb{C}^{n}\times \Bbb{C}.$ Then the
affine mapping \ref{line} send the complex line $l$ to another complex line $%
l^{*}$ given by 
\[
\left\{ (\kappa +b,\chi +2i\langle \kappa ,b\rangle +c+i\langle b,b\rangle
)+e(\mu ,\nu +2i\langle \mu ,b\rangle ):e\in \Bbb{C}\right\} . 
\]
Note that the complex line $l^{*}$ passes through the origin by taking 
\[
b=-\kappa ,\quad c=-\Re \chi . 
\]
Thus we reduce the discussion to complex lines passing through the origin.

Suppose that the complex line $l$ is tangent to the complex tangent
hyperplane at the origin so that $l$ is given by 
\[
\left\{ c(a,0):c\in \Bbb{C}\right\} 
\]
for some nonzero vector $a\in \Bbb{C}^{n}.$ Then each point in the
intersection of the real hyperquadric $M$ by the complex line $l$ satisfies 
\[
c\overline{c}\langle a,a\rangle =0. 
\]
Thus we obtain that, whenever $\langle a,a\rangle \neq 0,$%
\[
M\cap \left\{ c(a,0):c\in \Bbb{C}\right\} =\left\{ (0,0)\right\} 
\]
and, whenever $\langle a,a\rangle =0,$%
\[
M\cap \left\{ c(a,0):c\in \Bbb{C}\right\} =\left\{ c(a,0):c\in \Bbb{C}%
\right\} . 
\]
Suppose that the complex line $l$ is transversal to the complex tangent
hyperplane at the origin so that $l$ is given by 
\[
\left\{ c(a,1):c\in \Bbb{C}\right\} 
\]
for some vector $a\in \Bbb{C}^{n}.$ We claim that the complex tangent
hyperplanes of the real hyperquadric $M$ and the complex line $l$ are
transversal at each point of the intersection $\gamma $: 
\[
\gamma =M\cap \left\{ c(a,1):c\in \Bbb{C}\right\} . 
\]
Let $(ca,c),$ $c\neq 0,$ be a point of $M$ so that 
\begin{equation}
\frac{1}{2i}(c+\overline{c})=c\overline{c}\langle a,a\rangle  \tag*{(4.3)}
\label{rela}
\end{equation}
and $(\mu ,\nu )\in \Bbb{C}^{n}\times \Bbb{C}$ be a vector tangent to the
complex tangent hyperplane of $M$ at the point $(ca,c).$ Then we obtain 
\[
\nu -2i\langle \mu ,ca\rangle =0 
\]
so that 
\begin{eqnarray*}
(\mu ,\nu ) &=&\mu _{1}(1,0,\cdots ,0,0,2ie_{1}\overline{c}\overline{a}^{1})
\\
&&+\mu _{2}(0,1,\cdots ,0,0,2ie_{2}\overline{c}\overline{a}^{2}) \\
&&+\cdots \\
&&+\mu _{n}(0,0,\cdots ,0,1,2ie_{n}\overline{c}\overline{a}^{n})
\end{eqnarray*}
where 
\begin{gather*}
\langle z,z\rangle =e_{1}z^{1}\overline{z}^{1}+\cdots +e_{n}z^{n}\overline{z}%
^{n} \\
e_{1},\cdots ,e_{n}=\pm 1.
\end{gather*}
Thus the transversality at $\gamma $ is determined by the value: 
\[
\det \left( 
\begin{array}{ccccc}
1 & 0 & \cdots & 0 & 2ie_{1}\overline{c}\overline{a}^{1} \\ 
0 & 1 & \ddots & \vdots & 2ie_{2}\overline{c}\overline{a}^{2} \\ 
\vdots & \ddots & \ddots & 0 & \vdots \\ 
0 & \cdots & 0 & 1 & 2ie_{n}\overline{c}\overline{a}^{n} \\ 
a^{1} & a^{2} & \cdots & a^{n} & 1
\end{array}
\right) =1-2i\overline{c}\langle a,a\rangle . 
\]
Suppose that 
\[
1-2i\overline{c}\langle a,a\rangle =0. 
\]
Then the equality \ref{rela} yields $c=0.$ This is a contradiction to $c\neq
0$.

Thus the complex tangent hyperplanes of the real hyperquadric $M$ and the
complex line $l$ are transversal at each point of the intersection $\gamma .$
Therefore, the intersection $\gamma $ is a curve transversal to the complex
tangent hyperplanes of $M$ at each point of $\gamma $. This completes the
proof.\endproof

Let $M$ be a nondegenerate analytic real hypersurface in a complex manifold
and $p$ be a point of $M.$ Then we can take a normal coordinate with the
center at the point $p$ so that $M$ is defined by 
\begin{equation}
v=\langle z,z\rangle +\sum_{\min (s,t)\geq 2}F_{st}(z,\bar{z},u)  \nonumber
\end{equation}
where 
\[
\Delta F_{22}=\Delta ^{2}F_{23}=\Delta ^{3}F_{33}=0. 
\]
A connected open curve $\gamma $ on $M$ is called a chain if it is locally
putted into the $u$-curve of a normal coordinate at each point of $\gamma .$
A connected closed subarc of a chain $\gamma $ shall be called a
chain-segment.

\begin{lemma}
\label{chain-chain}Let $M,M^{\prime }$ be nondegenerate analytic real
hypersurfaces and $\phi $ be a biholomorphic mapping on an open neighborhood 
$U$ of a point $p\in M$ such that 
\[
\phi \left( U\cap M\right) \subset M^{\prime }. 
\]
Suppose that there is a chain $\gamma $ of $M$ passing through the point $p.$
Then the analytic curve $\phi \left( U\cap \gamma \right) $ is a chain of $%
M^{\prime }$.
\end{lemma}

\proof
Let $q\in \phi \left( U\cap \gamma \right) .$ Since $\gamma $ is a chain of $%
M,$ there exist an open neighborhood $V$ of the point $\phi ^{-1}(q)\in
\gamma $ and a normalization $\phi _{1}$ of $M$ such that $\phi _{1}$ is
biholomorphic $V$ and 
\[
\phi _{1}\left( V\cap \gamma \right) \subset \left\{ z=v=0\right\} . 
\]
Note that $\phi _{1}\circ \phi ^{-1}$ is a normalization of $M^{\prime }$
such that, for a sufficiently small open neighborhood $O$ of the point $q,$ $%
\phi _{1}\circ \phi ^{-1}$ is biholomorphic on $O$ and 
\[
\phi _{1}\circ \phi ^{-1}\left( O\cap \phi \left( V\cap \gamma \right)
\right) \subset \left\{ z=v=0\right\} . 
\]
Since $q$ is an arbitrary point of $\phi \left( U\cap \gamma \right) ,$ the
analytic curve $\phi \left( U\cap \gamma \right) $ is a chain. This
completes the proof.\endproof

\begin{lemma}
\label{trans2}Let $M$ be a real hyperquadric $v=\langle z,z\rangle $ and $%
p=(\kappa ,\chi )\neq 0$ be a point of $M.$ Then there is a chain-segment $%
\gamma :[0,1]\rightarrow M$ such that 
\[
\gamma (0)=0,\quad \gamma (1)=p 
\]
whenever 
\[
\Re \chi \neq 0\quad \text{or}\quad \langle \kappa ,\kappa \rangle \neq 0. 
\]
\end{lemma}

\proof
Let $\phi _{\sigma }$ be a local automorphism of a real hyperquadric with
initial value $\sigma \in H$. Then the inverse $\phi _{\sigma }^{-1}$ of the
local automorphism $\phi _{\sigma }$ is given by 
\begin{equation}
\phi _{\sigma }^{-1}:\left\{ 
\begin{array}{c}
z=\frac{C^{-1}(z^{*}+\rho ^{-1}Caw^{*})}{1-2i\langle z^{*},\rho
^{-1}Ca\rangle -w^{*}(-r\rho ^{-1}+i\langle \rho ^{-1}Ca,\rho ^{-1}Ca\rangle
)} \\ 
w=\frac{\rho ^{-1}w^{*}}{1-2i\langle z^{*},\rho ^{-1}Ca\rangle -w^{*}(-r\rho
^{-1}+i\langle \rho ^{-1}Ca,\rho ^{-1}Ca\rangle )}
\end{array}
\right. .  \tag*{(4.4)}  \label{inverse}
\end{equation}
Thus the chain passing through the origin and transversal to the complex
tangent hyperplane at the origin, $\gamma ,$ is given with a normal
parametrization by 
\[
\gamma =\phi ^{-1}\left( z^{*}=v^{*}=0\right) :\left\{ 
\begin{array}{c}
z=\frac{\rho ^{-1}au^{*}}{1-\rho ^{-1}u^{*}(-r+i\langle a,a\rangle )} \\ 
w=\frac{\rho ^{-1}u^{*}}{1-\rho ^{-1}u^{*}(-r+i\langle a,a\rangle )}
\end{array}
\right. . 
\]
By taking $r=0,$ we easily see that the chain $\gamma $ is the intersection
of $M$ and the complex line 
\[
\left\{ c(a,1):c\in \Bbb{C}\right\} . 
\]
By Lemma \ref{trans}, the chain $\gamma $ is transversal to the complex
tangent hyperplanes of $M$ at each point of $\gamma .$

Let $(\kappa ,\chi )\neq 0$ be a point in $\Bbb{C}^{n}\times \Bbb{C}$ on $M$
such that 
\[
\Im \chi =\langle \kappa ,\kappa \rangle . 
\]
Then we have $\chi \neq 0$ whenever 
\[
\Re \chi \neq 0\quad \text{or}\quad \langle \kappa ,\kappa \rangle \neq 0. 
\]
Note that the origin and the point $(\kappa ,\chi ),$ $\chi \neq 0,$ is
connected by the chain 
\[
\Gamma =M\cap \left\{ c(\chi ^{-1}\kappa ,1):c\in \Bbb{C}\right\} . 
\]
This completes the proof.\endproof

\begin{lemma}
\label{spherical}Let $M$ be a spherical analytic real hypersurface. Then $M$
is locally biholomorphic to a real hyperquadric.
\end{lemma}

In the paper \cite{Pa2}, we have proved Lemma \ref{spherical}.

\begin{theorem}
\label{core2-lemma}Let $M$ be a spherical analytic real hypersurface and $%
\gamma :[0,1]\rightarrow M$ be a curve such that $\gamma [0,\tau ]$ is a
chain-segment for each $\tau <1.$ Then $\gamma [0,1]$ is a chain-segment of $%
M.$
\end{theorem}

\proof
By Lemma \ref{spherical}, the real hypersurface $M$ is biholomorphic to a
real hyperquadric at the point $\gamma (1).$ Then, by Lemma \ref{spherical},
taking a normal coordinate with center at the point $\gamma (1)$ yields 
\[
v=\langle z,z\rangle 
\]
where the curve $\gamma [0,1]$ touches the origin and the part $\gamma (0,1)$
is a chain. By Lemma \ref{trans} and Lemma \ref{trans2}, there exist a chain 
$\Gamma $ passing through the origin, an open neighborhood $U$ of the
origin, and a normalization $\phi $ of $M$ such that $\phi $ is
biholomorphic on $U$ and 
\[
\phi \left( \Gamma \cap U\right) \subset \left\{ z=v=0\right\} . 
\]
Since $\gamma [0,1)\subset \Gamma $ and $\Gamma $ is a chain of $M,$ $\gamma
[0,1]$ is a chain-segment. This completes the proof.\endproof

\begin{proposition}
\label{cone}Let $M$ be a spherical analytic real hypersurface and $p$ be a
point of $M.$ Suppose that there are an open cone $V_{\theta }$ with its
vertex at the point $p$ and euclidean angle $\theta ,$ $0<\theta <\frac{\pi 
}{2},$ to the complex tangent hyperplane at the point $p,$ and an open
neighborhood $U$ of the point $p.$ Then there is a number $\delta >0$ such
that, for each given curve $\eta :[0,1]\rightarrow V_{\theta }\cap
B(p;\delta ),$ there is a continuous family of chain-segments 
\[
\gamma :[0,1]\times [0,1]\rightarrow U\cap M 
\]
where $\gamma (s,\cdot ):[0,1]\rightarrow U\cap M$ is a chain-segment of $M$
for each $s\in [0,1]$ satisfying 
\[
\gamma (s,0)=p\quad \text{and}\quad \gamma (s,1)=\eta (s)\quad \text{for all 
}s\in [0,1]. 
\]
\end{proposition}

\proof
By Lemma \ref{spherical}, there is a biholomorphic mapping $\phi $ of $M$
near the point $p$ to a real hyperquadric $v=\langle z,z\rangle $ such that $%
N(p)=0.$ Then we take a sufficiently small number $\varepsilon >0$ so that
each point $q\in \phi (V_{\theta })\cap B(0;\varepsilon )$ is connected by a
chain-segment $\gamma \subset \phi (U\cap M)$ to the origin$.$ Then we take
a number $\delta >0$ such that 
\[
\phi \left( V_{\theta }\cap B(p;\delta )\right) \subset \phi (V_{\theta
})\cap B(0;\varepsilon ). 
\]
By Lemma \ref{trans2}, there is a continuous family of complex line $l_{\tau
}$ such that the intersection $l_{\tau }\cap \phi (U\cap M)$ is a chain and
the chain $l_{\tau }\cap \phi (U\cap M)$ connects the point $p$ and the
point $\phi (\eta (\tau ))$ for each $\tau \in [0,1].$ Hence there is a
continuous family of chain-segments 
\[
\Gamma :[0,1]\times [0,1]\rightarrow \phi (U\cap M) 
\]
where $\Gamma (\tau ,\cdot ):[0,1]\rightarrow \phi (U\cap M)$ is a
chain-segment for each $\tau \in [0,1]$ satisfying 
\[
\Gamma (\tau ,0)=0\quad \text{and}\quad \Gamma (\tau ,1)=\phi (\eta (\tau
))\quad \text{for all }\tau \in [0,1]. 
\]
Then, by Lemma \ref{chain-chain}, the desired family of chain-segments on $M$
is given by 
\[
\gamma \equiv \phi ^{-1}\circ \Gamma :[0,1]\times [0,1]\rightarrow U\cap M. 
\]
This completes the proof.\endproof

\subsection{Chains on a nonspherical real hypersurface}

\begin{lemma}
\label{core1-lemma}Let $M$ be a nondegenerate analytic real hypersurface and 
$p$ be a point of $M.$ Suppose that there are an open cone $V_{\theta }$
with its vertex at the point $p$ and euclidean angle $\theta ,$ $0<\theta <%
\frac{\pi }{2},$ to the complex tangent hyperplane at the point $p,$ and an
open neighborhood $U$ of the point $p.$ Then there is a number $\delta >0$
such that, for each given curve $\eta :[0,1]\rightarrow V_{\theta }\cap
B(p;\delta ),$ there is a continuous family of chain-segments 
\[
\gamma :[0,1]\times [0,1]\rightarrow U\cap M 
\]
where $\gamma (s,\cdot ):[0,1]\rightarrow U\cap M$ is a chain-segment of $M$
for each $s\in [0,1]$ satisfying 
\[
\gamma (s,0)=p\quad \text{and}\quad \gamma (s,1)=\eta (s)\quad \text{for all 
}s\in [0,1]. 
\]
\end{lemma}

\proof
By translation and unitary transformation, if necessary, we may assume that
the point $p$ is at the origin and the real hypersurface $M$ is defined near
the origin by 
\[
v=F(z,\overline{z},u),\quad \left. F\right| _{0}=\left. F_{z}\right|
_{0}=\left. F_{\overline{z}}\right| _{0}=0 
\]
so that 
\[
F(z,\overline{z},u)=\sum_{s=2}^{\infty }F_{2}(z,\overline{z},u). 
\]

With a sufficiently small number $\varepsilon >0,$ we consider an analytic
family of real hypersurfaces $M_{\mu },$ $\left| \mu \right| \leq
\varepsilon ,$ defined near the origin by the equations: 
\[
v=F^{*}(z,\overline{z},u;\mu ) 
\]
where 
\[
F^{*}(z,\overline{z},u;\mu )=\sum_{s=2}^{\infty }\mu ^{k-2}F_{k}(z,\overline{%
z},u). 
\]
Note that the function $F^{*}(z,\overline{z},u;\mu )$ is analytic of $%
z,u,\mu $ and the real hypersurface $M_{0}$(i.e., $\mu =0$) is spherical.

Then we obtain an analytic family of ordinary differential equations 
\begin{equation}
p^{\prime \prime }=Q(\tau ,p,\overline{p},p^{\prime },\overline{p}^{\prime
};\mu )  \tag*{(4.5)}  \label{chain eq}
\end{equation}
so that each chain $\gamma $ passing through the origin on $M_{\mu }$ is
given by 
\[
\gamma :\left\{ 
\begin{array}{l}
z=p(\tau ) \\ 
w=\tau +iF^{*}\left( p(\tau ),\overline{p}(\tau ),\tau ;\mu \right)
\end{array}
\right. 
\]
where $p(\tau )$ is a solution of the equation \ref{chain eq}. The solution $%
p$ of the equation \ref{chain eq} is given as an analytic function of $\tau
,\mu ,a,$ where 
\[
a=p^{\prime }(0). 
\]
In fact, for a given real number $\nu \in \Bbb{R}^{+},$ there are real
numbers $\tau _{1},\varepsilon _{1}$ such that the analytic function 
\[
p=p(\tau ,\mu ,a) 
\]
converges absolutely and uniformly on the range 
\[
\left| a\right| \leq \nu ,\quad \left| \tau \right| \leq \tau _{1},\quad
\left| \mu \right| \leq \varepsilon _{1}. 
\]

Since $M_{0}$ is spherical, by Theorem \ref{cone}, for an open neighborhood $%
U_{0}$ of the origin and an open cone $V_{\theta _{0}}$ with its vertex at
the origin and euclidean angle $\theta _{0},$ $0<\theta _{0}<\frac{\pi }{2},$
to the complex tangent hyperplane at the origin$,$ there is a number $\delta
_{0}>0$ such that, for each given curve $\eta :[0,1]\rightarrow V_{\theta
_{0}}\cap B(0;\delta _{0}),$ there is a continuous family of chain-segments 
\[
\gamma _{0}:[0,1]\times [0,1]\rightarrow U_{0}\cap M_{0} 
\]
where $\gamma _{0}(s,\cdot ):[0,1]\rightarrow U_{0}\cap M_{0}$ is a
chain-segment of $M_{0}$ for each $s\in [0,1]$ satisfying 
\[
\gamma _{0}(s,0)=0\quad \text{and}\quad \gamma _{0}(s,1)=\eta (s)\quad \text{%
for all }s\in [0,1]. 
\]
Then, for an open neighborhood $U_{1}$ of the origin and an open cone $%
V_{\theta _{1}}$ with its vertex at the origin and euclidean angle $\theta
_{1},$ $0<\theta _{1}<\frac{\pi }{2},$ to the complex tangent hyperplane at
the origin, there exist real numbers $\mu _{1},\delta _{1}>0$ such that, for
each given curve $\eta :[0,1]\rightarrow V_{\theta _{1}}\cap B(0;\delta
_{1}),$ there is a continuous family of chain-segments 
\[
\gamma _{1}:[0,1]\times [0,1]\rightarrow U_{1}\cap M_{\mu _{1}} 
\]
where $\gamma _{1}(s,\cdot ):[0,1]\rightarrow U_{1}\cap M_{\mu _{1}}$ is a
chain-segment of $M_{\mu _{1}}$ for each $s\in [0,1]$ satisfying 
\[
\gamma _{1}(s,0)=0\quad \text{and}\quad \gamma _{1}(s,1)=\eta (s)\quad \text{%
for all }s\in [0,1]. 
\]

By the way, the real hypersurface $M_{\mu },$ $\mu \neq 0,$ is obtained from 
$M$ by the biholomorphic mapping: 
\[
\chi _{\mu }:\left\{ 
\begin{array}{l}
z^{*}=\mu ^{-1}z \\ 
w^{*}=\mu ^{-2}w
\end{array}
\right. . 
\]
For an open neighborhood $U$ of the origin and an open cone $V_{\theta }$
with its vertex at the origin and euclidean angle $\theta ,$ $0<\theta <%
\frac{\pi }{2},$ to the complex tangent hyperplane at the origin, we take $%
\theta _{1}$ and a real number $\delta >0$ such that 
\[
\chi _{\mu _{1}}\left( V_{\theta }\cap B(0;\delta )\right) \subset V_{\theta
_{1}}\cap B(0;\delta _{1}). 
\]
Then, for each given curve $\eta :[0,1]\rightarrow V_{\theta }\cap
B(0;\delta ),$ there is a continuous family of chain-segments 
\[
\gamma _{1}:[0,1]\times [0,1]\rightarrow U_{1}\cap M_{\mu _{1}} 
\]
where $\gamma _{1}(s,\cdot ):[0,1]\rightarrow \chi _{\mu _{1}}\left( U\cap
M\right) $ is a chain-segment of $M_{\mu _{1}}$ for each $s\in [0,1]$
satisfying 
\[
\gamma _{1}(s,0)=0\quad \text{and}\quad \gamma _{1}(s,1)=\chi _{\mu
_{1}}\left( \eta (s)\right) \quad \text{for all }s\in [0,1]. 
\]
Then, by Lemma \ref{chain-chain}, the desired family of chain-segments on $M$
is given by 
\[
\gamma \equiv \chi _{\mu _{1}}^{-1}\circ \gamma _{1}:[0,1]\times
[0,1]\rightarrow U\cap M. 
\]
This completes the proof.\endproof

\begin{theorem}
\label{core1}Let $M$ be an analytic real hypersurface with nondegenerate
Levi form and $U$ be an open neighborhood of a point $p$ of $M.$ Then there
are a number $\varepsilon >0$ and a point $q\in U\cap M$ such that 
\[
B(p;\varepsilon )\subset U 
\]
and, for each given curve $\eta :[0,1]\rightarrow B(p;\varepsilon )\cap M,$
there is a continuous family of chain-segments 
\[
\gamma :[0,1]\times [0,1]\rightarrow U\cap M 
\]
where $\gamma (s,\cdot ):[0,1]\rightarrow U\cap M$ is a chain-segment of $M$
for each $s\in [0,1]$ satisfying 
\[
\gamma (s,0)=q\quad \text{and}\quad \gamma (s,1)=\eta (s)\quad \text{for all 
}s\in [0,1]. 
\]
\end{theorem}

\proof
We take a point $q$ sufficiently near the point $p$ such that there are an
open cone $V_{\theta }$ and a number $\delta >0$ in Lemma \ref{core1-lemma}
satisfying 
\[
p\in V_{\theta }\cap B(q;\delta ). 
\]
Then we take a number $\varepsilon >0$ such that 
\[
B(p;\varepsilon )\in V_{\theta }\cap B(q;\delta ). 
\]
The desired result follows from Lemma \ref{core1-lemma}. This completes the
proof.\endproof

\subsection{Piecewise chain curve}

Let $M$ be an analytic real hypersurface with nondegenerate Levi form. Let $%
\gamma $ be a piecewise differentiable curve of $[0,1]$ into $M$ such that
there are disjoint open intervals $I_{i},$ $i=1,\cdots ,m,$ satisfying 
\[
\lbrack 0,1]=\bigcup_{i=1}^{m}\overline{I_{i}} 
\]
and each fraction $\gamma \left( I_{i}\right) ,$ $i=1,\cdots ,m,$ is a
chain-segment. Then $\gamma $ shall be called a piecewise chain curve.

\begin{lemma}
\label{pcc}Let $M$ be a connected analytic real hypersurface and $\Gamma $
be a continuous curve on $M$ connecting two points $p,q\in M.$ Then, for a
given number $\varepsilon >0,$ there is a piecewise chain curve $\gamma
:[0,1]\rightarrow M$ such that 
\begin{align*}
\gamma (0)=p,\quad \gamma (1)=q \\
\gamma [0,1]\subset \bigcup_{x\in \Gamma }B(x;\varepsilon ).
\end{align*}
\end{lemma}

\proof
Since the curve $\Gamma $ is compact, there are finitely many points $%
x_{i}\in \Gamma ,$ $i=1,\cdots ,l,$ such that 
\[
\Gamma \subset \bigcup_{i=1}^{l}B(x_{i};\varepsilon ). 
\]
Suppose that $x$ is a point on $\Gamma $ and $x\in B(x_{i};\varepsilon ).$
Then, by Lemma \ref{core1}, there is a number $\delta _{x}>0$ such that
every two points $y,z\in B(x;\delta _{x})$ are connected by at most $2$%
-pieced chain curve $\gamma \subset B(x_{i};\varepsilon ).$

Note that the set $\left\{ B(x;\delta _{x}):x\in \Gamma \right\} $ is an
open covering of $\Gamma .$ Since $\Gamma $ is compact, there is a finite
subcover, say, 
\[
\Gamma \subset \bigcup_{j=1}^{k}B(y_{j};\delta _{y_{j}}). 
\]
Then there is at most $2k$-pieced chain curve $\gamma :[0,1]\rightarrow M$
connecting the two point $p,q\in M$ such that 
\[
\gamma [0,1]\subset \bigcup_{x\in \Gamma }B(x;\varepsilon ). 
\]
This completes the proof.\endproof

\begin{lemma}
\label{pcc2}Let $M$ be a nondegenerate analytic real hypersurface$.$ Suppose
that there is a piecewise chain curve $\gamma $ connecting two points $%
p,q\in M.$ Then $M$ is biholomorphic to a real hyperquadric at the point $p$
if and only if $M$ is biholomorphic to a real hyperquadric at the point $q.$
\end{lemma}

\proof
Without loss of generality, we may assume that $p$ and $q$ are connected by
a chain-segment $\gamma :[0,1]\rightarrow M$. Then there is a chain $\Gamma $
of $M$ satisfying 
\[
\gamma [0,1]\subset \Gamma . 
\]
For each point $x\in \Gamma ,$ there are an open neighborhood $U_{x}$ of $x$
and a biholomorphic mapping $N_{x}$ such that 
\begin{eqnarray*}
N_{x}(x) &=&0 \\
N_{x}(U_{x}\cap \Gamma ) &\subset &\left\{ z=v=0\right\} .
\end{eqnarray*}
Since the subset $\gamma [0,1]$ is compact, there is a finite subcover, say, 
\[
\left\{ U_{x_{i}}:i=1,\cdots ,m\right\} . 
\]
Suppose that the normalization $N_{x_{i}}$ transforms $M\cap U_{x_{i}}$ to
the real hypersurface $M_{x_{i}}^{\prime }$ defined near the origin by 
\[
v=\langle z,z\rangle +\sum_{s,t\geq 2}F_{st}^{i}(z,\overline{z},u). 
\]
Note that the functions $F_{st}^{i}(z,\overline{z},u)$ are analytic of $u$
on the set $N_{x_{i}}(U_{x_{i}}\cap \Gamma ).$ Thus 
\[
F_{st}^{i}(z,\overline{z},u)\equiv 0 
\]
whenever there is an open subset $U\subset U_{x_{i}}$ satisfying 
\[
F_{st}^{i}(z,\overline{z},u)=0\quad \text{for }u\in N_{x_{i}}(U\cap \Gamma
). 
\]
Thus the desired result follows. This completes the proof.\endproof

\begin{theorem}
\label{remove2}Let $M$ be a connected nondegenerate analytic real
hypersurface$.$ Then $M$ is not biholomorphic to a real hyperquadric at each
point of $M$ whenever there is a point $p$ of $M$ at which $M$ is not
biholomorphic to a real hyperquadric.
\end{theorem}

\proof
The contrapositive may be stated as follows: $M$ is locally biholomorphic to
a real hyperquadric at each point of $M$ whenever $M$ is biholomorphic to a
real hyperquadric at a point $p$ of $M.$ Suppose that there is a point $p$
of $M$ at which $M$ is biholomorphic to a real hyperquadric. By Lemma \ref
{pcc}, each point $q$ of $M$ is connected to $p$ by a piecewise chain curve.
Then, by Lemma \ref{pcc2}, $M$ is biholomorphic to a real hyperquadric at
the point $q$ as well. Since $M$ is connected, this completes the proof.\endproof

\begin{lemma}
\label{umbilic}Let $M$ be an analytic real hypersurface and $U$ be an open
neighborhood of a point $p\in M.$ Suppose that $U\cap M$ consists of umbilic
points. Then $U\cap M$ is locally biholomorphic to a real hyperquadric.
\end{lemma}

In the paper \cite{Pa2}, we have given the proof of Lemma \ref{umbilic}.

\begin{theorem}
Let $M$ be a nondegenerate analytic real hypersurface and $p$ be a point of $%
M.$ Suppose that $M$ is not biholomorphic to a real hyperquadric at the
point $p.$ Then there is a normalization $\phi $ near the point $p$ such
that $\phi \left( M\right) $ is defined by the equation, for $\dim M=3,$%
\[
v=\langle z,z\rangle +\sum_{\min (s,t)\geq 2,\max (s,t)\geq 4}F_{st}(z,%
\overline{z},u) 
\]
where 
\[
F_{24}(z,\overline{z},u)\neq 0, 
\]
and, for $\dim M\geq 5,$%
\[
v=\langle z,z\rangle +\sum_{\min (s,t)\geq 2}F_{st}(z,\overline{z},u) 
\]
where 
\[
F_{22}(z,\overline{z},u)\neq 0. 
\]
\end{theorem}

\proof
Suppose that the assertion is not true. Then $M$ is umbilic on all points of
all chains passing through the point $p.$ Then, by Theorem \ref{core1-lemma}%
, there is an open set $U$ such that every point of $U\cap M$ is connected
to $p$ by a chain of $M.$ Hence $U\cap M$ consists of umbilic points so
that, by Lemma \ref{umbilic}, $U\cap M$ is locally biholomorphic to a real
hyperquadric. By Lemma \ref{pcc2}, $M$ is biholomorphic to a real
hyperquadric at the point $p$ as well. This is a contradiction. This
completes the proof.\endproof

\subsection{Global straightening of a chain}

Let $M$ be a nondegenerate analytic real hypersurface defined near the
origin by the equation 
\[
v=\frac{1}{4\alpha }\ln \frac{1}{1-4\alpha \langle z,z\rangle }+\sum_{\min
(s,t)\geq 2}F_{st}(z,\bar{z},u) 
\]
where $\alpha $ is a given real number and 
\[
\Delta F_{22}=\Delta ^{2}F_{23}=\Delta ^{3}F_{33}=0. 
\]
By using the expansion 
\[
-\ln \left( 1-x\right) =\sum_{m=1}^{\infty }\frac{x^{m}}{m}, 
\]
the defining equation of $M$ comes to 
\[
v=\langle z,z\rangle +\sum_{\min (s,t)\geq 2}F_{st}^{*}(z,\bar{z},u) 
\]
where 
\begin{eqnarray*}
\Delta F_{22}^{*}(z,\overline{z},u) &=&4\alpha (n+1)\langle z,z\rangle \\
\Delta ^{2}F_{23}^{*}(z,\overline{z},u) &=&0 \\
\Delta ^{3}F_{33}^{*}(z,\overline{z},u) &=&32\alpha ^{2}n(n+1)(n+2).
\end{eqnarray*}
We may require the maximal analytic extension along the $u$-curve on the
real hypersurface $M$ in Moser-Vitushkin normal form.

\begin{lemma}
\label{linear-linear}Let $M$ be a real hypersurface defined near the origin
by 
\[
v=\langle z,z\rangle +\sum_{\min (s,t)\geq 2}F_{st}(z,\bar{z},u) 
\]
where 
\[
\Delta F_{22}=\Delta ^{2}F_{23}=\Delta ^{3}F_{33}=0. 
\]
Let $\frak{L}$ be the mapping 
\[
\frak{L}:\left\{ 
\begin{array}{l}
z^{*}=\frac{z}{1-i\alpha w} \\ 
w^{*}=\frac{1}{2i\alpha }\ln \frac{1+i\alpha w}{1-i\alpha w}=\frac{1}{\alpha 
}\tan ^{-1}\alpha w
\end{array}
\right. . 
\]
Then $\frak{L}\left( M\right) $ is defined near the origin by the equation 
\[
v=\frac{1}{4\alpha }\ln \frac{1}{1-4\alpha \langle z,z\rangle }+\sum_{\min
(s,t)\geq 2}^{\infty }F_{st}^{*}(z,\bar{z},u) 
\]
where 
\[
\Delta F_{22}^{*}=\Delta ^{2}F_{23}^{*}=\Delta ^{3}F_{33}^{*}=0. 
\]
\end{lemma}

\proof
Suppose that the real hypersurface $M$ is in Chern-Moser normal form is
defined by the equation 
\[
v=\langle z,z\rangle +\sum_{\min (s,t)\geq 2}F_{st}(z,\bar{z},u) 
\]
where 
\begin{equation}
\Delta F_{22}=\Delta ^{2}F_{23}=\Delta ^{3}F_{33}=0.  \tag*{(4.6)}
\label{un}
\end{equation}
Let $\frak{L}$ be a normalization of $M$ to Moser-Vitushkin normal form
leaving the $u$-curve invariant. We require the identity initial value on
the normalization $\frak{L}$ so that $\frak{L}$ is necessarily of the
form(cf. the proof of Theorem \ref{exi-uni}) 
\[
\frak{L}:\left\{ 
\begin{array}{l}
z^{*}=\sqrt{q^{\prime }(w)}U(w)z \\ 
w^{*}=q(w)
\end{array}
\right. 
\]
where 
\[
\langle U(u)z,U(u)z\rangle =\langle z,z\rangle \quad \mathrm{and}\quad
U(0)=id_{n\times n}. 
\]
Suppose that $\frak{L}$ transforms $M$ to a real hypersurface $M^{\prime }$
defined by 
\[
v=\langle z,z\rangle +\sum_{\min (s,t)\geq 2}F_{st}^{*}(z,\bar{z},u). 
\]
Then we obtain 
\begin{align}
F_{22}(z,\overline{z},u)& =q^{\prime }(u)F_{22}^{*}(U(u)z,\overline{U(u)z}%
,q(u))  \nonumber \\
& -i\langle z,z\rangle \langle z,U(u)^{-1}\{U^{\prime }(u)+\frac{1}{2}%
q^{\prime }(u)^{-1}q^{\prime \prime }(u)U(u)\}z\rangle  \nonumber \\
& +i\langle z,z\rangle \langle U(u)^{-1}\{U^{\prime }(u)+\frac{1}{2}%
q^{\prime }(u)^{-1}q^{\prime \prime }(u)U(u)\}z,z\rangle  \tag*{(4.7)}
\label{deux}
\end{align}
and 
\begin{equation}
F_{23}(z,\overline{z},u)=q^{\prime }(u)\sqrt{\left| q^{\prime }(u)\right| }%
F_{23}^{*}(U(u)z,\overline{U(u)z},q(u)).  \tag*{(4.8)}  \label{trois}
\end{equation}
The condition $\Delta ^{2}F_{23}=0$ in \ref{un} yields 
\[
\Delta ^{2}F_{23}^{*}(z,\overline{z},u)=0 
\]
so that the $u$-curve is a chain of $M^{\prime }.$ We require that $%
M^{\prime }$ is in Moser-Vitushkin normal form so that 
\begin{eqnarray*}
\Delta F_{22}^{*}(z,\overline{z},u) &=&4\alpha (n+1)\langle z,z\rangle \\
\Delta ^{2}F_{23}^{*}(z,\overline{z},u) &=&0 \\
\Delta ^{3}F_{33}^{*}(z,\overline{z},u) &=&32\alpha ^{2}n(n+1)(n+2).
\end{eqnarray*}
Then, in \ref{deux}, we require the following condition 
\begin{eqnarray*}
\Delta F_{22}(z,\overline{z},u) &=&0 \\
\Delta F_{22}^{*}(z,\overline{z},u) &=&4\alpha (n+1)\langle z,z\rangle
\end{eqnarray*}
so that 
\begin{align}
& 4\alpha (n+1)q^{\prime }(u)\langle z,z\rangle  \nonumber \\
& +i(n+2)\{\langle U(u)^{-1}U^{\prime }(u)z,z\rangle -\langle
z,U(u)^{-1}U^{\prime }(u)z\rangle \}  \nonumber \\
& +i\langle z,z\rangle \{\mathrm{Tr}(U(u)^{-1}U^{\prime }(u))-\overline{%
\mathrm{Tr}(U(u)^{-1}U^{\prime }(u))}\}  \nonumber \\
& =0.  \tag*{(4.9)}  \label{quatre}
\end{align}
From the condition $\langle U(u)z,U(u)z\rangle =\langle z,z\rangle ,$ we
obtain 
\begin{eqnarray*}
\langle U(u)^{-1}U^{\prime }(u)z,z\rangle +\langle z,U(u)^{-1}U^{\prime
}(u)z\rangle &=&0 \\
\mathrm{Tr}(U(u)^{-1}U^{\prime }(u))+\overline{\mathrm{Tr}%
(U(u)^{-1}U^{\prime }(u))} &=&0.
\end{eqnarray*}
The equality \ref{quatre} comes to 
\[
2\alpha i(n+1)q^{\prime }(u)id_{n\times n}=(n+2)U(u)^{-1}U^{\prime }(u)+%
\mathrm{Tr}(U(u)^{-1}U^{\prime }(u))id_{n\times n}, 
\]
which yields 
\[
\mathrm{Tr}(U(u)^{-1}U^{\prime }(u))=\alpha niq^{\prime }(u). 
\]
Hence we obtain 
\[
U^{\prime }(u)=\alpha iq^{\prime }(u)U(u). 
\]
Thus the function $U(u)$ is given by 
\[
U(u)=\exp \alpha iq(u). 
\]

Then the mapping $\frak{L}$ is necessarily of the form: 
\[
\frak{L}:\left\{ 
\begin{array}{l}
z^{*}=\sqrt{q^{\prime }(w)}z\exp \alpha iq(w) \\ 
w^{*}=q(w)
\end{array}
\right. . 
\]
Then we have 
\begin{align}
F_{33}(z,\overline{z},u)& =q^{\prime }(u)\left| q^{\prime }(u)\right|
F_{33}^{*}(z,\overline{z},q(u))  \nonumber \\
& -6\alpha q^{\prime }(u)^{2}\langle z,z\rangle F_{22}^{*}(z,\overline{z}%
,q(u))  \nonumber \\
& +\left\{ -\frac{q^{\prime \prime \prime }(u)}{3q^{\prime }(u)}+\frac{1}{2}%
\left( \frac{q^{\prime \prime }(u)}{q^{\prime }(u)}\right) ^{2}+6\alpha
^{2}q^{\prime }(u)^{2}\right\} \langle z,z\rangle ^{3}.  \tag*{(4.10)}
\label{cinq}
\end{align}
We have the following identities: 
\begin{eqnarray*}
\Delta ^{3}\langle z,z\rangle ^{3} &=&6n(n+1)(n+2) \\
\Delta ^{3}\left\{ F_{33}^{*}(z,\overline{z},q(u))\right\} &=&\Delta
^{3}F_{33}^{*}(z,\overline{z},q(u)) \\
\Delta ^{3}\left\{ \langle z,z\rangle F_{22}^{*}(z,\overline{z}%
,q(u))\right\} &=&3(n+2)\Delta ^{2}F_{22}^{*}(z,\overline{z},q(u)).
\end{eqnarray*}
Then, by requiring in \ref{cinq} the following conditions 
\[
\Delta ^{3}F_{33}(z,\overline{z},u)=0 
\]
and 
\begin{eqnarray*}
\Delta F_{22}^{*}(z,\overline{z},u) &=&4\alpha (n+1)\langle z,z\rangle \\
\Delta ^{3}F_{33}^{*}(z,\overline{z},u) &=&32\alpha ^{2}n(n+1)(n+2),
\end{eqnarray*}
we obtain 
\[
\frac{q^{\prime \prime \prime }(u)}{3q^{\prime }(u)}-\frac{1}{2}\left( \frac{%
q^{\prime \prime }(u)}{q^{\prime }(u)}\right) ^{2}+\frac{2\alpha ^{2}}{3}%
\cdot q^{\prime }(u)^{2}=0. 
\]
We easily check that the solution $q(u)$ with the initial value 
\[
q(0)=q^{\prime \prime }(0)=0\quad \mathrm{and}\quad q^{\prime }(0)=1 
\]
is given by 
\[
q(w)=\frac{1}{2i\alpha }\ln \frac{1+i\alpha w}{1-i\alpha w}=\frac{1}{\alpha }%
\tan ^{-1}\alpha w. 
\]
Then, with this $q(w),$ we easily check as well that 
\[
\sqrt{q^{\prime }(w)}\exp \alpha iq(w)=\frac{1}{1-i\alpha w}. 
\]
This completes the proof.\endproof

\begin{lemma}[Ezhov]
\label{MV}Let $M$ be an analytic real hypersurface in Moser-Vitushkin normal
form and $\phi $ be a normalization of $M$ to Moser-Vitushkin normal form
such that $\phi $ leaves the $u$-curve invariant. Then the mapping $\phi $
is given by 
\[
\phi :\left\{ 
\begin{array}{l}
z^{*}=\sqrt{\mathrm{sign\{}q^{\prime }(0)\mathrm{\}}q^{\prime }(w)}%
Uze^{\alpha i(q(w)-w)} \\ 
w^{*}=q(w)
\end{array}
\right. 
\]
where 
\[
\langle Uz,Uz\rangle =\mathrm{sign\{}q^{\prime }(0)\mathrm{\}}\langle
z,z\rangle 
\]
and the function $q(u)$ is a solution of the equation: 
\[
\frac{q^{\prime \prime \prime }}{3q^{\prime }}-\frac{1}{2}\left( \frac{%
q^{\prime \prime }}{q^{\prime }}\right) ^{2}+\frac{2\alpha ^{2}}{3}\left(
q^{\prime 2}-1\right) =0. 
\]
\end{lemma}

\proof
Suppose that the real hypersurface $M$ is in Moser-Vitushkin normal form is
defined by the equation 
\[
v=\langle z,z\rangle +\sum_{\min (s,t)\geq 2}F_{st}(z,\bar{z},u) 
\]
where 
\begin{eqnarray*}
\Delta F_{22}(z,\overline{z},u) &=&4\alpha (n+1)\langle z,z\rangle \\
\Delta ^{2}F_{23}(z,\overline{z},u) &=&0 \\
\Delta ^{3}F_{33}(z,\overline{z},u) &=&32\alpha ^{2}n(n+1)(n+2).
\end{eqnarray*}
Let $\phi $ be a normalization of $M$ to Moser-Vitushkin normal form leaving
the $u$-curve invariant. Then the mapping $\phi $ is necessarily of the
form(cf. the proof of Theorem \ref{exi-uni}) 
\[
\phi :\left\{ 
\begin{array}{l}
z^{*}=\sqrt{\mathrm{sign\{}q^{\prime }(0)\mathrm{\}}q^{\prime }(w)}U(w)z \\ 
w^{*}=q(w)
\end{array}
\right. 
\]
where 
\[
\langle U(u)z,U(u)z\rangle =\mathrm{sign\{}q^{\prime }(0)\mathrm{\}}\langle
z,z\rangle . 
\]
Suppose that $\phi $ transforms $M$ to a real hypersurface $M^{\prime }$
defined by 
\[
v=\langle z,z\rangle +\sum_{\min (s,t)\geq 2}F_{st}^{*}(z,\bar{z},u). 
\]
Then we obtain 
\begin{align*}
F_{22}(z,\overline{z},u)& =q^{\prime }(u)F_{22}^{*}(U(u)z,\overline{U(u)z}%
,q(u)) \\
& -i\langle z,z\rangle \langle z,U(u)^{-1}\{U^{\prime }(u)+\frac{1}{2}%
q^{\prime }(u)^{-1}q^{\prime \prime }(u)U(u)\}z\rangle \\
& +i\langle z,z\rangle \langle U(u)^{-1}\{U^{\prime }(u)+\frac{1}{2}%
q^{\prime }(u)^{-1}q^{\prime \prime }(u)U(u)\}z,z\rangle \\
F_{23}(z,\overline{z},u)& =q^{\prime }(u)\sqrt{\left| q^{\prime }(u)\right| }%
F_{23}^{*}(U(u)z,\overline{U(u)z},q(u)).
\end{align*}
Then we easily see that 
\[
\Delta ^{2}F_{23}^{*}(z,\overline{z},u)=0. 
\]
We require the following condition 
\begin{eqnarray*}
\Delta F_{22}(z,\overline{z},u) &=&\Delta F_{22}^{*}(z,\overline{z},u) \\
&=&4\alpha (n+1)\langle z,z\rangle
\end{eqnarray*}
so that 
\begin{align}
& 4\alpha (n+1)(q^{\prime }(u)-1)\langle z,z\rangle  \nonumber \\
& +i(n+2)\{\langle U(u)^{-1}U^{\prime }(u)z,z\rangle -\langle
z,U(u)^{-1}U^{\prime }(u)z\rangle \}  \nonumber \\
& +i\langle z,z\rangle \{\mathrm{Tr}(U(u)^{-1}U^{\prime }(u))-\overline{%
\mathrm{Tr}(U(u)^{-1}U^{\prime }(u))}\}  \nonumber \\
& =0.  \tag*{(4.11)}  \label{equ 22}
\end{align}
From the equality 
\[
\langle U(u)z,U(u)z\rangle =\mathrm{sign\{}q^{\prime }(0)\mathrm{\}}\langle
z,z\rangle , 
\]
we have identities 
\begin{eqnarray*}
\langle U(u)^{-1}U^{\prime }(u)z,z\rangle +\langle z,U(u)^{-1}U^{\prime
}(u)z\rangle &=&0 \\
\mathrm{Tr}(U(u)^{-1}U^{\prime }(u))+\overline{\mathrm{Tr}%
(U(u)^{-1}U^{\prime }(u))} &=&0.
\end{eqnarray*}
The equality \ref{equ 22} comes to 
\[
2\alpha i(n+1)(q^{\prime }(u)-1)id_{n\times n}=(n+2)U(u)^{-1}U^{\prime }(u)+%
\mathrm{Tr}(U(u)^{-1}U^{\prime }(u))id_{n\times n}, 
\]
which yields 
\[
\mathrm{Tr}(U(u)^{-1}U^{\prime }(u))=\alpha ni(q^{\prime }(u)-1). 
\]
Hence we obtain 
\[
U^{\prime }(u)=\alpha i(q^{\prime }(u)-1)U(u). 
\]
Thus the function $U(u)$ is given by 
\[
U(u)=U(0)e^{\alpha i(q(u)-u)}. 
\]

Then the mapping $\phi $ is necessarily of the form: 
\[
\phi :\left\{ 
\begin{array}{l}
z^{*}=\sqrt{\mathrm{sign\{}q^{\prime }(0)\mathrm{\}}q^{\prime }(w)}U(0)z\exp
\alpha i(q(w)-w) \\ 
w^{*}=q(w)
\end{array}
\right. 
\]
where 
\[
\langle U(0)z,U(0)z\rangle =\mathrm{sign\{}q^{\prime }(0)\mathrm{\}}\langle
z,z\rangle . 
\]
Then we have 
\begin{align}
F_{33}(z,\overline{z},u)& =q^{\prime }(u)\left| q^{\prime }(u)\right|
F_{33}^{*}(U(0)z,\overline{U(0)z},q(u))  \nonumber \\
& -6\alpha q^{\prime }(u)(q^{\prime }(u)-1)\langle z,z\rangle
F_{22}^{*}(U(0)z,\overline{U(0)z},q(u))  \nonumber \\
& +\left\{ -\frac{q^{\prime \prime \prime }(u)}{3q^{\prime }(u)}+\frac{1}{2}%
\left( \frac{q^{\prime \prime }(u)}{q^{\prime }(u)}\right) ^{2}+3\alpha
^{2}\left( q^{\prime }(u)-1\right) ^{2}\right\} \langle z,z\rangle ^{3}. 
\tag*{(4.12)}  \label{f33}
\end{align}
We have the following identities: 
\begin{eqnarray*}
\Delta ^{3}\langle z,z\rangle ^{3} &=&6n(n+1)(n+2) \\
\Delta ^{3}\left\{ F_{33}^{*}(U(0)z,\overline{U(0)z},q(u))\right\} &=&%
\mathrm{sign\{}q^{\prime }(0)\mathrm{\}}\Delta ^{3}F_{33}^{*}(z,\overline{z}%
,q(u)) \\
\Delta ^{3}\left\{ \langle z,z\rangle F_{22}^{*}(U(0)z,\overline{U(0)z}%
,q(u))\right\} &=&3(n+2)\Delta ^{2}F_{22}^{*}(z,\overline{z},q(u)).
\end{eqnarray*}
Then, requiring in \ref{f33} the following condition 
\begin{eqnarray*}
\Delta ^{3}F_{33}(z,\overline{z},u) &=&\Delta ^{3}F_{33}^{*}(z,\overline{z}%
,u) \\
&=&32\alpha ^{2}n(n+1)(n+2),
\end{eqnarray*}
we obtain 
\begin{eqnarray*}
\frac{q^{\prime \prime \prime }(u)}{3q^{\prime }(u)}-\frac{1}{2}\left( \frac{%
q^{\prime \prime }(u)}{q^{\prime }(u)}\right) ^{2} &=&6\alpha ^{2}\left(
q^{\prime }(u)-1\right) ^{2}+\frac{16\alpha ^{2}}{3}\left( q^{\prime
}(u)^{2}-1\right) \\
&&-12\alpha ^{2}q^{\prime }(u)\left( q^{\prime }(u)-1\right) \\
&=&-\frac{2\alpha ^{2}}{3}\left( q^{\prime }(u)^{2}-1\right) .
\end{eqnarray*}
This completes the proof.\endproof

\begin{lemma}[Vitushkin]
\label{Vit}Let $q(u)$ be an analytic solution of the equation 
\begin{equation}
\frac{q^{\prime \prime \prime }(u)}{3q^{\prime }(u)}-\frac{1}{2}\left( \frac{%
q^{\prime \prime }(u)}{q^{\prime }(u)}\right) ^{2}+\frac{2\alpha ^{2}}{3}%
\left( q^{\prime }(u)^{2}-1\right) =0.  \tag*{(4.13)}  \label{schwarz}
\end{equation}
Then the function $q(u)$ is given by the relation 
\[
e^{2\alpha iq(u)}=e^{i\lambda }\frac{e^{2\alpha iu}+\kappa }{1+\overline{%
\kappa }e^{2\alpha iu}} 
\]
where 
\[
\lambda \in \Bbb{R},\Bbb{\quad }\kappa \in \Bbb{C},\Bbb{\quad }\left| \kappa
\right| \neq 1. 
\]
Further, the function $q(u)$ satisfies the relation 
\[
\left[ \frac{q(u_{2})-q(u_{1})}{\pi \alpha ^{-1}}\right] =\mathrm{sign\{}%
q^{\prime }(0)\mathrm{\}}\left[ \frac{u_{2}-u_{1}}{\pi \alpha ^{-1}}\right]
. 
\]
\end{lemma}

\proof
Let $\frak{L}$ be the mapping in Lemma \ref{linear-linear}. Then the
normalization $\phi =(f^{*},g^{*})$ to Moser-Vitushkin normal form is given
by the relation 
\[
\phi =\frak{L}\circ \varphi \circ \frak{L}^{-1} 
\]
where $\varphi =(f,g)$ is a normalization to Chern-Moser normal form.
Explicitly, we obtain 
\[
\left\{ 
\begin{array}{l}
f^{*}(z,w)=\frac{f\left( z(1-i\tan \alpha w),\alpha ^{-1}\tan \alpha
w\right) }{1-i\alpha g\left( z(1-i\tan \alpha w),\alpha ^{-1}\tan \alpha
w\right) } \\ 
g^{*}(z,w)=\alpha ^{-1}\tan ^{-1}\alpha g\left( z(1-i\tan \alpha w),\alpha
^{-1}\tan \alpha w\right)
\end{array}
\right. . 
\]
Here we take 
\begin{equation}
\varphi :\left\{ 
\begin{array}{c}
z^{*}=\frac{Cz}{1-rw} \\ 
w^{*}=\frac{\rho w}{1-rw}
\end{array}
\right.  \nonumber
\end{equation}
so that 
\begin{equation}
\phi :\left\{ 
\begin{array}{l}
z^{*}=\sqrt{\mathrm{sign\{}Q^{\prime }(0)\mathrm{\}}Q^{\prime }(\exp 2\alpha
iw)}Cz \\ 
w^{*}=\frac{1}{2\alpha i}\ln Q(\exp 2\alpha iw)
\end{array}
\right.  \nonumber
\end{equation}
where 
\[
Q(w)=e^{i\lambda }\frac{w+\kappa }{1+\overline{\kappa }w} 
\]
and 
\begin{gather*}
e^{i\lambda }=\frac{\alpha (1+\rho )+ir}{\alpha (1+\rho )-ir},\quad \kappa =%
\frac{\alpha (1-\rho )-ir}{\alpha (1+\rho )+ir}, \\
\rho =Q^{\prime }(0)\neq 0,\quad \rho ,r\in \Bbb{R}.
\end{gather*}

Then the solution $q(u)$ of the equation \ref{schwarz} is given by 
\[
q(u)=\frac{1}{2\alpha i}\ln Q(\exp 2\alpha iu) 
\]
so that 
\[
e^{2\alpha iq(u)}=e^{i\lambda }\frac{e^{2\alpha iu}+\kappa }{1+\overline{%
\kappa }e^{2\alpha iu}} 
\]
where 
\[
\lambda \in \Bbb{R}\quad \text{and}\quad \left| \kappa \right| \neq 1. 
\]

Finally, note that the mapping 
\[
w^{*}=e^{i\lambda }\frac{w+\kappa }{1+\overline{\kappa }w} 
\]
is an automorphism of the circle $S^{1}.$ Thus we obtain 
\[
\left[ \frac{q(u_{2})-q(u_{1})}{\pi \alpha ^{-1}}\right] =\mathrm{sign\{}%
q^{\prime }(0)\mathrm{\}}\left[ \frac{u_{2}-u_{1}}{\pi \alpha ^{-1}}\right]
. 
\]
This completes the proof.\endproof

\begin{theorem}[Vitushkin]
\label{embedding}Let $M$ be a nondegenerate analytic real hypersurface and $%
\gamma $ be a chain passing through the point $p.$ Suppose that there are an
open neighborhood $U$ of $p$ and a normalizing mapping $\phi $ of $M$ to
Moser-Vitushkin normal form such that $\phi $ translates the point $p$ to
the origin and 
\[
\phi (\gamma \cap U)\subset \left\{ z=v=0\right\} . 
\]
Then the biholomorphic mapping $\phi $ of $M$ is biholomorphically continued
along $\gamma $ such that $\gamma $ is mapped by the mapping $\phi $ into
the $u$-curve in Moser-Vitushkin normal form.
\end{theorem}

\proof
Let $M^{\prime }$ be a real hypersurface in Moser-Vitushkin normal form such
that $M^{\prime }$ is maximally extended along the $u$-curve to the interval 
$(u_{-},u_{+}),$ where 
\[
-\infty \leq u_{-}<0<u_{+}\leq \infty . 
\]
Let $\phi $ be a normalizing mapping of $M$ to $M^{\prime }$ such that $%
\gamma $ is mapped by $\phi $ into the $u$-curve and the point $p$ is mapped
by $\phi $ to the origin. Then we claim that the mapping $\phi $ is
biholomorphically continued along $\gamma $ so that 
\begin{equation}
\phi (\gamma )\subset (u_{-},u_{+}).  \tag*{(4.14)}  \label{included}
\end{equation}
Suppose that the assertion is not true. Then there is a chain-segment $%
\lambda :[0,1]\rightarrow \gamma $ such that $\lambda (0)=p$ and $\phi $ is
analytically continued along all subpath $\lambda [0,\tau ],$ $\tau <1,$ but
not the whole path $\lambda [0,1].$ Let $q=\lambda (1).$ Since $\gamma $ is
a chain and $q$ is an interior point of $\gamma ,$ there are an open
neighborhood $V$ of the point $q$ and a normalizing mapping $h$ of $M$ to
Moser-Vitushkin normal form satisfying 
\begin{eqnarray*}
h(q) &=&0 \\
h(V\cap \gamma ) &\subset &\left\{ z=v=0\right\} .
\end{eqnarray*}
We take a point $x$ on $\lambda [0,1]$ such that 
\[
x\in \lambda [0,1)\cap U. 
\]
Then, by Lemma \ref{MV}, there are an open neighborhood $W$ of the point $x$
and a biholomorphic mapping $k$ satisfying 
\[
\phi =k\circ h\quad \text{on }W\cap V. 
\]
By Lemma \ref{MV} and Lemma \ref{Vit}, the mapping $k$ is biholomorphically
extended to an open neighborhood of the whole $u$-curve. Thus passing to an
open subset of $U$ containing $\lambda [0,1]\cap U,$ if necessary, the
following mapping 
\[
k\circ h\quad \text{on }V 
\]
is an analytic continuation of $\phi $ over the point $\lambda (1).$ Then
necessarily we have 
\[
k\circ h(\lambda [0,1]\cap V)\subset (u_{-},u_{+}). 
\]
This completes the proof.\endproof

\subsection{Extension of a chain}

\begin{lemma}
\label{conti-family}Let $M$ be a nondegenerate analytic real hypersurface
and $\gamma :[0,1]\rightarrow M$ be a continuous curve. Then there exist a
continuous family of real hypersurfaces $M_{\tau },$ $\tau \in [0,1],$ in
normal form and a continuous family of biholomorphic mappings $\phi _{\tau
}, $ $\tau \in [0,1],$ such that $\phi _{\tau }$ translates the point $%
\gamma (\tau )$ to the origin and transforms the germ $M$ at $\gamma (\tau )$
to the germ $M_{\tau }$ at the origin for each $\tau \in [0,1]$ and the
radius of convergence of the mapping $\phi _{\tau }$ at the origin depends
only on $M$ and the point $\gamma (\tau ).$
\end{lemma}

\proof
Without loss of generality, we may assume that the point $\gamma (0)$ is the
origin and the real hypersurface $M$ is defined near the origin by 
\[
v=F\left( z,\overline{z},u\right) ,\quad \left. F\right| _{0}=\left.
F_{z}\right| _{0}=\left. F_{\overline{z}}\right| _{0}=0 
\]
and the curve $\gamma [0,1]$ is given by some continuous functions $p(\tau )$
and $q(\tau )$ via the equation 
\[
\gamma :\left\{ 
\begin{array}{l}
z=p(\tau ) \\ 
w=q(\tau )+iF\left( p(\tau ),\overline{p}(\tau ),q(\tau )\right)
\end{array}
\right. \quad \text{for }\tau \in [0,1] 
\]
where 
\[
\overline{q(\tau )}=q(\tau ). 
\]

Let $\varphi _{\tau }$ be a biholomorphic mapping defined by 
\[
\varphi _{\tau }:\left\{ 
\begin{array}{l}
z^{*}=z-p(\tau ) \\ 
w^{*}=w-q(\tau )-2i\sum_{\alpha =1}^{n}z^{\alpha }\left( \frac{\partial F}{%
\partial z^{\alpha }}\right) \left( p(\tau ),\overline{p}(\tau ),q(\tau
)\right)
\end{array}
\right. . 
\]
Then we obtain the real hypersurfaces $\varphi _{\tau }\left( M\right) ,\tau
\in [0,1],$ defined at the origin by the equation 
\[
v=F^{\tau }\left( z,\overline{z},u\right) ,\quad \left. F^{\tau }\right|
_{0}=\left. F_{z}^{\tau }\right| _{0}=\left. F_{\overline{z}}^{\tau }\right|
_{0}=0 
\]
where 
\begin{eqnarray*}
F^{\tau }\left( z,\overline{z},u\right) &=&F\left( z+p(\tau ),\overline{z}+%
\overline{p}(\tau ),u+q^{*}(\tau )\right) \\
&&-F\left( p(\tau ),\overline{p}(\tau ),q(\tau )\right) \\
&&-\sum_{\alpha =1}^{n}z^{\alpha }\left( \frac{\partial F}{\partial
z^{\alpha }}\right) \left( p(\tau ),\overline{p}(\tau ),q(\tau )\right) \\
&&-\sum_{\beta =1}^{n}\overline{z}^{\beta }\left( \frac{\partial F}{\partial 
\overline{z}^{\beta }}\right) \left( p(\tau ),\overline{p}(\tau ),q(\tau
)\right)
\end{eqnarray*}
and 
\begin{eqnarray*}
q^{*}(\tau ) &=&q(\tau )+i\sum_{\alpha =1}^{n}z^{\alpha }\left( \frac{%
\partial F}{\partial z^{\alpha }}\right) \left( p(\tau ),\overline{p}(\tau
),q(\tau )\right) \\
&&-i\sum_{\beta =1}^{n}\overline{z}^{\beta }\left( \frac{\partial F}{%
\partial \overline{z}^{\beta }}\right) \left( p(\tau ),\overline{p}(\tau
),q(\tau )\right) .
\end{eqnarray*}

Let $\psi _{\tau }=\left( f^{\tau },g^{\tau }\right) $ be a normalization of
the germs $\varphi _{\tau }\left( M\right) ,\tau \in [0,1],$ with identity
initial value such that the real hypersurface $\psi _{\tau }\circ \varphi
_{\tau }\left( M\right) $ is defined by the equation 
\[
v=\langle z,z\rangle +\sum_{k=4}^{\infty }F_{k}^{*\tau }(z,\bar{z},u). 
\]
By Theorem \ref{main}, the functions 
\[
\left\{ 
\begin{array}{l}
\tau \longmapsto \left( f^{\tau },g^{\tau }\right) \\ 
\tau \longmapsto \sum_{k=4}^{\infty }F_{k}^{*\tau }(z,\bar{z},u)
\end{array}
\right. 
\]
are conditnuous. Then the mappings $\phi _{\tau }=\psi _{\tau }\circ \varphi
_{\tau }$ and the real hypersurfaces $M_{\tau }=\psi _{\tau }\circ \varphi
_{\tau }\left( M\right) $ for each $\tau \in [0,1]$ satisfy all the required
conditions. This cmpletes the proof.\endproof

\begin{lemma}
\label{u-estimate}Let $M$ be a nondegenerate analytic real hypersurface and $%
\gamma :[0,1]\rightarrow M$ be a curve such that $\gamma [0,\tau ]$ is a
chain-segment for each $\tau <0.$ Let $U$ be an open set satisfying $\gamma
[0,1)\subset U$ and $\phi $ be a normalization of $M$ on $U$ to
Moser-Vitushkin normal form. Suppose that there is a chain-segment $\lambda $
on $\phi \left( M\right) $ in the $u$-curve satisfying 
\[
\phi (\gamma [0,1))\subset \lambda . 
\]
Then 
\[
\sup_{0\leq \tau <1}\left\{ 
\begin{array}{ll}
\frac{\left\| \left( \left. \frac{\partial f}{\partial z}\right| _{\phi
\circ \gamma (\tau )}\right) \right\| }{\left( \left| \left. \frac{\partial g%
}{\partial w}\right| _{\phi \circ \gamma (\tau )}\right| \right) ^{\frac{1}{2%
}}},\quad & \frac{\left\| \left( \left. \frac{\partial f}{\partial z}\right|
_{\phi \circ \gamma (\tau )}\right) ^{-1}\right\| }{\left( \left| \left. 
\frac{\partial g}{\partial w}\right| _{\phi \circ \gamma (\tau )}\right|
\right) ^{-\frac{1}{2}}}
\end{array}
\right\} <\infty 
\]
where 
\[
\phi ^{-1}=\left( f,g\right) . 
\]
\end{lemma}

\proof
Suppose that the real hypersurface $M$ is defined on an open neighborhood $U$
of the origin by the equation 
\[
v=F\left( z,\overline{z},u\right) ,\quad \left. F\right| _{0}=\left.
dF\right| _{0}=0 
\]
and the curve $\gamma [0,1]\subset M\cap U$ is passing through the origin.
Then there is a biholomorphic mapping 
\[
\phi :\left\{ 
\begin{array}{l}
z=z^{*}+D\left( z^{*},w^{*}\right) \\ 
w=w^{*}+g\left( z^{*},w^{*}\right)
\end{array}
\right. 
\]
where 
\[
D_{z}\left( 0,u\right) =0,\quad \Re g(0,u)=0 
\]
such that the mapping $\phi $ straightens $\gamma [0,1]$ into the $u$-curve
and transforms $M$ to a real hypersurface $\phi \left( M\right) $ defined by 
\[
v=F_{11}^{*}\left( z,\overline{z},u\right) +\sum_{s,t\geq 2}F_{st}^{*}\left(
z,\overline{z},u\right) 
\]
where 
\[
\left( \mathrm{tr}\right) ^{2}F_{23}=0. 
\]
Note that $F_{11}^{*}\left( z,\overline{z},u_{\tau }\right) $ is the Levi
form of $M$ at the point $\gamma \left( \tau \right) \in U\cap M$ for $\tau
\in [0,1],$ where 
\[
\left( 0,u_{\tau }\right) =\phi \circ \gamma \left( \tau \right) . 
\]
Since $\lambda $ is a chain-segment, $F_{11}^{*}\left( z,\overline{z}%
,u_{\tau }\right) $ may be finite for all $\tau \in [0,1].$ Thus we can take
a matrix $E_{1}\left( u\right) $ and a real number $c>0$ such that 
\[
F_{11}^{*}\left( z,\overline{z},u\right) =\langle E_{1}\left( u\right)
z,E_{1}\left( u\right) z\rangle 
\]
and 
\[
\sup_{\tau \in [0,1]}\left\{ \left\| E_{1}\left( u_{\tau }\right) \right\|
,\quad \left\| E_{1}\left( u_{\tau }\right) ^{-1}\right\| \right\} \leq
c<\infty . 
\]

We shall show 
\[
\sup_{0\leq \tau <1}\left\{ \left\| E\left( u_{\tau }\right) \right\| ,\text{
}\left\| E\left( u_{\tau }\right) ^{-1}\right\| \right\} <\infty 
\]
where $\phi ^{-1}=\left( f,g\right) $ and 
\[
E(u)=\frac{\left( \left. \frac{\partial f}{\partial z}\right|
_{z=v=0}\right) }{\sqrt{\left| \left. \frac{\partial g}{\partial w}\right|
_{z=v=0}\right| }}. 
\]
Here the function $E(u)$ satisfies the following ordinary differential
equation(cf. \cite{Pa1}) 
\begin{eqnarray*}
&&F_{11}^{*}\left( E\left( u\right) ^{-1}E^{\prime }\left( u\right) z,%
\overline{z},u\right) \\
&=&-\frac{2i}{n+1}\cdot \mathrm{tr}F_{22}^{*}\left( z,\overline{z},u\right) +%
\frac{1}{2}\left( \frac{\partial F_{11}^{*}}{\partial u}\right) \left( z,%
\overline{z},u\right) \\
&&+\frac{i}{(n+1)(n+2)}\cdot \left( \mathrm{tr}\right) ^{2}F_{22}^{*}\times
F_{11}^{*}\left( z,\overline{z},u\right) .
\end{eqnarray*}
We easily see that there is a real number $e>0$ such that 
\[
\sup_{0\leq \tau <1}\left\| E(u_{\tau })^{-1}E^{\prime }(u_{\tau })\right\|
\leq e<\infty . 
\]
Notice that 
\[
E(u)^{-1}E^{\prime }(u)=-\left( E(u)^{-1}\right) ^{\prime }E(u). 
\]
Because $\lambda $ is a chain-segment, we have 
\[
\int_{\phi \circ \gamma [0,1]}du\leq \int_{\lambda }du<\infty . 
\]
Hence we obtain the following estimates 
\begin{eqnarray*}
\left\| E(u_{\tau })\right\| &\leq &\left\| E(u_{0})\right\| \exp \int_{\phi
\circ \gamma [0,1]}edu<\infty \\
\left\| E(u_{\tau })^{-1}\right\| &\leq &\left\| E(u_{0})^{-1}\right\| \exp
\int_{\phi \circ \gamma [0,1]}edu<\infty
\end{eqnarray*}
where 
\[
\left( 0,u_{0}\right) =\phi \circ \gamma \left( 0\right) . 
\]

So the condition $F_{11}^{*}\left( z,\overline{z},u\right) =\langle E\left(
u\right) z,E\left( u\right) z\rangle $ determines the matrix $E\left(
u\right) $ such that 
\[
E\left( u\right) =U\left( u\right) E_{1}\left( u\right) 
\]
where $U\left( u\right) $ is any matrix satisfying 
\[
\langle U\left( u\right) z,U\left( u\right) z\rangle =\langle z,z\rangle . 
\]
Hence we have the following relation 
\begin{eqnarray*}
c^{-1}\left\| E\left( u_{\tau }\right) \right\| &\leq &\left\| U\left(
u_{\tau }\right) \right\| \leq c\left\| E\left( u_{\tau }\right) \right\| \\
c^{-1}\left\| E\left( u_{\tau }\right) ^{-1}\right\| &\leq &\left\| U\left(
u_{\tau }\right) ^{-1}\right\| \leq c\left\| E\left( u_{\tau }\right)
^{-1}\right\|
\end{eqnarray*}
for all $\tau \in [0,1].$ Therefore, we also have showed 
\[
\sup_{\tau \in [0,1)}\left\{ \left\| U\left( u_{\tau }\right) \right\|
,\quad \left\| U\left( u_{\tau }\right) ^{-1}\right\| \right\} <\infty . 
\]
This completes the proof.\endproof

\begin{lemma}
\label{lowest}Let $M$ be an analytic real hypersurface in normal form
defined by 
\[
v=\langle z,z\rangle +F^{*}\left( z,\overline{z},u\right) 
\]
where 
\[
F^{*}\left( z,\overline{z},u\right) =\sum_{k=4}^{\infty }F_{k}^{*}\left( z,%
\overline{z},u\right) . 
\]
Suppose that $M$ is not a real hyperquadric. Then there is an integer $l\geq
4$ such that 
\begin{align}
F_{k}^{*}\left( z,\overline{z},u\right) =0\quad \text{for all }k\leq l-1 
\nonumber \\
F_{l}^{*}\left( z,\overline{z},u\right) \neq 0  \nonumber  \label{integer 1}
\end{align}
for any value of $U,a,\rho ,r$.
\end{lemma}

In the paper \cite{Pa2}, we have given the proof of Lemma \ref{lowest}.

\begin{theorem}
\label{core2}Let $M,M^{\prime }$ be nonspherical analytic real hypersurfaces
and $\gamma :[0,1]\rightarrow M$ be a curve such that $\gamma [0,\tau ]$ is
a chain-segment for each $\tau <1.$ Let $U$ be an open neighborhood of $%
\gamma [0,1)$ and $\phi $ be a biholomorphic mapping on $U$ such that $\phi $
transforms $M$ to a real hypersurface $M^{\prime }$ satisfying 
\[
\phi (U\cap M)\subset M^{\prime } 
\]
and there is a chain-segment $\lambda :[0,1]\rightarrow M^{\prime }$
satisfying 
\[
\phi (\gamma [0,1))\subset \lambda . 
\]
Suppose that there is a real number $c\geq 1$ such that 
\begin{equation}
\sup_{0\leq \tau \leq 1}\sup_{\left( U_{\tau },0,\rho _{\tau },r_{\tau
}\right) \in H_{\lambda (\tau )}(M^{\prime })}\left\| U_{\tau }\right\| \leq
c<\infty  \tag*{(4.15)}  \label{U-auto}
\end{equation}
where $H_{\lambda (\tau )}(M^{\prime })$ is the local automorphism group of $%
M^{\prime }$ at the point $\lambda (\tau )$ in a normal coordinate with the
chain-segment $\lambda $ on the $u$-curve. Then there exists a chain $\Gamma 
$ on $M$ satisfying 
\[
\gamma [0,1]\subset \Gamma , 
\]
i.e., $\gamma [0,1]$ is a chain-segment.
\end{theorem}

\proof
Without loss of generality, we may assume that the real hypersurface $%
M^{\prime }$ is in Moser-Vitushkin normal form with the chain-segment $%
\lambda $ on the $u$-curve so that $M^{\prime }$ is defined by the equation 
\[
v=\frac{1}{4\alpha }\ln \frac{1}{1-4\alpha \langle z,z\rangle }%
+\sum_{k=4}^{\infty }G_{k}\left( z,\overline{z},u\right) . 
\]
Here we assume $\alpha \neq 0$ and later we shall take a sufficiently small
value for $\alpha .$

There exists a continuous function $\tau \mapsto u_{\tau }$ for $\tau \in
[0,1]$ such that 
\[
(0,u_{\tau })=\phi (\gamma (\tau ))\subset \lambda \quad \text{for }\tau \in
[0,1). 
\]
Since the chain-segment $\lambda $ is compact, there is a real number $u_{1}$
such that 
\[
(0,u_{1})=\lim_{\tau \rightarrow 1}\phi (\gamma (\tau ))\in \lambda . 
\]
Then we obtain a continuous family of analytic real hypersurfaces $M_{\tau
}^{\prime },$ $\tau \in [0,1],$ defined near the origin by 
\[
v=\frac{1}{4\alpha }\ln \frac{1}{1-4\alpha \langle z,z\rangle }+G^{\tau
}\left( z,\overline{z},u\right) 
\]
where, for $\tau \in [0,1],$%
\begin{eqnarray*}
G^{\tau }\left( z,\overline{z},u\right) &=&\sum_{k=4}^{\infty }G_{k}\left( z,%
\overline{z},u+u_{\tau }\right) \\
&=&\sum_{k=4}^{\infty }G_{k}^{\tau }\left( z,\overline{z},u\right) .
\end{eqnarray*}

By Lemma \ref{conti-family}, we obtain a continuous family of analytic real
hypersurfaces $M_{\tau },$ $\tau \in [0,1],$ in normal form and a continuous
family of biholomorphic mappings $\varphi _{\tau }$ for the real
hypersurface $M$ and the curve $\gamma :[0,1]\rightarrow M.$ Then, for each $%
\tau \in [0,1),$ there exist an open neighborhood $V_{\tau }$ of the origin
and a chain $\gamma _{\tau }$ on $M_{\tau }$ passing through the origin such
that 
\[
\varphi _{\tau }^{-1}\left( V_{\tau }\cap \gamma _{\tau }\right) \subset
\gamma [0,1). 
\]
Suppose that $M_{\tau },$ $\tau \in [0,1],$ is defined in normal form by 
\[
v=\langle z,z\rangle +\sum_{k=4}^{\infty }F_{k}^{\tau }\left( z,\overline{z}%
,u\right) . 
\]
By Lemma \ref{lowest}, there is a well-defined integer $m_{\tau },$ $\tau
\in [0,1],$ such that 
\[
\left\{ 
\begin{array}{l}
F_{k}^{\tau }\left( z,\overline{z},u\right) =0\quad \text{for }k\leq m_{\tau
}-1 \\ 
F_{m_{\tau }}^{\tau }\left( z,\overline{z},u\right) \neq 0
\end{array}
\right. 
\]
because $M_{\tau }$ is nonspherical.

Let $\phi _{\tau }$ be a normalization of $M_{\tau }$ for each $\tau \in
[0,1)$ to Moser-Vitushkin normal form such that the initial value $\sigma $
of the mapping $\phi _{\tau }$ is given by 
\[
\sigma =\left( id_{n\times n},a_{\tau },1,0\right) 
\]
where $a_{\tau }$ is determined by the condition 
\[
\phi _{\tau }\left( \gamma _{\tau }\cap M_{\tau }\right) \subset \left\{
z=v=0\right\} 
\]

Suppose that $\phi _{\tau }\left( M_{\tau }\right) ,\tau \in (0,1),$ is
defined near the origin by the equation 
\[
v=\frac{1}{4\alpha }\ln \frac{1}{1-4\alpha \langle z,z\rangle }+F^{*\tau
}\left( z,\overline{z},u;a_{\tau }\right) 
\]
where 
\[
F^{*\tau }\left( z,\overline{z},u;a_{\tau }\right) =\sum_{k=4}^{\infty
}F_{k}^{*\tau }\left( z,\overline{z},u;a_{\tau }\right) . 
\]
Notice that the function $\tau \longmapsto a_{\tau }$ is continuous on $%
[0,1) $ and, by Theorem \ref{main}, the function 
\[
\tau \longmapsto F^{*\tau }\left( z,\overline{z},u;a\right) 
\]
is continuous on $[0,1)$ for a fixed $a\in \Bbb{C}^{n}$.

Note that the two real hypersurfaces $\phi _{\tau }\left( M_{\tau }\right) $
and $M_{\tau }^{\prime }$ are in Moser-Vitushkin normal form and
biholomorphic to each other at the origin for all $\tau \in [0,1)$ by a
biholomorphic mapping leaving the $u$-curve invariant. Thus there is a
mapping 
\[
\psi _{\tau }:\left\{ 
\begin{array}{l}
z^{*}=\sqrt{q_{\tau }^{\prime }(w)}U_{\tau }z\exp \alpha i(q_{\tau }(w)-w)
\\ 
w^{*}=q_{\tau }(w)
\end{array}
\right. 
\]
such that 
\[
\phi _{\tau }\left( M_{\tau }\right) =\psi _{\tau }\left( M_{\tau }^{\prime
}\right) \quad \text{for all }\tau \in [0,1). 
\]
Then the function $q_{\tau }(u)$ is a solution of the ordinary differential
equation 
\[
\frac{q^{\prime \prime \prime }}{3q^{\prime }}-\frac{1}{2}\left( \frac{%
q^{\prime \prime }}{q^{\prime }}\right) ^{2}+\frac{2\alpha ^{2}}{3}\left(
q^{\prime 2}-1\right) =0 
\]
with the initial conditions 
\[
\Re q(0)=0,\quad \Re q^{\prime }(0)=\rho _{\tau }\in \Bbb{R}^{+},\quad \Re
q^{\prime \prime }(0)=2\rho _{\tau }r_{\tau }\in \Bbb{R}. 
\]

Suppose that $\psi _{\tau }\left( M_{\tau }^{\prime }\right) ,\tau \in
[0,1), $ is defined by the equation 
\[
v=\frac{1}{4\alpha }\ln \frac{1}{1-4\alpha \langle z,z\rangle }+G^{*\tau
}\left( z,\overline{z},u\right) 
\]
where 
\[
G^{*\tau }\left( z,\overline{z},u\right) =\sum_{k=4}^{\infty }G_{k}^{*\tau
}\left( z,\overline{z},u;U_{\tau },\rho _{\tau },r_{\tau }\right) . 
\]
Since $\phi _{\tau }\left( M_{\tau }\right) =\psi _{\tau }\left( M_{\tau
}^{\prime }\right) $ for $\tau \in [0,1),$ we obtain 
\begin{equation}
F_{k}^{*\tau }\left( z,\overline{z},u;a_{\tau }\right) =G_{k}^{*\tau }\left(
z,\overline{z},u;U_{\tau },\rho _{\tau },r_{\tau }\right) \quad \text{for }%
k\geq 4.  \tag*{(4.16)}  \label{first}
\end{equation}
We take a sequence $\tau _{j},j\in \Bbb{N},$ such that 
\[
\tau _{j}\in [0,1)\quad \text{and}\quad \tau _{j}\nearrow 1. 
\]
Then there exist a matrix $U_{\tau _{j}}$ and a function $q_{\tau _{j}}(u)$
such that 
\[
\psi _{\tau _{j}}:\left\{ 
\begin{array}{l}
z^{*}=\sqrt{q_{\tau _{j}}^{\prime }(w)}U_{\tau _{j}}z\exp \alpha i(q_{\tau
_{j}}(w)-w) \\ 
w^{*}=q_{\tau _{j}}(w)
\end{array}
\right. . 
\]
Lemma \ref{u-estimate} and the condition \ref{U-auto} allow us to assume 
\[
\sup_{j}\left\| U_{\tau _{j}}\right\| <\infty 
\]
so that, passing to a subsequence, if necessary, there exist a matrix $U$
satisfying 
\[
U=\lim_{j\rightarrow \infty }U_{\tau _{j}}. 
\]
By Lemma \ref{Vit}, all the functions $q_{\tau _{j}}(u)$ satisfy the
following estimate 
\begin{eqnarray*}
\left| q_{\tau _{j}}(u)\right| &=&\left| q_{\tau _{j}}(u)-q_{\tau
_{j}}(0)\right| \\
&\leq &\pi \left| \alpha \right| ^{-1}\left\{ \left[ \frac{\left| q_{\tau
_{j}}(u)-q_{\tau _{j}}(0)\right| }{\pi \left| \alpha \right| ^{-1}}\right]
+1\right\} \\
&\leq &\pi \left| \alpha \right| ^{-1}\left\{ \left[ \frac{\left| u\right| }{%
\pi \left| \alpha \right| ^{-1}}\right] +2\right\} \\
&\leq &\left| u\right| +2\pi \left| \alpha \right| ^{-1}.
\end{eqnarray*}
Because $\lambda $ is a chain-segment on $M^{\prime },$ the functions $%
q_{\tau _{j}}(u)$ are bounded in the range which we have interested in.
Further, notice that 
\[
q_{\tau _{j}}\left( \pi \alpha ^{-1}\right) =\pm \pi \alpha ^{-1}\quad \text{%
for all }j. 
\]
Then, passing to a subsequence, if necessary, Montel theorem and Hurwitz
theorem allow us to have a function $q(u)$ such that 
\[
q(u)=\lim_{j\rightarrow \infty }q_{\tau _{j}}(u) 
\]
and 
\[
q\left( 0\right) =0\quad \text{and}\quad q^{\prime }\left( 0\right) \neq 0. 
\]
Hence, passing to a subsequence, if necessary, there is a real number $e>0$
such that 
\[
\sup_{j}\left\{ \left\| U_{\tau _{j}}\right\| ,\left\| U_{\tau
_{j}}^{-1}\right\| ,\left| \rho _{\tau _{j}}\right| ,\left| \rho _{\tau
_{j}}^{-1}\right| ,\left| r_{\tau _{j}}\right| \right\} \leq e<\infty . 
\]

By the definition of the integer $m_{1},$ we have 
\[
\left\{ 
\begin{array}{l}
\lim_{\tau \rightarrow 1}F_{k}^{\tau }\left( z,\overline{z},u\right)
=0,\quad k=4,\cdots ,m_{1}-1, \\ 
\lim_{\tau \rightarrow 1}F_{m_{1}}^{\tau }\left( z,\overline{z},u\right)
\neq 0.
\end{array}
\right. 
\]
Then the funtion $F_{m_{1}+1}^{*\tau }\left( z,\overline{z},u;a_{\tau
}\right) $ may be decomposed to three parts as follows 
\begin{equation}
F_{m_{1}+1}^{*\tau }\left( z,\overline{z},u;a_{\tau }\right)
=F_{m_{1}+1}^{\tau }\left( z,\overline{z},u\right) +H_{m_{1}+1}^{\tau
}\left( z,\overline{z},u;a_{\tau }\right) +L_{m_{1}+1}^{\tau }\left( z,%
\overline{z},u;a_{\tau }\right)  \tag*{(4.17)}  \label{second}
\end{equation}
where

\begin{enumerate}
\item[(1)]  the function $H_{m_{1}+1}^{\tau }\left( z,\overline{z}%
,u;a\right) $ is determined by the function $F_{m_{1}}^{\tau }\left( z,%
\overline{z},u\right) $,

\item[(2)]  the function $L_{m_{1}+1}^{\tau }\left( z,\overline{z}%
,u;a\right) $ is determined by the functions $F_{k}^{\tau }\left( z,%
\overline{z},u\right) ,$ $k\leq m_{1}-1,$

\item[(3)]  the function $H_{m_{1}+1}^{\tau }\left( z,\overline{z}%
,u;a\right) $ is linear with respect to $a$ and the mapping 
\[
a\longmapsto \lim_{\alpha \rightarrow 0}\lim_{\tau \rightarrow
1}H_{m_{1}+1}^{\tau }\left( z,\overline{z},u;a\right) 
\]
is injective(cf. Lemma \ref{Theo1}), where $\alpha $ is the parameter of
Moser-Vitushkin normal form,

\item[(4)]  the function $L_{m_{1}+1}^{\tau }\left( z,\overline{z}%
,u;a\right) $ depends polynomially on the parameter $a$ and 
\[
\lim_{\alpha \rightarrow 0}\lim_{\tau \rightarrow 1}L_{m_{1}+1}^{\tau
}\left( z,\overline{z},u;a\right) =0\quad \text{for any fixed }a\in \Bbb{C}%
^{n}. 
\]
\end{enumerate}

\noindent Notice that there is a real number $\varepsilon _{1}>0$ such that
the mapping 
\[
a\longmapsto H_{m_{1}+1}^{\tau }\left( z,\overline{z},u;a\right) 
\]
is injective for all $\left| \alpha \right| \leq \varepsilon _{1}$ and all $%
\tau \geq 1-\varepsilon _{1}$. We take a value for the parameter $\alpha $
such that $0<\left| \alpha \right| \leq \varepsilon _{1}.$

Then the equalities \ref{first} and \ref{second} yields 
\begin{eqnarray*}
&&H_{m_{1}+1}^{\tau _{j}}\left( z,\overline{z},u;a_{\tau _{j}}\right)
+L_{m_{1}+1}^{\tau _{j}}\left( z,\overline{z},u;a_{\tau _{j}}\right) \\
&=&G_{m_{1}+1}^{*\tau _{j}}\left( z,\overline{z},u;U_{\tau _{j}},\rho _{\tau
_{j}},r_{\tau _{j}}\right) -F_{m_{1}+1}^{\tau _{j}}\left( z,\overline{z}%
,u\right) .
\end{eqnarray*}
By taking smaller $\varepsilon _{1}>0,$ if necessary, the injectivity of the
mapping $a\longmapsto H_{m_{1}+1}^{\tau }\left( z,\overline{z},u;a\right) $
allows to take an estimate of $a_{\tau _{j}}$ such that there is a real
number $c_{1}>0$ and 
\[
\left| a_{\tau _{j}}\right| \leq c_{1}<\infty \quad \text{for all }\tau _{j},%
\text{ }j\in \Bbb{N}. 
\]
Notice that the function $\tau \rightarrow a_{\tau }$ is continuous. Thus
there exists a real number $c>0$ such that 
\[
\left| a_{\tau }\right| \leq c<\infty \quad \text{for all }\tau \in [0,1). 
\]

Therefore, there exists a sufficiently small real number $\delta >0$
independent of $\tau \in [0,1)$ such that the real hypersurface $M_{\tau }$
and the chain $\gamma _{\tau }\subset M_{\tau }$ extend to 
\begin{eqnarray*}
M_{\tau } &=&\phi _{\tau }^{-1}\left( \psi _{\tau }\left( M_{\tau }^{\prime
}\right) \cap B(0;\delta )\right) \\
\gamma _{\tau } &=&\phi _{\tau }^{-1}\left( \left\{ z=v=0\right\} \cap
B(0;\delta )\right)
\end{eqnarray*}
and the mappings $\varphi _{\tau }^{-1}$ extends biholomorphically on $%
B(0;\delta ).$

Then there is a sufficiently small real numbers $\varepsilon >0$ such that 
\[
\gamma (1)\in \varphi _{1-\varepsilon }^{-1}\left( B(0;\delta )\right) 
\]
so that the following curve $\Gamma $ defined by 
\[
\Gamma =\gamma [0,1)\cup \varphi _{1-\varepsilon }^{-1}\left( \gamma
_{1-\varepsilon }\right) 
\]
is a chain on $M$ such that 
\[
\gamma [0,1]\subset \Gamma . 
\]
This completes the proof.\endproof

Note that the condition \ref{U-auto} is trivially satisfied if the Levi form
on the real hypersurface $M$ is definite.

\begin{theorem}
Let $M$ be a strongly pseudoconvex analytic real hypersurface and $\gamma $
be a chain on $M.$ Let $\Gamma :(0,1)\rightarrow M$ be the maximally
extended connected open analytic curve on $M$ containing the chain $\gamma $
and $\Gamma _{0}$ be a maximal subarc of $\Gamma $ such that $\Gamma _{0}$
contains the chain $\gamma $ and $\Gamma _{0}$ is transversal to the complex
tangent hyperplane of $M$ at each point of $\Gamma _{0}.$ Then $\Gamma _{0}$
is a chain, i.e., for each point $p\in \Gamma _{0},$ there exist an open
neighborhood $U$ of the point $p$ and a biholomorphic mapping $\phi $ on $U$
such that 
\[
\phi (U\cap \Gamma )\subset \left\{ z=v=0\right\} 
\]
and the mapping $\phi $ translates the point $p$ to the origin and
transforms the germ $M$ at the point $p$ to normal form.
\end{theorem}

\proof
Let $M^{\prime }$ be a real hypersurface in Moser-Vitushkin normal form such
that $M^{\prime }$ is maximally extended along the $u$-curve to the interval 
$(u_{-},u_{+}),$ where 
\[
-\infty \leq u_{-}<0<u_{+}\leq \infty . 
\]
Suppose that there is a normalizing mapping $\phi $ of $M$ to $M^{\prime }$
such that $\gamma $ is mapped by $\phi $ into the $u$-curve. Then, by
Theorem \ref{embedding}, the mapping $\phi $ is biholomorphically continued
along $\gamma $ so that 
\[
\phi (\gamma )\subset (u_{-},u_{+}). 
\]
By Theorem \ref{core2}, the chain $\gamma $ can be extended on $M,$ say, to
a chain $\Gamma \subset M,$ whenever an end limit of $\gamma $ exists on $M$
and the corresponding end limit of $\phi (\gamma )$ is an interior point of $%
(u_{-},u_{+}).$ By Theorem \ref{embedding}, the mapping $\phi $ is
biholomorphically continued along $\Gamma $ so that 
\[
\phi (\Gamma )\subset (u_{-},u_{+}). 
\]
Hence there exists a unique chain $\Gamma :(0,1)\rightarrow M$ maximally
extended from the chain $\gamma $ such that 
\[
\lim_{\tau \rightarrow 0}\Gamma (\tau )\notin M\quad \text{or}\quad
\lim_{\tau \rightarrow 0}\phi (\Gamma (\tau ))\in \{u_{-},u_{+}\} 
\]
and 
\[
\lim_{\tau \rightarrow 1}\Gamma (\tau )\notin M\quad \text{or}\quad
\lim_{\tau \rightarrow 1}\phi (\Gamma (\tau ))\in \{u_{-},u_{+}\}. 
\]

Suppose that 
\[
\lim_{\tau \rightarrow 0}\Gamma (\tau )\in M\quad \text{and}\quad \lim_{\tau
\rightarrow 0}\phi (\Gamma (\tau ))\in \{u_{-},u_{+}\}. 
\]
We claim that the analytic curve $\Gamma :(0,1)\rightarrow M$ is not
analytically continued over the limit 
\[
q=\lim_{\tau \rightarrow 0}\Gamma (\tau )\in M 
\]
transversely to the complex tangent hyperplane of $M$ at the point $q\in M.$
Otherwise, there exist an open neighborhood $U$ of the point $q$ and an
analytic curve $\lambda :(-1,1)\rightarrow M$ such that 
\[
U\cap \lambda (0,1)\subset \Gamma (0,1)\quad \text{and}\quad \lambda (0)=q 
\]
and $\lambda $ is transversal to the complex tangent hyperplanes of $M$ at
each point of $\lambda .$ Then there exist a sufficiently small real number $%
\varepsilon >0$ such that the analytic curve 
\[
\lambda (-\varepsilon ,\varepsilon )\cup \Gamma (0,1) 
\]
is an chain as well. Then $M^{\prime }$ is analytically extended along the $%
u $-curve over the point $u_{-}$ or $u_{+}.$ This is a contradiction to the
definition of the point $u_{-}$ and $u_{+}$.

Therefore, the chain $\Gamma $ is the maximally extended connected open
analytic curve containing the chain $\gamma $ which is transversal to the
complex tangent hyperplanes of $M$ at each point of $\Gamma .$ This
completes the proof.\endproof

\section{Analytic continuation of a biholomorphic mapping}

\subsection{On a spherical real hypersurface}

\begin{theorem}[Pinchuk, Chern-Ji]
\label{conti}Let $M$ be a spherical analytic real hypersurfaces with
definite Levi form in a complex manifold, $U$ be a connected neighborhood of
a point $p\in M,$ and $\phi $ be a biholomorphic mapping such that $\phi
(U\cap M)\subset S^{2n+1}.$ Then the mapping $\phi $ continues
holomorphically along any path in $M$ as a locally biholomorphic mapping.
\end{theorem}

\proof
Suppose that the assertion is not true. Then there would exists a path $%
\gamma [0,1]$ such that a biholomorphic mapping $\phi $ at the point $%
p=\gamma (0)$ can be biholomorphically continued along all subpath $\gamma
[0,\tau ]$ with $\tau <1,$ but not along the whole path. We set $q=\gamma
(1).$ Since $M$ is spherical, every point of $M$ is umbilic. By Lemma \ref
{umbilic}, there is an open subset $U_{q}$ of the point $q$ and a
biholomorphic mapping $h_{q}$ on $U_{q}$ such that 
\[
h_{q}(U_{q}\cap M)\subset S^{2n+1} 
\]
and we can take $\tau $ satisfying $\gamma (t)\in U_{q}\cap M$ for all $t\in
[\tau ,1]$. Then there are an open neighborhood $U$ of the point $\gamma
(\tau )$ and a unique automorphisms $\varphi $ of $S^{2n+1}$ such that 
\[
\phi =\varphi \circ h_{q}\quad \text{on }U\cap U_{q}. 
\]
By a classical theorem of Poincar\'{e}, $\varphi $ is biholomorphic on an
open neighborhood of $S^{2n+1}.$ Thus passing to an open subset of $U_{q}$
containing $\gamma [0,1]\cap U_{q},$ if necessary, $\varphi \circ h_{q}$ is
an analytic continuation of $\phi $ on $U_{q}.$ This is a contradiction.
This completes the proof.\endproof

\begin{theorem}[Pinchuk]
\label{pcj}Let $D$ be a bounded strongly pseudoconvex domain in $\Bbb{C}%
^{n+1}$ with simply connected real-analytic boundary $\partial D$. Suppose
that $\partial D$ is a spherical analytic real hypersurface. Then there is a
biholomorphic mapping $\phi $ of $D$ onto $B^{n+1}.$
\end{theorem}

\proof
By Lemma \ref{umbilic}, $\partial D$ is locally biholomorphic to $S^{2n+1}.$
We take a point $p\in \partial D$ and an open neighborhood $U$ of $p$ such
that there is a biholomorphic mapping $\phi $ on $U$ satisfying $\phi (U\cap
\partial D)\subset S^{2n+1}.$ Then, by Theorem \ref{conti}, the mapping $%
\phi $ extends along any path on $\partial D$ as a local biholomorphic
mapping. Since $\partial D$ is simply connected, the monodromy theorem
yields a unique biholomorphic extension $\phi ,$ by keeping the same
notation, on an open neighborhood of $\partial D.$

Note that $\phi :\partial D\rightarrow S^{2n+1}$ is an open mapping because $%
\phi $ is biholomorphic on an open neighborhood of $\partial D.$ Since $%
\partial D$ is compact, the mapping $\phi :\partial D\rightarrow S^{2n+1}$
is a covering map. Further, since $S^{2n+1}$ is simply connected, the
mapping $\phi :\partial D\rightarrow S^{2n+1}$ is a simple covering map so
that there exists a biholomorphic inverse $\phi ^{-1}:S^{2n+1}\rightarrow
\partial D.$ By Hartogs extension theorem, the mappings $\phi ,\phi ^{-1}$
extend to open neighborhoods respectively of $\overline{D}$ and $\overline{%
B^{n+1}}.$ Thus the mapping $\phi $ induces a biholomorphic mapping of $D$
onto $B^{n+1}.$ This completes the proof.\endproof

\begin{theorem}
Let $D$ be simply connected open set in a complex manifold with compact
simply connected real-analytic boundary $\partial D$ and compact closure $%
\overline{D}$. Suppose that $\partial D$ is a spherical analytic real
hypersurface. Then there is a biholomorphic mapping $\phi $ of $D$ onto $%
B^{n+1}.$
\end{theorem}

\proof
By the same argument, there is a biholomorphic mapping $\phi $ on an open
neighborhood of the boundary $\partial D$ such that $\phi :\partial
D\rightarrow S^{2n+1}$ is a simple covering map. Thus there exists a
biholomorphic inverse $\phi ^{-1}:S^{2n+1}\rightarrow \partial D.$ By
Hartogs extension theorem, the mapping $\phi ^{-1}$ extends to the open ball 
$B^{n+1}$ as a local biholomorphic mapping. Since $\overline{B^{n+1}}$ and $%
\overline{D}$ are compact, the mapping $\phi ^{-1}:B^{n+1}\rightarrow D$ is
a covering mapping. Since $D$ is simply connected, $\phi
^{-1}:B^{n+1}\rightarrow D$ is a simple covering map. Thus the mapping $\phi 
$ induces a biholomorphic mapping of $D$ onto $B^{n+1}.$ This completes the
proof.\endproof

Let $Q$ be a real hyperquadric in $\Bbb{CP}^{n+1}$ which is defined in a
homogeneous coordinate 
\[
\left( \eta ,\zeta ^{1},\cdots ,\zeta ^{n},\xi \right) \in \Bbb{C}^{n+2} 
\]
by the equation 
\[
\frac{1}{2i}\left( \xi \overline{\eta }-\eta \overline{\xi }\right) =\langle
\zeta ,\zeta \rangle 
\]
where 
\[
\langle \zeta ,\zeta \rangle \equiv \zeta ^{1}\overline{\zeta ^{1}}+\cdots
+\zeta ^{e}\overline{\zeta ^{e}}-\cdots -\zeta ^{n}\overline{\zeta ^{n}}. 
\]
Then the real hyperquadric $Q$ is given in the inhomogeneous coordinate chart%
$\simeq \Bbb{C}^{n+1}$ by the equation 
\[
\frac{1}{2i}\left( w-\overline{w}\right) =\langle z,z\rangle 
\]
where 
\[
z=\left( \frac{\zeta ^{1}}{\eta },\cdots ,\frac{\zeta ^{n}}{\eta }\right)
,\quad w=\frac{\xi }{\eta }. 
\]

\begin{lemma}[Chern-Moser]
\label{extension}Let $Q$ be a real hyperquadric in $\Bbb{CP}^{n+1}$ and $U$
be an open neighborhood of a point $p\in Q.$ Suppose that there is a
biholomorphic mapping $\phi $ on $U$ such that $\phi (U\cap Q)\subset Q.$
Then the mapping $\phi $ extends to be an automorphism of $Q$ which is
biholomorphic on an open neighborhood of $Q.$
\end{lemma}

\proof
By composing a linear mapping of $\Bbb{CP}^{n+1}$ to $\phi ,$ if
necessarily, we may assume that $\phi $ has a fixed point $q\in Q.$ Further,
passing to an inhomogeneous coordinate chart, $\phi $ is a local
automorphism of the real hyperquadric $v=\langle z,z\rangle $ in $\Bbb{C}%
^{n+1}.$ By Theorem \ref{exi-uni}, the mapping $\phi $ is necessarily to be
a fractional linear mapping as follows: 
\begin{eqnarray*}
z^{*} &=&\frac{C(z-aw)}{1+2i\langle z,a\rangle -w(r+i\langle a,a\rangle )} \\
w^{*} &=&\frac{\rho w}{1+2i\langle z,a\rangle -w(r+i\langle a,a\rangle )}.
\end{eqnarray*}
Thus the mapping $\phi $ extends to be a linear mapping in $\Bbb{CP}^{n+1}.$
This completes the proof.\endproof

\begin{theorem}
Let $M$ be a spherical analytic real hypersurface with nondefinite Levi form
in a complex manifold, $U$ be a connected neighborhood of a point $p\in M,$
and $\phi $ be a biholomorphic mapping on $U$ such that 
\[
\phi (U\cap M)\subset Q 
\]
where $Q$ is a real hyperquadric in $\Bbb{CP}^{n+1}.$ Then the mapping $\phi 
$ continues holomorphically along any path on $M$ as a locally biholomorphic
mapping.
\end{theorem}

\proof
Suppose that the assertion is not true. Then there would exists a path $%
\gamma [0,1]$ with $p=\gamma (0)$ such that a biholomorphic mapping $\phi $
on the neighborhood $U$ of $p$ can be biholomorphically continued along all
subpath $\gamma [0,\tau ]$ with $\tau <1,$ but not along the whole path. We
set $q=\gamma (1).$ Since $M$ is spherical, by Lemma \ref{umbilic}, there is
an open subset $U_{q}$ of the point $q$ and a biholomorphic mapping $h_{q}$
on $U_{q}$ such that 
\[
h_{q}(U_{q}\cap M)\subset Q. 
\]
We take $\tau $ satisfying $\gamma (t)\in U_{q}\cap M$ for all $t\in [\tau
,1]$. Then there are an open neighborhood $U$ of the point $\gamma (\tau )$
and a unique automorphisms $\varphi $ of $Q$ such that 
\[
\phi =\varphi \circ h_{q}\quad \text{on }U\cap U_{q}. 
\]
By Lemma \ref{extension}, $\varphi $ is biholomorphic on an open
neighborhood of $Q.$ Then passing to an open subset of $U_{q}$ containing $%
\gamma [0,1]\cap U_{q},$ if necessary, $\varphi \circ h_{q}$ is an analytic
continuation of $\phi $ on $U_{q}.$ This is a contradiction. This completes
the proof.\endproof

\subsection{On a nonspherical real hypersurface}

\begin{lemma}
\label{continuation}Let $M,M^{\prime }$ be nonspherical analytic real
hypersurfaces and $U$ be an open neighborhood of a point $p\in M$. Suppose
that $M^{\prime }$ is compact and the local automorphism group of $M^{\prime
}$ at each point $q\in M^{\prime }$ is compact. Let $\phi $ be a
biholomorphic mapping of $M$ such that $\phi (U\cap M)\subset M^{\prime }.$
Then $\phi $ is analytically continued along any chain $\gamma $ passing
through the point $p.$
\end{lemma}

\proof
Suppose that the assertion is not true. Then there is a chain-segment $%
\gamma :[0,1]\rightarrow M$ such that $\gamma (0)=p$ and $\phi $ can be
biholomorphically continued along all subpath $\gamma [0,\tau ]$ with $\tau
<1,$ but not along the whole path.

Because $M^{\prime }$ is compact, there exists the limit 
\[
q\equiv \lim_{\tau \rightarrow 1}\phi (\gamma (\tau ))\in M^{\prime }. 
\]
By Lemma \ref{chain-chain}, the subarc $\phi \circ \gamma :[0,\tau
]\rightarrow M^{\prime }$ is a chain-segment for all $\tau <1.$ Then, by
Theorem \ref{core2}, there exists a chain $\Gamma ^{\prime }$ on $M^{\prime
} $ such that 
\[
\lim_{\tau \rightarrow 1}\phi (\gamma [0,\tau ])\subset \Gamma ^{\prime }, 
\]
where the condition \ref{U-auto} in Theroem \ref{core2} is satisfied because
the local automorphism group of $M^{\prime }$ at each point $q\in M^{\prime
} $ is compact.

Without loss of generality, we may assume that $M^{\prime }$ is in
Moser-Vitushkin normal form with the chain $\Gamma ^{\prime }$ in the $u$%
-curve. Since $\gamma :[0,1]\rightarrow M$ is a chain-segment, there is a
chain $\Gamma $ on $M$ such that 
\[
\gamma [0,1]\subset \Gamma . 
\]
Note that $\phi (U\cap \gamma [0,1])\subset \Gamma ^{\prime }.$ Then, by
Theorem \ref{embedding}, the mapping $\phi $ is biholomorphically continued
along the chain $\Gamma .$ Since the point $\gamma (1)$ is an interior point
of $\Gamma ,$ $\phi $ is biholomorphically continued on an open neighborhood
of the point $\gamma (1).$ This is a contradiction. This completes the proof.\endproof

\begin{theorem}[Pinchuk, Ezhov-Kruzhilin-Vitushkin]
\label{nonspherical}Let $M,$ $M^{\prime }$ be nonspherical connected
analytic real hypersurfaces in complex manifolds such that $M^{\prime }$ is
compact and every local automorphism group of $M^{\prime }$ at each point is
compact. Suppose that there exist an open neighborhood $U$ of a point $p$ of 
$M$ and a biholomorphic mapping $\phi $ on $U$ such that $\phi (U\cap
M)\subset M^{\prime }.$ Then the mapping $\phi $ is biholomorphically
continued along any path in $M$.
\end{theorem}

\proof
Suppose that the assertion is not true. Then there is a path $\gamma
:[0,1]\rightarrow M$ such that $\gamma (0)=p$ and the mapping $\phi $ can be
biholomorphically continued along all subpath $\gamma [0,\tau ]$ with $\tau
<1,$ but not along the whole path.

Let $V$ be an open neighborhood of the point $q=\gamma (1).$ Then, by
Theorem \ref{core1}, there are a real number $\delta >0$ and a point $x\in
V\cap M$ such that $B(q;\delta )\subset V$ and, for each given curve $\eta
:[0,1]\rightarrow B(q;\delta )\cap M,$ there is a continuous family of
chain-segments 
\[
\Gamma :[0,1]\times [0,1]\rightarrow V\cap M 
\]
where $\Gamma (s,\cdot ):[0,1]\rightarrow V\cap M$ is a chain-segment of $M$
for each $s\in [0,1]$ satisfying 
\[
\Gamma (s,0)=q\quad \text{and}\quad \Gamma (s,1)=\eta (s)\quad \text{for all 
}s\in [0,1]. 
\]
We take $\tau $ such that $\tau <1$ and $\gamma [\tau ,1]\subset B(q;\delta
)\cap M.$ Then there is a continuous family of chain-segments 
\[
\Gamma :[\tau ,1]\times [0,1]\rightarrow V\cap M 
\]
such that 
\[
\Gamma (s,0)=x\quad \text{and}\quad \Gamma (s,1)=\gamma (s)\quad \text{for
all }s\in [\tau ,1]. 
\]

By Lemma \ref{continuation}, the germ $\phi _{\gamma _{\tau }}$ at the point 
$\gamma (\tau )$ is analytically continued to a germ $\phi _{x}$ at the
point $x$ along the chain-segment $\Gamma (\tau ,\cdot ).$ Then, by Lemma 
\ref{continuation}, the germ $\phi _{x}$ is analytically continued to a germ 
$\phi _{\gamma _{s}}$ at each point $\gamma _{s}\in \gamma [\tau ,1]$ along
the chain-segments $\Gamma (s,\cdot ),$ $s\in [\tau ,1].$

We claim that the germs $\phi _{\gamma _{s}},$ $s\in [\tau ,1],$ are the
analytic continuations of the germ $\phi _{\gamma _{\tau }}$ at the point $%
\gamma (\tau )$ along the subarc $\gamma [\tau ,1].$ Otherwise, there would
exist a number $r,$ $\tau <r\leq 1,$ such that the germs $\phi _{\gamma
_{s}},$ $s\in [\tau ,r),$ are analytic continuations of the germ $\phi
_{\gamma _{\tau }}$ at the point $\gamma (\tau ),$ but the germ $\phi
_{\gamma _{r}}$ is not an analytic continuation of the germ $\phi _{\gamma
_{\tau }}$. By the way, the germ $\phi _{\gamma _{r}}$ is an analytic
continuation of the germ $\phi _{x}.$ Note that the chain-segment $\Gamma
(r,\cdot )$ is compact. Thus there is a number $\varepsilon >0$ such that
each germ $\phi _{\Gamma (r,t)},$ $0\leq t\leq 1,$ at the point $\Gamma
(r,t) $ converges absolutely and uniformly on the open ball $B(\Gamma
(r,t);\tau ). $ Then we can find a number $r_{1},$ $\tau <r_{1}<r,$ such
that 
\[
\Gamma (r_{1},[0,1])\subset \bigcup_{0\leq t\leq 1}B(\Gamma (r,t);\tau ). 
\]
Note that the germ $\phi _{\gamma _{r_{1}}}$ is an analytic continuation of $%
\phi _{x}$ along the chain-segment $\Gamma (r_{1},\cdot )$ and, at the same
time, it is an analytic continuation of $\phi _{\gamma _{r}}$ along the
subarc $\gamma [r_{1},r].$ Then, necessarily, $\phi _{\gamma _{r}}$ is an
analytic continuation of $\phi $ at the point $\gamma (\tau )$ along the
path $\gamma [\tau ,r].$ This contradiction proves our claim.

Therefore, the mapping $\phi $ is analytically continued to an open
neighborhood of the point $q=\gamma (1)$ along the path $\gamma [0,1].$ This
is a contradiction. This completes the proof.\endproof

\begin{theorem}[Pinchuk]
\label{pin}Let $D,D^{\prime }$ be bounded strongly pseudoconvex domains in $%
\Bbb{C}^{n+1}$ with simply connected real-analytic boundaries. Suppose that
there is a connected neighborhood $U$ of a point $p\in \partial D$ and a
biholomorphic mapping $\phi $ on $U$ such that $\phi (U\cap \partial
D)\subset \partial D^{\prime }.$ Then $\phi $ extends to a biholomorphic
mapping of $D$ onto $D^{\prime }.$
\end{theorem}

\proof
Suppose that $\partial D$ is a nonspherical real hypersurface. Then, by
Theorem \ref{remove2}, $\partial D$ is nonspherical as well. By Theorem \ref
{nonspherical}, $\phi $ analytically extends along any path on $\partial D.$
Since $\partial D$ is simply connected, by the monodromy theorem, $\phi $
analytically extend to an open neighborhood of $\partial D$ as a local
biholomorphic mapping. Since $\partial D$ is compact, $\phi :\partial
D\rightarrow \partial D^{\prime }$ is a covering map. Since $\partial
D^{\prime }$ is simply connected, $\phi :\partial D\rightarrow \partial
D^{\prime }$ is a simple covering map so that there is a biholomorphic
inverse $\phi ^{-1}:\partial D^{\prime }\rightarrow \partial D.$ Then, by
Hartogs extension theorem, $\phi ,\phi ^{-1}$ analytically extend to open
neighborhoods respectively of $\overline{D},\overline{D^{\prime }}.$

Suppose that $\partial D$ is a spherical real hypersurface. Then, by Theorem 
\ref{remove2} and Lemma \ref{umbilic}, $\partial D^{\prime }$ is spherical
as well. By Theorem \ref{pcj}, the domains $D,D^{\prime }$ are both
biholomorphic to an open ball $B^{n+1}$ so that $D$ is biholomorphic to $%
D^{\prime }.$ This completes the proof.\endproof

\end{document}